\documentclass[a4paper,11pt]{amsart}

\usepackage{amsmath, amsthm, amssymb}
\usepackage{enumerate}
\usepackage{hyperref}

\usepackage{color}

\usepackage[ansinew]{inputenc}
\usepackage{graphicx}

\usepackage[margin=3.2 cm]{geometry}

\newtheorem{theorem}{Theorem}[section]
\newtheorem*{teoA}{Theorem A}
\newtheorem*{teoB}{Theorem B} 

\newtheorem{cor}[theorem]{Corollary}
\newtheorem{teo}[theorem]{Theorem}
\newtheorem{lemma}[theorem]{Lemma}
\newtheorem{lema}[theorem]{Lemma}

\newtheorem{prop}[theorem]{Proposition}

\newtheorem*{claim}{Claim}

\theoremstyle{definition}
\newtheorem{definition}[theorem]{Definition}
\newtheorem{remark}[theorem]{Remark}

\newtheorem{convention}[theorem]{Convention}

\newcommand{\eps}{\varepsilon}

\newcommand{\wt}[1]{\tilde{#1}}
\newcommand{\R}{\mathbb R}
\newcommand{\RR}{\R}
\newcommand{\ZZ}{\mathbb Z}

\newcommand{\cT}{\mathcal{T}}

\newcommand{\cW}{\mathcal{W}}
\newcommand{\cF}{\mathcal{F}}

\newcommand{\cG}{\mathcal{G}}
\newcommand{\cI}{\mathcal{I}}

\newcommand{\wcF}{\widetilde \cF}

\newcommand{\wcW}{\widetilde \cW}
\newcommand{\wcG}{\widetilde \cG}

\newcommand{\cH}{\mathcal{H}}

\newcommand{\cL}{\mathcal{L}}

\newcommand{\cS}{\mathcal{S}}

\newcommand{\mt}{\widetilde M}

\newcommand{\TT}{\mathbb{T}}

\newcommand{\wwt}{\widetilde T}

\title[Classification of transitive partially hyperbolic diffeomorphisms in dimension 3]
{Partially hyperbolic diffeomorphisms, ergodicity, and transverse foliations in dimension 3}

\author{Sergio R.\ Fenley} 
\address{Florida State University, Tallahassee, FL 32306}
\email{sfenley@fsu.edu}

\author{Rafael Potrie} 
\address{Centro de Matem\'atica, Universidad de la Rep\'ublica, Uruguay and IRL-IFUMI CNRS}
\email{rpotrie@cmat.edu.uy}
\urladdr{http://www.cmat.edu.uy/~rpotrie/}

\thanks{S.F. was partially supported by NSF DMS-2054909. R. P. was partially supported by CSIC. This article is in some sense a culmination of a project started over 10 years ago. We would like to thank the input, help and encouragement of many colleagues and friends that have participated in several ways and which are too many to list here. We single out A. Hammerlindl for several comments including pointing out an oversight in \S 6 of the first version of this paper.  The oversight was solved as presented here, but in further work together with him we were able to substantially improve the results on transverse foliations. For this reason, this paper will remain a permanent preprint, and it will be incorporated into the forthcoming article \cite{FHP}. S. Crovisier and S. Martinchich sent us feedback which was important to improve the presentation. R.P. would like to thank S. Llavayol for her help with some figures.}

\begin{document}

\begin{abstract}
We give a complete topological classification of (chain-)transitive partially hyperbolic diffeomorphisms in 3-manifolds in terms of Anosov flows, completing a program proposed by Pujals. In particular, this also allows to give a full answer to the ergodicity conjecture of Hertz-Hertz-Ures for partially hyperbolic diffeomorphisms in dimension 3. This is achieved by showing a result about pairs of transverse $2$-dimensional foliations in 3-manifolds with Gromov hyperbolic leaves (Theorem B),  which may be of independent interest. This paper will be superseded by the forthcoming article \cite{FHP}, where we obtain a much more general version of Theorem B of this paper. This will allow to remove the assumption of chain transitivity of the partially hyperbolic diffeomorphism and will provide a full classification of partially hyperbolic diffeomorphisms in dimension $3$. For this reason, this paper will remain a permanent preprint. 
\bigskip

\end{abstract}

\maketitle

\section{Introduction}
\subsection{Presentation of results}
The main goals of this paper are to complete the topological classification of transitive partially hyperbolic diffeomorphisms in dimension 3, and to obtain a strong consequence for the ergodicty of such systems. In 2001, E.Pujals proposed a conjecture to address the classification problem which was formalized in \cite{BonattiWilkinson}. The proposal consisted in comparing the topology of those systems to Anosov systems. In most cases the conjecture stated that a transitive partially hyperbolic diffeomorphism is a variable time map of an Anosov flow. The conjecture turned out to be false (see \cite{BGP,BGHP}), but it was later understood that the counterexamples still had some connection with Anosov flows (see \cite[Theorem A]{BFP}).  This allowed for a reformulation of the conjecture which involved introducing the notion of collapsed Anosov flows: see for instance \cite[Question 1]{BFP}. More discussion on the conjecture can be found in \cite{BDV, HassPes, HP-survey, CRHRHU, PotICM,BFP,FP-dt} and references therein (historical developments are summarized in \S~\ref{ss.history}).

To make our main statement more concise, we will specialize to manifolds whose fundamental group has exponential growth. The result below solves \cite[Question 1]{BFP},  as well as the similar questions proposed in e.g. \cite{HP-survey,CRHRHU, BGP, BGHP, PotICM} regarding the classification of transitive partially hyperbolic diffeomorphisms in dimension 3 up to the center direction. 

\begin{teoA} 
Let $f: M \to M$ be a (chain-)transitive partially hyperbolic diffeomorphism in a closed 3-manifold $M$ with fundamental group of exponential growth. Then, $f$ is a collapsed Anosov flow. 
\end{teoA}

We will give precise definitions of the concepts appearing in the statement in \S~\ref{s.preliminary}. For now, let us say that if $f:M \to M$ is a transitive collapsed Anosov flow, then, there is a transitive Anosov flow $\varphi_t:M \to M$ and $\beta:M \to M$ a self orbit equivalence of $\varphi_t$ (i.e. a homeomorphism of $M$ sending orbits to orbits) so that $f$ and $\beta$ are semiconjugated. That is, there is a continuous surjective map $h: M \to M$ so that $f \circ h= h \circ \beta$. Note that in particular, this implies that $M$ admits an Anosov flow (see \cite[Question 1]{BGP}, \cite[Question 5.7]{CRHRHU} or \cite[Questions 7 and 9 ]{PotICM}). Our result provides additional information about the map $h$ and how it can fail injectivity (see Definition \ref{def-scaf}). 

We note here that recent results ( \cite[Appendix]{BowdenMassoni}) imply that every self orbit equivalence of an (oriented) Anosov flow can be realized by a collapsed Anosov flow (see \S~\ref{ss.CAF} for precise definitions). Therefore, this gives a complete classification of (oriented) partially hyperbolic diffeomorphisms modulo Anosov flows and their self-orbit equivalences, which are subject of intense study (see \cite{BarthelmeMann}).  For a classification when $\pi_1(M)$ does not have exponential growth see \cite{HP-survey} and \S~\ref{ss.history}. These manifolds are simpler and do not admit Anosov flows: partially hyperbolic diffeomorphisms in them have a different behavior. 

Theorem A also has some direct dynamical consequences. In particular we can give a full proof of the ergodicity conjecture proposed in \cite{HHP2} by Rodriguez Hertz, Rodriguez Hertz and Ures (see also \cite{CRHRHU}): 

\begin{cor}\label{coro-ergo}
Let $f: M \to M$ be a volume preserving partially hyperbolic diffeomorphism of a closed 3-manifold with non-solvable fundamental group. Then, $f$ is accessible, and if it is $C^{1+}$ then it is ergodic (and in fact a $K$-system). 
\end{cor}

This is an immediate consequence of Theorem A thanks to the main result of \cite{FP-acc}. We explain more about ergodicity and this result in \S \ref{ss.ergodic}. We also refer the reader to the paper cited above, and to \cite{FengUres,FengUres2,CRHRHU} and references therein for more detailed accounts of previous results in the direction of this conjecture. Note that it has been previously shown that ergodicity is an abundant property among partially hyperbolic diffeomorphisms with one dimensional center \cite{HHP1,HHP2} (see also \cite{BurnsWilkinson}), but we emphasize that here we are showing that {\emph{all}} conservative partially hyperbolic diffeomorphisms in 3-manifolds with non-solvable fundamental group are ergodic. 

The proof of Theorem A relies on a modification of a strategy proposed in \cite{FPt1s} which reduces the classification problem to a question about transverse pairs of foliations in 3-manifolds. The key in proving the above results is the following\footnote{In a forthcoming paper \cite{FHP}, that will supersede this one, we are going to show a stronger version of Theorem B, without the assumption of Gromov hyperbolicity of leaves, but rather just assuming non existence of toral leaves. The conclusion will be that either the intersection is leafwise Hausdorff (when lifted to the universal cover) or there is a generalized Reeb surface. In the case where $M$ does not have solvable fundamental group and leaves are Gromov hyperbolic, this implies that the second possibility cannot happen.}:

\begin{teoB}
Let $\cF_1, \cF_2$ be transverse foliations by Gromov hyperbolic leaves, then, one of the following holds for $\cG =\cF_1 \cap \cF_2$ the one dimensional foliation obtained by intersection. Then, one of the following holds:
\begin{itemize}
\item either $\cG$ has a generalized Reeb surface in a leaf of $\cF_1$ or $\cF_2$, 
\item or there is a non-empty closed sublamination $\Lambda_1$ of $\cF_1$ and a sublamination $\Lambda_2$ of $\cF_2$ made of \emph{bubble leaves} of $\cG$,
\item or for every leaf $L \in \wt{\cF_1}$ (or $\wt{\cF_2}$) one has that leaves of $\wt{\cG}$ in $L$ are uniform quasgeodesics in $L$.  
\end{itemize}
Moreover, the second option cannot happen if $M$ is either 
hyperbolic, or Seifert. Finally if $\cF_1$ and $\cF_2$ are minimal, 
then the last option can only happen if $\cF_1$ and $\cF_2$ are both the weak stable foliation of a suspension Anosov flow. 
\end{teoB}

\begin{figure}[htbp]
\begin{center}
\includegraphics[scale=1.0]{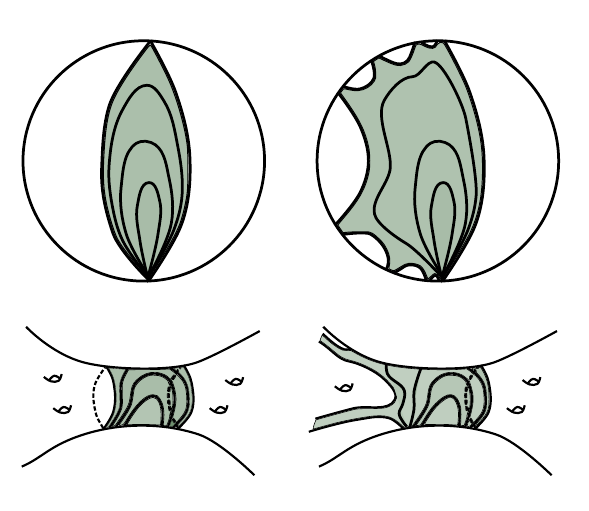}
\begin{picture}(0,0)
\end{picture}
\end{center}
\vspace{-0.5cm}
\caption{{\small In the left a Reeb annulus and its lift to the universal cover of the leaf. In the right, a Reeb crown with its corresponding lift to the universal cover of the leaf. In the orientable case, these are the two instances of generalized Reeb surfaces.}\label{f.reebcrown}}
\end{figure}

 If the third option of Theorem B occurs, we say that $\cG$ is \emph{leafwise quasigeodesic} (see \S~\ref{ss.GHQG}). A \emph{generalized Reeb surface} in $\cF_i$ is a $\cG$-saturated surface $S$ with at least one closed leaf in the boundary and so that boundary components are pairwise non-separated from each other when lifted to the universal cover (see Figure \ref{f.reebcrown} and \S~\ref{ss.Crowns}). The simplest example is a Reeb annulus where each boundary component is a closed curve of $\cG$ and in the interior the leaves are lines
spiralling to the boundary components in the same direction. Examples of transverse foliations for which such Reeb surfaces exist can be found in \cite{MatsumotoTsuboi} (see also \cite[\S 7]{FPt1s}) and in \cite{BBP}. 

Finally, a leaf of $\cF_i$ is a \emph{bubble leaf} if its lift to the universal cover $L$ (notice that $L$ is Gromov hyperbolic and canonically
compactified with a circle $S^1(L)$) satisfies that there is a point $\xi \in S^1(L)$ so that every leaf of $\wt{\cG}|_{L}$ (the lift of $\cG$ restricted
to $L$) verifies that both rays converge to $\xi$ in $L \cup S^1(L)$, see \S~\ref{ss.bubble}. In \cite[\S 6]{FP-hsdff} an example of transverse minimal foliations intersecting in bubble laminations is presented.

We note here that Theorem B is new even if one assumes that both foliations are Anosov foliations, or if one assumes that both foliations are $\RR$-covered\footnote{Recall that a foliation $\cF$ is $\RR$-covered if the leaf space of  the foliation $\wcF$ in the universal cover is homeomorphic to the reals.} (see \cite{BarbotFP} for the case where both assumptions are met). We point out that we obtain a stronger result than Theorem B (see Theorems  \ref{teo-injective} and \ref{teo-noninjective}).  

Up to now some important classes of 3-manifolds have been covered by the classification program of partially hyperbolic diffeomorphisms: those whose fundamental group is solvable \cite{HP-survey},  hyperbolic 3-manifolds \cite{FP-dt}, and unit tangent\footnote{Note however that the proof in unit tangent bundles relies on a strong description of horizontal foliations, unavailable for instance in general Seifert manifolds.} bundles \cite{FPt1s}. However, the only cases where it was possible to treat the case of foliations without additional structure (structure coming from the restrictions imposed by the manifold) was when $f$ is homotopic to the identity (see \cite{BFFP,BFFP2,FP-dt}). In particular those results relied crucially on the fact that a map homotopic to the identity preserving a minimal non-$\RR$-covered foliation has a lift that fixes every leaf in the universal cover and commutes with deck transformations (see \cite[\S 3]{BFFP}). Without assumptions on the isotopy class of $f$, we were not able to get a definite result even in the $\RR$-covered case. In this article we follow a completely different approach, which does not rely on the dynamics except at a small point at the very end of the analysis, to exclude generalized Reeb surfaces.  Note that Anosov flows with non-$\RR$-covered foliations are abundant (see \cite{BoI}) but still their leaf spaces have some well studied structure (see \cite{Fenley-bran}). Here we cannot use any structure of the foliations in the universal cover since it is unknown a priori. To illustrate this fact, we note that a non-trivial\footnote{See \cite[Questions 13 and 14]{PotICM}.} consequence of our results is that chain-transitive partially hyperbolic diffeomorphisms verify that all the leaves of their (branching) foliations have cyclic fundamental group.

\medskip

The link between Theorem A and Theorem B relies on a fundamental result of Burago and Ivanov \cite{BI} who gave the first topological obstructions for the existence of partially hyperbolic diffeomorphisms in 3-manifolds. 
In \cite{BI} it is shown that (under orientabilty conditions) there are $f$-invariant branching foliations which can be blown up to transverse foliations. It is shown in \cite{FPt1s} that the existence of Reeb surfaces in the intersection foliation is incompatible with foliations obtained by partially hyperbolic diffeomorphisms. Note that here we are not able to obtain Reeb surfaces in every case, so we introduce the notion of Reeb crowns that are enough to conclude. It is also shown in \cite{FP-erg} that (chain)-transitive partially hyperbolic diffeomorphisms in manifolds with fundamental group of exponential growth have blown up foliations which are minimal and by Gromov hyperbolic leaves (this is the only place where chain-transitivity is used). Finally in \cite{BFP} it is shown that, under orientability conditions, the leafwise quasigeodesic behavior implies the collapsed Anosov flow property. 

The proof of Theorem B, which comprises the bulk of this paper, uses a result from \cite{FP-hsdff} showing that if the leaf space of $\widetilde{\cG}$ is Hausdorff, then, the foliation $\cG$ is leafwise quasigeodesic in
leaves of both $\cF_1$ and $\cF_2$. Therefore the goal here is to show that non-Hausdorfness leads to the existence of (generalized) Reeb surfaces. 

This study is divided into two parts: the first part assumes that the intersection between any leaf of $\widetilde{\cF_1}$ and any leaf of $\widetilde{\cF_2}$ is connected. This case involves extending and improving results from \cite{FP-hsdff} to produce some geometric properties of the leafwise geometry of the intersected foliation, which is then used to show that non-Hausdorfness allows to produce either bubble sublaminations, or some generalized Reeb surfaces. Part of this analysis relies on importing ideas from \cite{FP-hsdff} (in particular, that one can \emph{push} properties from one leaf to nearby leaves, and extend this to our situation). But it also requires several new ideas, in particular, on how to use the geometric properties of the leaves to show that curves which are non-separated are asymptotic  to closed curves (if there are
no bubble sublaminations) when projected to the manifold. 

The second part of the analysis is the case that there is a pair of leaves $L$ of $\wcF_1$ and $E$ of $\wcF_2$, so that $L \cap E$ is not connected. This part is  more subtle and requires doing an analysis which is  more $3$-dimensional to prove that a configuration similar to the one in the examples from \cite{MatsumotoTsuboi} can be obtained. In particular in this case we produce an actual Reeb surface. Even if this second part requires more delicate arguments, it works with fewer assumptions, and the conclusion about Reeb surfaces is also stronger. The starting point for this study relies on a result  from our recent paper with Barbot \cite[\S 8]{BarbotFP}. 

We believe that Theorem B provides hope in attacking the general problem of obtaining a description of the geometry of transverse foliations in 3-manifolds. In \cite{FHP} we will prove a stronger version of Theorem B which does not assume Gromov hyperbolicity of leaves (just that the foliations are \emph{taut}) and which is enough to prove also a stronger version of Theorem A where chain transitivity is removed. The paper \cite{FHP} will have substantial overlap with the current paper, so this will remain a permanent preprint. 

The orientability conditions are used to construct branching foliations in \cite{BI} and to blow them up into actual foliations. In this article, under the conditions of Theorem A (or even under more general assumptions, see Theorem \ref{teo.mainorient}), we are able to remove the orientability condition: this is done in \S~\ref{s.uniquelyintegrable}. As a byproduct, this allows us to deduce that branching foliations exist even without assuming orientability and to obtain some unique integrability properties of the bundles that also solve several questions raised in \cite{BFP}. It is to be noted that differently from the case of virtually solvable fundamental group, where some quotients are possible (and also classified, see \cite{HP-Solv, HaPo-tori}) here we do not need to consider finite covers to semiconjugate to a model, and in particular, an interesting consequence of our results is that if $f$ is a (chain)-transitive partially hyperbolic diffeomorphism in a closed 3-manifold with non-virtually solvable fundamental group, then, the center bundle is orientable. 

\subsection{Historical remarks and overall strategy} 
\label{ss.history}
Partially hyperbolic systems were introduced in the 70's by Brin-Pesin \cite{BrinPesin} and Hirsch-Pugh-Shub \cite{HiPS} as an extension of Anosov systems. A partially hyperbolic diffeomorphism is one whose tangent map preserves two cone fields, one for the future and one for the past with some general position and expansion/contraction properties. This implies that it preserves stable, unstable, and center bundles, but also that it is a robust property, known to hold for many dynamical systems which present persistent dynamical behavior \cite{BDV,DPU}. See \S~\ref{ss.PH} for precise definitions. 

A general way to obtain a partially hyperbolic diffeomorphism
is the following: start with an Anosov flow $\varphi_t$ and consider
$f = \varphi_1$, the time one map of $\varphi_t$. Then $f$ is partially
hyperbolic with the stable and unstable bundles of $f$ being
the same as the stable and unstable bundles of $\varphi_t$ and the
center bundle of $f$ being generated by the flow direction. If $\tau: M \to \RR_{>0}$ is a smooth function with derivative close to $1$, then $f$ defined by $f(x) := \varphi_{\tau(x)}(x)$ is also partially hyperbolic.
These are examples of what is called a {\em discretized Anosov flow}, see
\cite{BFFP,BFFP2,Martinchich}.

Besides discretized Anosov flows, there are two other obvious families of examples in 3-manifolds. These arise from algebraic constructions: Linear automorphisms
of $\TT^3$ (associated to a matrix $A$ having eigenvalues with distinct
absolute values); and skew products which preserve a foliation by circles and the quotient dynamics is an Anosov map of $\TT^2$ (also viewed as modeled in some fiber preserving automorphisms of nilmanifolds).
In 2001 Pujals asked whether these are all
the possibilities of transitive partially hyperbolic diffeomorphisms in dimension $3$. Probably this was preceded by the observation made in \cite{DPU} that Reeb components of center stable foliations or center unstable foliations were disallowed in some cases. This was later extended and precised in \cite{BBI,BI} where the first topological obstructions for admitting partially hyperbolic diffeomorphisms were devised. 
 
The Pujals question was studied and formalized by Bonatti and Wilkinson in
\cite{BonattiWilkinson}. 
It became known as the Pujals' conjecture and was also reinterpreted in slight variations in many other works (see \cite{BDV,CRHRHU,HassPes}, and \cite[\S 20]{PughShub}). Some positive results were obtained: see \cite{CRHRHU,HP-survey} for surveys. In manifolds with (virtually) solvable fundamental group, the conjecture was proved in \cite{HP-Solv}: every (chain)-transitive partially hyperbolic diffeomorphism in such a manifold is leaf conjugate to an algebraic model as the ones described above\footnote{In such manifolds, when one does not assume chain-transitivity, one needs to take into account the possibility of having $cs$ or $cu$-tori, that is tori tangent to $E^c \oplus E^s$ or $E^c \oplus E^u$, see \cite{HHU3,HaPo-tori}.}.  

The full conjecture was eventually disproved (see \cite{BGP,BGHP}). We describe one of the most simple
examples. Start with the geodesic flow $\varphi_t$ in $M = T^1 S$ where
$S$ is a hyperbolic surface. Suppose that $S$ has a closed geodesic
$\gamma$ which is very short, so that Dehn twist along
$\gamma$ can be effected by a diffeomorphism which distorts
the geometry of $S$ very little. Taking the derivative, this induces
a smooth map $g$ on $M$ which distorts the geometry very little, but
no iterate of it is homotopic to the identity. Now 
consider $t_1, t_2 > 0$ arbitrarily large and let 

$$f \ := \ \varphi_{t_2} \circ g \circ \varphi_{t_1}.$$

\noindent
In \cite{BGP} it was 
proved that $f$ is partially hyperbolic and can be made volume preserving and ergodic. 
No map which is a lift of a finite iterate of $f$ is homotopic to the
identity, so this cannot be a discretized Anosov flow, even
up to iterates or lifts. It cannot be one of the other two types
described in Pujals conjecture either, as the fundamental group of $M$ has exponential growth. In \cite{BGHP} this was  generalized to show that in $M = T^1 S$,
given any mapping class of $S$, there
are partially hyperbolic diffeomorphisms of $M = T^1S$, inducing
that mapping class in $S$. Other examples were also constructed.

However there is much more structure in the counterexamples described above, that was not initially noticed. The construction of those examples required an Anosov flow to start with, but the link with such an Anosov flow was somewhat unclear. This is the genesis of the new structure, called a {\em collapsed
Anosov flow} and introduced in \cite{BFP}. The properties verified by the new examples are varied, so, in \cite{BFP} several possible definitions were proposed (see \S~\ref{ss.CAF} for a more detailed description).  

The goal of this article is to show that the collapsed Anosov flow
property is as general as possible: any partially hyperbolic diffeomorphism which is chain-transitive in a $3$-manifold with $\pi_1(M)$ of exponential growth is a collapsed Anosov flow in the strongest possible sense. 
One lateral contribution of the current work is that it clarifies the many different definitions of collapsed Anosov flows proposed in \cite{BFP} by showing uniqueness of branching foliations (compare with \cite[\S 10]{BFP} and see \S~\ref{s.uniquelyintegrable}).

To accomplish the classification we prove the quasigeodesic property
of center leaves in center stable and center unstable (branching) leaves.
This has an important antecedent: we proved in 
\cite{FP-dt} (relying strongly in  \cite{BFFP,BFFP2}), that this property holds
for partially hyperbolic diffeomorphisms when the underlying
manifold $M$ is hyperbolic. 

Following the introduction of collapsed Anosov flows in 
\cite{BFP},  a strategy was outlined in \cite{FPt1s} to prove the quasigeodesic
property in the more general context of transverse 
two dimensional foliations in $3$-manifolds. Recently, this
strategy was implemented in \cite{FP-hsdff} assuming that the intersected foliation has Hausdorff leaf space in the universal cover.
It is exactly this Hausdorff property that we analyze in this
article. We prove that if the
Hausdorff property fails, then we can construct an object called
a generalized Reeb surface. The generalized Reeb surface
leads to a contradiction if the transverse foliations came
from a partially hyperbolic diffeomorphism. We note that this paper is formally independent from  \cite{BFFP,BFFP2,FP-dt}, though it owes intellectually to those developments. 

\subsection {Consequences on Ergodicity}\label{ss.ergodic} 
Here we discuss Corollary \ref{coro-ergo}.  The conjecture proposed in \cite{HHP2} was motivated by previous work in \cite{HHP1} on the abundance of ergodicity for partially hyperbolic diffeomorphisms with one dimensional center bundle. A conjecture by Pugh and Shub (in turn, motivated by  the work of Grayson-Pugh-Shub \cite{GPS} on stable ergodicity, see \cite{PS-icm,PughShub}) states that among volume preserving partially hyperbolic diffeomorphisms, those which are ergodic form an open and dense subset. 

When the center dimension is one, this was proved in \cite{HHP1} (see also \cite{BurnsWilkinson, Wi} for more discussion on the conjectures in all dimensions). It was noticed that in dimension 3, it could be that non-ergodic partially hyperbolic diffeomorphisms could be described completely. This lead to the conjecture in \cite{HHP2} stating that the only obstruction to ergodicity would be the existence of a torus tangent to the $E^s \oplus E^u$ distribution. Note that the only 3-manifolds which can admit such tori have fundamental group which is solvable (see \cite{CRHRHU}). 

There has been intense work in this conjecture, surveyed in \cite{CRHRHU}, but many advances have also appeared after that survey was published. To describe briefly the status up to Corollary \ref{coro-ergo} which settles completely the conjecture, let us start by saying that in the case where the fundamental group is virtually solvable the conjecture follows from previous work of \cite{HHP2, Hammerlindl, HU, GS} with completely different techniques. Given that, the remaining case of the conjecture is then to show that ergodicity holds unconditionally for conservative partially hyperbolic diffeomorphisms in 3-manifolds with non-virtually solvable fundamental group. 

In \cite{HRHU, FP-erg} the conjecture was settled in some specific 3-manifolds or isotopy classes of diffeomorphisms including hyperbolic 3-manifolds and some isotopy classes in Seifert manifolds. Under the assumption of non-existence of periodic orbits, the conjecture was settled in \cite{FengUres} and more recently, an important progress was made in \cite{FengUres2} assuming that the dynamics is homotopic to the identity but without conditions on the manifold.  We refer the reader to the introduction of that paper for an updated account on the conjecture and more results in the direction of the conjecture. 

In this paper we prove the full conjecture as a consequence of our classification result. Note that volume preserving diffeomorphisms are chain-recurrent, so our Theorem A applies to volume preserving partially hyperbolic diffeomorphisms, and we have established in \cite{FP-acc} that volume preserving collapsed Anosov flows are ergodic (in fact K-systems, thanks to \cite{BurnsWilkinson}) so this gives Corollary \ref{coro-ergo}. 

The strategy in \cite{FP-acc} consists in showing accessibility of such systems\footnote{In fact, it is shown that non-wandering collapsed Anosov flows are accessible.}.
Accessibility means that any two points may be connected by a
path which is a concatenation of paths in stable and unstable leaves, see
\cite{HHP2}. The result in \cite{FP-acc} builds upon \cite{HHP2} where they show that a non-accessible volume preserving partially hyperbolic diffeomorphism must preserve a lamination $\Lambda^{su}$ tangent to $E^s \oplus E^u$, and it is not hard to show (see e.g. \cite{FP-erg}) that $\Lambda^{su}$ can be completed to a foliation transverse to the center direction. 

The study of the interaction
of the lamination $\Lambda^{su}$ with the property of collapsed Anosov flows 
leads to a contradiction \cite{FP-acc}. This shows the accessibility
property for collapsed Anosov flows, and consequently
 the ergodicity property in the conservative setting.
Very roughly the proof in \cite{FP-acc} uses that Anosov flows
in dimension $3$ (if not orbitally equivalent to suspensions) 
admit a pair of periodic orbits which are freely
homotopic to the inverse of each other. 
They project by the collapsing map $h$ to a pair of curves
tangent to the center bundle and which are transverse to the
$su$ lamination $\Lambda^{su}$ and going in opposite directions.
The study of this situation
and the projection by $h$ of the ``lozenge" between the periodic
orbits leads to a contradiction. We refer the reader to
\cite{FP-acc} for more details.

We will not mention ergodicity or accessibility further in this article.

\subsection{Organization of the paper}  In \S~\ref{s.preliminary} we give precise definitions of the objects that are involved in the results and proofs. We also state (and whenever necessary prove) some general some results that are also important for the proofs. In \S~\ref{s.precise} we state the precise versions of our results which imply Theorem's A and B. In \S~\ref{s.teoA} we deduce Theorem A (under orientability conditions) and its stronger version from Theorem B. In \S~\ref{s.teoB} we explain the strategy of the proof of Theorem B and reduce it to two statements that are proved in \S~\ref{s.B1} and \S~\ref{s.B2} respectively after some preliminary considerations are made in \S~\ref{s.lemas}. Section \ref{s.uniquelyintegrable} is somewhat independent on the rest of the arguments, and shows uniqueness of branching foliations under certain conditions, which is then used to show that one can remove the orientability assumptions which we make in the other parts of the paper.


\section{Some preliminary concepts and results}\label{s.preliminary}  
In this paper $M$ will always denote a closed 3-manifold (compact, connected, without boundary). We will denote $\mt$ the universal cover of $M$.

\subsection{Transverse foliations}\label{ss.transverse}
In this article a \emph{foliation} $\cF$ will be a partition of $M$ by immersed surfaces tangent to a continuous distribution $T\cF$ of planes. This is sometimes called $C^{0,1+}$-foliation (see \cite{CandelConlon}). Each surface of the partition is called a \emph{leaf} and we denote by $\cF(x)$ the leaf through the point $x$. We will always assume that $\cF$ has no Reeb components\footnote{Though some results do not need this assumption.}: a Reeb component is a foliation of the solid
torus or Klein bottle so that the boundary is a leaf, and leaves in the interior
are planes \cite{CandelConlon}.
In addition we will always assume that $\cF$ has no sphere or
projective plane leaves.
In this article we will call such 
foliations \emph{Reebless}.
By classical results in foliation theory a Reebless foliation implies the following properties that we will use (see \cite{CandelConlon,Calegari-book}): 

\begin{itemize}
\item the universal cover $\mt$ of $M$ is homeomorphic to $\RR^3$,  
\item each leaf $L \in \widetilde{\cF}$ (the foliation in $\mt$ induced by lifting $\cF$) is a plane separating $\mt$ into two connected components homeomorphic to open balls, in particular, a curve transverse to $\widetilde{\cF}$ cannot be closed, and every leaf $L \in \cF$ is $\pi_1$-injective. 
\end{itemize}

We will say that two (Reebless) foliations $\cF_1, \cF_2$ are \emph{transverse} if the plane fields $T\cF_1$ and $T\cF_2$ are transverse at every point. This implies that the foliations intersect in a one dimensional foliation $\cG = \cF_1 \cap \cF_2$ tangent to the line-field $T\cG = T\cF_1 \cap T\cF_2$. 
We call such $\cG$ a \emph{subfoliation}, in this case of both
$\cF_1$ and $\cF_2$.

When $\cF_1$ and $\cF_2$ are transverse foliations, there is some constant $\eps_0>0$ of \emph{ local product structure} where there are simultaneous trivializing charts which contain an $\eps_0$-neighborhood of every point. 

\begin{remark}
Note that if transverse  foliations have compact leaves, then these can only be tori or Klein bottles since they admit a continuous line-field tangent to them. One reason why we exclude Reeb components is that if they are allowed, many examples can be made that make it difficult to find obstructions. In particular, using Reeb components one can produce transverse foliations in any 3-manifold. In fact, every orientable 3-manifold admits a triple of pairwise transverse foliations (a \emph{total foliation}), see \cite{Hardop}. See also \cite{ST,ADN} for more discussion on total foliations in the presence of Reeb components. We refer also to \cite[\S 7]{Thurston} where this is remarked and the problem of understanding transverse taut foliations\footnote{Taut foliations are those which admit a transversal intersecting every leaf. Being taut (and having no spherical or projective plane leaves) implies the foliation is Reebless. Every foliation without compact leaves is taut.} is presented. 
\end{remark}

\subsection{Leaf spaces} 

The properties of Reebless foliations imply that for each such foliation $\cF$ the leaf space $\cL = \mt/_{\wcF}$ (with the 
quotient topology) is a simply connected, possibly non-Hausdorff, one dimensional manifold. The action of $\pi_1(M)$ on $\cL$ determines the topology of the leaves in many ways, in particular, the fundamental group of a leaf $A \in \cF$ is exactly the stabilizer of a given lift $L$ of $A$ in $\cL$. 

There is a natural dichotomy for a given Reebless foliation $\cF$ of $M$. Either the leaf space is Hausdorff (in which case it is homeomorphic to $\RR$ and we say that $\cF$ is $\RR$-covered), or it is not. In the latter case, importance is given to \emph{non separated leaves}, meaning leaves $L \in \wcF$ (or $\cL$) so that there is a sequence $L_n$ converging both to $L$ and some other leaf in the leaf space
\footnote{We warn the reader that the term {\em branching} is also used to denote the non separated leaves in the literature. In this article we also deal with branching foliations which are associated with partially hyperbolic diffeomorphisms and will be explained later. Because of that we do not use the term branching for leaves of a foliation which are not separated from another leaf.}.

We say that two leaves $L$ and $L'$ are non-separated if there is a sequence $L_n$ converging both to $L$ and $L'$. More generally, given a collection of leaves of $\wcG$ we say that
this collection is \emph{non-separated} if elements of this
collection are pairwise non separated from each other in
the leaf space of $\wcF$. Clearly if $L$ and $L'$ are non-separated, then, both are non-separated leaves. It can be that the set of leaves non-separated from a non-separated leaf $L$ is finite or (countably) infinite.

\begin{figure}[htbp]
\begin{center}
\includegraphics[scale=0.56]{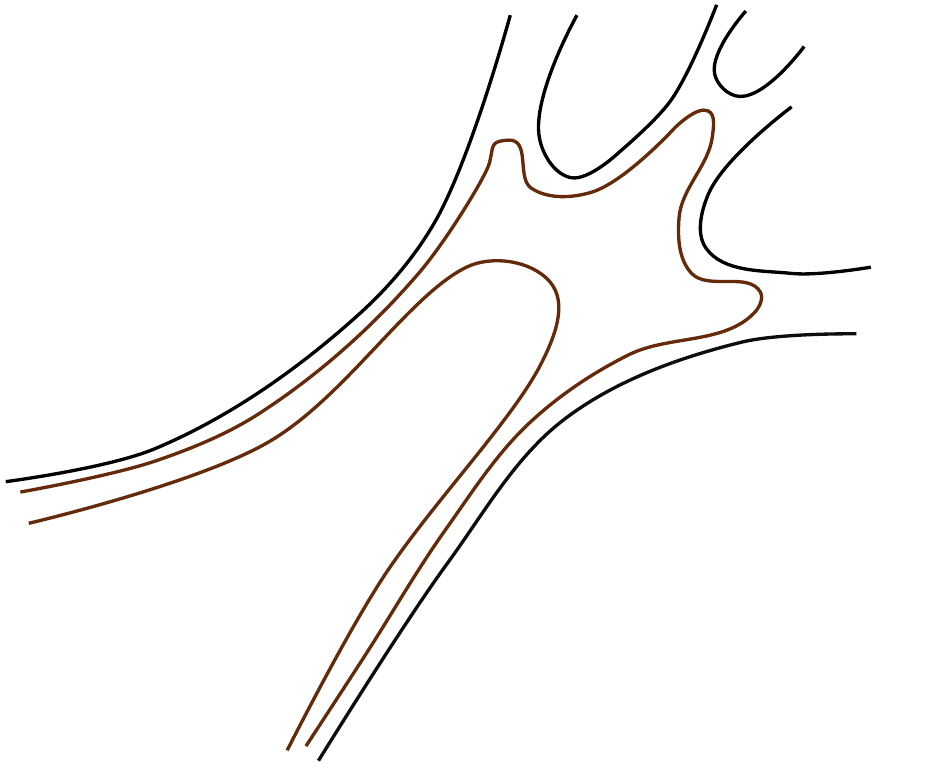}
\begin{picture}(0,0)
\put(-150,26){$r_1$}
\put(-210,100){$r_2$}
\put(-106,190){$\ell_3$}
\put(-55,205){$\ell_4$}
\put(-26,140){$\ell_5$}
\end{picture}
\end{center}
\vspace{-0.5cm}
\caption{{\small Non separated rays $r_1,r_2$ in a one dimensional foliation inside a leaf $L \in \wcF$. The leaves $\ell_3,\ell_4, \ell_5$ are also non separated from the rays $r_1,r_2$. (See definition of $NS(r_1,r_2)$ in \S \ref{ss.returns}.)}\label{f.nsr}}
\end{figure}
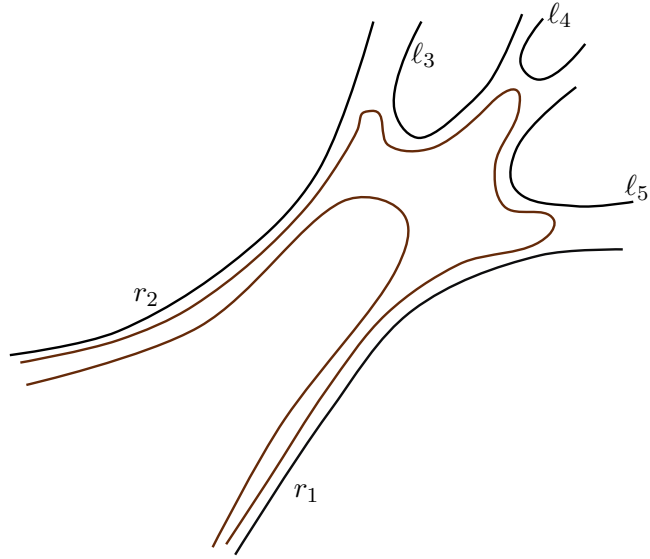

Note that for a one-dimensional subfoliation $\cG$ of a foliation $\cF$ in $M$, we can consider, for each $L \in \wcF$, the foliations $\cG_L = \wcG|_{L}$ and look at the corresponding leaf space. In particular, it makes sense to talk about non-separated leaves of $\cG_L$. For one-dimensional foliations, if $\ell_1, \ell_2$ are two leaves of $\cG_L$ which are non-separated from each other, then it makes sense to say that the non-separation happens in a given \emph{ray} of the leaves and we get non-separated rays $r_1 \subset \ell_1$ and $r_2 \subset \ell_2$. By this, we mean that if $\tau_i$ is a transversal to $r_i$ and if $c_n$ is a sequence of leaves converging to both $\ell_1$ and $\ell_2$, then, the segment of $c_n$ from $\tau_1$ to $\tau_2$ as $n\to \infty$ converges to rays $r_1$ and $r_2$ that we call non-separated. Note that the rays depend on the pair of non-separated leaves (see Figure \ref{f.nsr}). 

\subsection{Reeb surfaces}

Let $\cG$ be a one dimensional subfoliation of a foliation $\cF$.

A \emph{Reeb band} in $L$ leaf of $\wcF$  is a closed region $B \subset L$ saturated by $\cG_L$ whose boundary consists of two leaves $\ell_1, \ell_2 \in \cG_L$ and so that there is a sequence of leaves $\ell_n$ in the interior of $B$ with points $x_n, y_n \in \ell_n$ so that $x_n \to x_\infty \in \ell_1$ and $y_n \to y_\infty \in \ell_2$ (equivalently, the boundary leaves $\ell_1,\ell_2$ are \emph{non-separated} in the leaf space restricted to the band). This is depicted in the upper left side of Figure \ref{f.reebcrown}. 

A \emph{Reeb surface} is a compact connected surface with boundary $R$ contained in a leaf $S$ of $\cF$ so that each lift of $R$ to the universal cover of $S$ is a Reeb band. Here, when $R$ is connected, we say that $\hat R$ is \emph{a lift} of $R$ if it is a connected component of the preimage of $R$ by the universal covering projection. As a consequence, one has that boundary leaves of a Reeb surface $R$ are closed leaves of $\cG$ and that every flow transverse to $\cG$, when entering $R$ cannot escape from $R$. By Euler characteristic reasons a Reeb surface must be either an annulus or a M\"{o}bius band of some leaf of $\cF$.   

\subsection{Reeb crowns and generalized Reeb surfaces}\label{ss.Crowns}
Let $\cG$ be a one-dimensional subfoliation of a foliation $\cF$ of $M$. 

A \emph{Reeb crown} is a connected sub-surface $C$ of a leaf $S$ of $\cF$ which is saturated by $\cG$ and verifies the following properties: 

\begin{itemize}
\item If $L$ is the universal cover of $S$ and $\hat C$ denotes a lift of $C$ to $L$ we have that the boundary leaves $\partial \hat C$ are non-separated in the leaf space of the lifted foliation $\wcG$. 
\item The boundary $\partial C$ of $C$ contains a (unique) closed leaf $\eta$ and at least one other component which is non-compact. 
\item The fundamental group of $C$ is cyclic (in fact, its interior is topologically an open annulus). 
\end{itemize} 

Equivalently, one can define a Reeb crown by saying that it is the projection of a $\wcG$ saturated subset $V$ of a leaf $L \in \wcF$ whose boundary leaves are non-separated, such that the stabilizer of $V$ is cyclic and generated by a deck transformation $\gamma$ and the boundary of $V/_{\langle \gamma \rangle}$ has a unique closed boundary leaf. (See the right side in Figure \ref{f.reebcrown}.)

We define a \emph{generalized Reeb surface} to be either a Reeb surface or a Reeb crown.

We first show a general result about closed curves which lift to non-separated leaves that will be useful for detecting Reeb crowns and Reeb surfaces. 

\begin{lema}\label{lem-closedasympt}
Let $S$ be a leaf of $\cF$ containing a $\cG$ saturated connected set $C$ whose boundary contains a closed curve $c$ and such that when lifted to the universal cover $L$ of $S$ one has that if $\tilde C$ is a lift of $C$ to $L$ (i.e. a connected component of the preimage), then, the boundary leaves of $\tilde C$ are non-separated. Then, if $\tilde c$ is the lift of $c$ in $\tilde C$ it follows that there is $\eps$ so that every leaf of $\wcG$ in $\tilde C$ which has a point at distance less than $\eps$ from $\tilde c$ is asymptotic in one direction to one of the rays of $\tilde c$. 
\end{lema}

\begin{proof}
This follows from the fact that the holonomy in a neighborhood of $\tilde c$ from a point to its image by $\gamma \in \pi_1(S)$ which fixes $\tilde C$ must be expanding or contracting (depending on whether $\gamma$ sends the points in the direction of non-separation or the other) and therefore there is one direction where it is contracting. This is because otherwise holonomy of
$c$ (in $M$) would have fixed points arbitrarily near $c$, which contradicts that the boundary components of $\tilde C$ are non separated from each other.\  This implies that every curve $\ell \in \wcG$ intersecting some transversal of size $\eps$ in $\tilde C$ (note that the $\eps$ is independent on the point because $\tilde c$ projects to a compact curve) must be asymptotic to the side opposite to non-separation to $\tilde c$. 
\end{proof}

The following result  gives a way to detect generalized Reeb surfaces.

\begin{prop}\label{prop-genRS}
Let $C$ be a subsurface of a leaf $S \in \cF$ which is saturated by $\cG$ and such that if $\tilde C$ denotes a lift of $C$ to the universal cover $L$ of $S$, then  the boundary leaves of $\tilde C$ are non-separated. Assume that there is a boundary curve of $\tilde C$ which projects to a closed leaf. Then $C$ is a generalized Reeb surface. Moreover, if more than one boundary curve of $\tilde C$ projects to a closed leaf (which could a priori be the same
curve in $M$), then, $C$ is a Reeb surface. 
\end{prop}
\begin{proof}
We first assume that there are two curves $\ell_1, \ell_2 \in \tilde C$ which project to closed curves in $M$. In this case we will show that $\tilde C$ is a Reeb band, that is, the boundary of $\tilde C$ consists exactly of $\ell_1$ and $\ell_2$. For this, it is enough to show that if $\alpha_i$ is the deck transformation sending $\ell_i$ to itself which is primitive (i.e. it is not a power of a deck transformation that fixes $\ell_i$, and preserves
transversal orientation to $\ell_i$), then, $\alpha_1= \alpha_2^{\pm}$. 

Since $\ell_1$ and $\ell_2$ are non-separated, we can consider a small transversal $\tau: [0,\eps) \to \tilde C$ so that $\tau(0) \in \ell_1$ and so that if $\ell_t$ denotes the leaf through $x_t = \tau(t)$, then $\ell_t$ also converges to $\ell_2$ when $t \to 0$.

The key fact we will use is the consequence of Lemma \ref{lem-closedasympt}: since $\ell_i$ projects to a closed curve, it follows that the curve $\ell_t$ (at least for small $t$) must be asymptotic to $\ell_i$ in one direction. 

Now, assume by contradiction that $\alpha_1$ and $\alpha_2$ do not belong to the same cyclic group (note that since both project to simple closed curves this is the same as assuming that $\alpha_1 \neq \alpha_2^{\pm}$), in particular, we can assume that $\alpha_1(\ell_1)= \ell_1$ and $\alpha_1(\ell_2)\neq \ell_2$. 
It follows that, up to changing $\alpha_1$ for $\alpha_1^{-1}$ and reducing
$\tau$ if necessary, then $\alpha_1$ sends the collection
$\{ \ell_t, 0 <  t < \eps \}$ to $\{ \ell_t,  0 < t < \eps_1 \}$ for some $\eps_1 < \eps$. This shows that this collection also limits to $\ell_2$, in other words
$\alpha_1(\ell_2)$ is non separated from $\ell_2$. 
Notice that both $\ell_2$ and $\alpha_1(\ell_2)$ projects
to closed curves of $\cG$. 
This contradicts that $\ell_t$ has to be asymptotic to $\ell_2$
and also to $\alpha_1(\ell_2)$ in some direction.
Since $\alpha_1(\ell_2) = \ell_2$ then the region 
between $\ell_1$ and $\ell_2$ project to a Reeb surface,
which can be an annulus or a M\"{o}bius band.

To complete the proof we need to show that in case only one boundary component of $\tilde C$ projects into a closed curve, the resulting projection is a Reeb crown surface. For this, it is enough to show that the stabilizer of $\tilde C$ in $\pi_1(S)$ is cyclic. Consider boundary leaves $\ell_1, \ell_2$ of $\tilde C$ so that $\alpha(\ell_1)=\ell_1$ for some primitive $\alpha$. We want to show that the fundamental group of $C$ is generated by $\alpha$. Note that if some $\beta \notin \langle \alpha \rangle$ fixes $\tilde C$ then one has that $\beta(\ell_1) \neq \ell_1$ but one must have that $\beta(\ell_1)$ projects to a closed curve, contradicting the assumption. 
\end{proof}

To show incompatibility with partial hyperbolicity we will use: 

\begin{prop}\label{prop-closedcurve} Let $S$ be a complete surface with a pair of transverse 1-dimensional foliations $\cT_1$ and $\cT_2$ so that the foliation $\cT_1$ contains a generalized Reeb surface. Then, $\cT_2$ contains a closed curve. 
\end{prop}

\begin{proof}
Lift everything to a finite orientable cover. Then, one has that the flow induced by $\cT_2$ can be oriented to be pointing \emph{inward} of the generalized Reeb surface in $\cT_1$. With this property now find a compact annulus $A$ contained in the interior of the generalized Reeb surface so that $\cT_2$ is transverse to the boundary of
$A$ and points inwards to $A$.

As there are no singularities, applying Poincar\'e-Bendixon's theorem,  one deduces that $\cT_2$ must have a periodic orbit giving rise to the posited closed curve. 
\end{proof}


\subsection{Gromov hyperbolicity of leaves and leafwise quasigeodesic subfoliations}\label{ss.GHQG} 
Let $\cF$ be a foliation on $M$. Given a metric on $M$, it induces a metric in each of the leaves of $\cF$. By that we mean a Riemannian metric in $M$
which induces an inner product in tangent spaces of leaves of $\cF$, 
and the path metric in leaves of $\cF$ originating from that. We say that $\cF$ is by \emph{Gromov hyperbolic leaves} if there is $\delta>0$ so that every leaf $L \in \wcF$ is $\delta$-hyperbolic in the sense of Gromov (see \cite{Gr,Gh-Ha}). Note that this is independent of the choice of metric in $M$ (but $\delta$ may change). An important result of Candel (see \cite{Can} or \cite[\S 8]{Calegari-book}) states that for such foliations, there is always a metric in $M$ (smooth along leaves, but only continuous transversally) for which every leaf is of constant negative curvature. See also \cite[Appendix A.3]{BFP} for some remarks in the setting of partially hyperbolic foliations. 

It is known that every foliation without compact leaves in a closed hyperbolic 3-manifold is by Gromov hyperbolic leaves (see \cite{Thurston}). Minimal foliations on manifolds with exponential growth of fundamental group are also by Gromov hyperbolic leaves (see \cite[\S 5]{FP-erg}), and also horizontal foliations 
on Seifert manifolds with base having negative Euler characteristic (see \cite{HP-seif}), just to mention some interesting classes of examples. The obstrution to Gromov hyperbolicity, as stated in Candel's theorem is the existence of a transverse invariant measure with zero Euler characteristic.

When $\cG$ is a one-dimensional subfoliation of $\cF$ 
we now discuss $\cG$ being \emph{leafwise quasigeodesic} (with respect to $\cF$ or in leaves of $\cF$) to mean that when lifted to the universal cover, leaves $\ell \in \wcG$ are quasigeodesic inside the corresponding leaves $L \in \wcF$. We are considering $C^{0,1+}$ foliations, so it makes sense to measure distances with respect to the path metric induced by the
inner product on tangent spaces. Quasigeodesics also make sense and the notion is independent on the choice of metric in $M$ (but the constants may change). For a leaf $\ell \in \wcG$ contained in a leaf $L \in \wcF$ we say that $\ell$ is $k$-quasigeodesic in $L$ if for every $x,y \in \ell$ we have:

$$  d_{\ell} (x,y) \leq k d_L(x,y) + k. $$

\noindent (Note that $d_L(x,y) \leq d_{\ell}(x,y)$ with equality only if $\ell$ is a length minimizing geodesic between $x$ and $y$.) 

We refer the reader to \cite{BFP,FP-hsdff} for more discussion on the notion
of leafwise quasigeodesic foliations. Note in particular that being leafwise quasigeodesic implies that given a metric in $M$ there is a uniform constant on which all leaves are quasigeodesic with this constant, see \cite[Remark 6.2]{BFP}. When $\cF$ has Gromov hyperbolic leaves, then being leafwise quasigeodesic has strong implications, for example the Morse lemma (\cite{Gr,Gh-Ha}) implies that in the universal cover, the leaves of $\wcG$ are at a bounded distance
from length minimizing curves in the corresponding leaves
of $\wcF$.

\subsection{Bubble laminations}\label{ss.bubble} 

If $\cF$ has Gromov hyperbolic leaves and $\cG$ is a one dimensional subfoliation we say that a leaf $L \in \wcF$ is a {\em bubble leaf} if there is
a point $p$ in $S^1(L)$ so that all
leaves of $\wcG$ in $L$ only accumulate on $p$.
The set of leaves which lift to bubble leaves is closed \cite{FP-hsdff} and we say that $\Lambda \subset \cF$ is a \emph{bubble lamination} if every leaf of $\Lambda$ lifts to a bubble leaf.

\subsection{Partially hyperbolic diffeomorphisms}\label{ss.PH}
A diffeomorphism $f:M \to M$ is \emph{partially hyperbolic} if its tangent map $Df$ preserves a splitting into $1$-dimensional bundles $TM= E^s \oplus E^c \oplus E^u$ with the property that there is $N>1$ such that for every $x \in M$ and unit vectors $v^\sigma \in E^\sigma(x)$ (with $\sigma=s,c,u$) one has that: 

$$ 2 \|Df^N v^s \| \leq \min \{ 1, \|Df^N v^c \| \} \leq \max \{ 1, \|Df^N v^c \| \} \leq \frac{1}{2} \|Df^N v^u\|. $$

It is well known that the bundles $E^s$ and $E^u$ integrate uniquely into $f$-invariant foliations, that we will denote by $\cW^s$ and $\cW^u$.  They are called the \emph{strong stable} and \emph{strong unstable} foliations (see \cite{HiPS,CrovPot}). An important property about these foliations is the following:

\begin{remark}\label{rem-noclosed}
The foliations $\cW^s$ and $\cW^u$ do not have closed leaves. In fact since these foliations are contracted either by forward or backward iterations, a closed leaf would imply the existence of a closed leaf in a foliated neighborhood which is impossible. 
\end{remark}

We will say that a diffeomorphism is \emph{chain-transitive} if there is no proper open set $U \subset M$ (i.e. $U\notin \{ M, \emptyset \}$) such that $f(\overline U) \subset U$. There are many equivalences with this definition (see \cite{CrovPot}).  It is implied by two classical assumptions on dynamics:  if $f$ is transitive  (i.e. it has a dense orbit) then it is chain-transitive; also, if $f$ is volume preserving or non-wandering then $f$ is chain-transitive. We will use the following standard result:

\begin{prop}\label{prop-liftcover}
Let $f: X \to X$ be a chain-transitive homeomorphism of a compact connected metric space, then, if $g: \hat X \to \hat X$ is the lift of an iterate of $f$ to compact metric space $\hat X$ (finitely) covering $X$, then, $g$ is chain-transitive. 
\end{prop}
\begin{proof} We use (see \cite{CrovPot}) that a homeomorphism $h:Y \to Y$ is chain-transitive if and only if for every $y \in Y$ and $\eps>0$ there is a \emph{closed} $\eps$-\emph{chain} containing $y$, that is, a finite circularly ordered set $F=\{x_0,x_1, \ldots, x_k\}$
such that for every $0\leq i \leq k$ 
one has $d(h(x_i), x_{i+1})<\eps$ (here by convention $k+1=1$). The closed $\eps$-chain can always be chosen to be minimal in the sense that one cannot take a smaller circularly ordered set contained in $F$ which is a closed $\eps$-chain, and this implies that for every $i$ one has that $h(x_i)$ is only $\eps$-close to $x_{i+1}$. Indeed, if $d(h(x_i), x_j)< \eps$ for some  $j \neq i+1$ then one can eliminate $x_{i+1}, \ldots, x_{j-1}$ to get a smaller circularly ordered set. 

First note that if $f:X \to X$ is chain transitive and $n>0$ we have that $f^n$ is also chain-transitive because given $x \in X$ and $\eps>0$, using the uniform continuity of $f$, one can choose $\delta \ll \eps$ so that if $x_0, \ldots, x_n$ verifies that $d(f(x_i),x_{i+1})\leq \delta$ for all $i$, then one has that $d(f^n(x_0),x_n)< \eps$. This implies that a closed $\delta$-chain for $f$ is a union of at most $n$ closed $\eps$-chains for $f^n$ showing that $f^n$ is also chain-transitive. 

In order to deal with \footnote{We thank Sylvain Crovisier for simplifying our argument to deal with finite covers.} with finite covers we do the following: let $g: \hat X \to \hat X$ be a lift  (let us denote $d>0$ the degree of the cover $\hat X \to X$) of an iterate $f^n$ of $f$ and $\eps > 0$.

Let $\hat x \in \hat X$ and consider $x \in X$ the projection of $\hat x$, then  since $f^n$ (and $f^{-n}$) is chain-transitive, the point $x$ belongs to a closed $\eps$-chain $F$ for $f^n$ (and $f^{-n}$).  We can assume that the $\eps$-chain is minimal so that $f^{n}(z)$ and $f^{-n}(z)$ are $\eps$-close to at most one point of the chain. If $\eps$ is sufficiently small (smaller than half the size of well covered neighborhoods) one can lift 
the closed chain to $\hat X$ and by the above choices
for any $y$ in the lift of the closed chain, there is a unique
$z$ in the lift with $d(g(y),z) < \eps$ and a unique $w$ with
$d(g^{-1}(y),w) < \eps$. Therefore the lift is made up of a finite
number of closed $\eps$-chains, and one of them contains $\hat x$.
\end{proof}

\subsection{Branching foliations}\label{ss.branchingfol}
In general the bundle $E^c$ does not integrate into a foliation and this is one of the important challenges in the classification problem (see \cite{CRHRHU,HHU3,BGHP,PotICM,BFFP3,FP-dt,BFP}). However, since the bundle $E^c$ is one-dimensional and continuous, some structure can be obtained which was produced in the important work of Burago-Ivanov \cite{BI}. 

We will say that a partially hyperbolic diffeomorphism $f$ is \emph{oriented} if the bundles $E^s, E^c, E^u$ are oriented and $Df$ preserves their orientation. 

\begin{teo}(\cite{BI}) \label{teo-BI1}
Any oriented partially hyperbolic diffeomorphism preserves branching foliations $\cW^{cs}$ and $\cW^{cu}$ tangent respectively to $E^{s} \oplus E^c$ and $E^c \oplus E^u$. 
\end{teo}

Branching foliations are collections of immersed surfaces covering $M$ which verify some completeness and separation properties that can be expressed in $\mt$ as follows. A branching foliation $\cW$ of $M$ tangent to a 2-dimensional bundle $E$ is a collection of immersed surfaces which when lifted to the universal cover verify: 

\begin{itemize}
\item each leaf $L \in \wcW$ is everywhere tangent to $E$, is a properly embedded plane and separates $\mt$ in two connected components, 
\item if $L' \neq L$ then $L'$ is disjoint from one of the components of $\mt \setminus L$ (but could intersect $L$), 
\item every $x \in \mt$ belongs to at least one surface of the collection, 
\item if $L_n \in \wcW$ is a sequence of leaves and $x_n \in L_n$ verifies that $x_n  \to x_\infty,$ then, there is a leaf $L_\infty$ through $x_\infty$, such that for every compact disk $D$ in $L_\infty$, up to taking a subsequence $n_k$, there is a sequence of compact disks in $L_{n_k}$ converging in the $C^1$-topology to $D$. 
\end{itemize}

More discussion can be found in \cite[\S 3]{BFP} where we revisit the proof of \cite[Theorem 7.2]{BI} to get some injectivity of the blow up:
\footnote{Andy Hammerlindl \cite{Hammerlindl2} has shown us that it is possible to remove the transverse orientability assumption in Theorem \ref{teo-BI2}. This could give an alternative proof of some of our arguments in
this article.}

\begin{teo}\label{teo-BI2}
Given an branching foliation $\cW$ whose tangent bundle $E$ is orientable and transversally orientable there is a sequence of foliations $\cW_n$ tangent to bundles $E_n$ converging uniformly to $E$ with the property that for every $n>0$ there is a continuous map $h_n: M \to M$ which is $1/n$-$C^0$-close to the identity and is such that: 
\begin{itemize}
\item $h_n$ restricted to each leaf $L \in \cW_n$ is a $C^1$ immersion to a leaf $L' \in \cW$, 
\item for each $L' \in \cW$ there is a unique leaf $L \in \cW_n$ so that $h_n$ maps $L$ to $L'$. 
\end{itemize}
\end{teo}

It is also the case that in many contexts one can ensure that the leaves of the branching foliations are by Gromov hyperbolic leaves, in particular, to prove Theorem A, we need the following consequence of \cite[\S 5]{FP-erg}:

\begin{teo}\label{teo.chaintransitive}
Let $f$ be a chain transitive partially hyperbolic diffeomorphism preserving branching foliations $\cW^{cs}$ and $\cW^{cu}$ in a manifold with fundamental group of exponential growth. Then, for every sufficiently close approximating foliations, these must be minimal and by Gromov hyperbolic leaves. 
\end{teo}
\begin{proof}
In \cite[\S 5]{FP-erg} it is shown that a minimal foliation in a closed manifold with fundamental group of exponential growth is by Gromov hyperbolic leaves. If $\Lambda$ is a proper $f$-invariant closed $\cW^{cu}$-saturated set, then, for some small $\eps>0$ we have that the union of the $\eps$-strong stable manifolds through $\Lambda$ forms a proper open set whose closure is mapped in its interior by some iterate, contradicting chain transitivity (a symmetric argument holds for $\cW^{cs}$, see also \cite[Proposition 5.1]{HPS}). Therefore, the approximating foliations for a chain-transitive partially hyperbolic diffeomorphisms must be minimal (see \cite[Appendix B]{BFFP2} and \cite[Appendix A.3]{BFP}). 
\end{proof}


\subsection{Partial hyperbolicity and generalized Reeb surfaces}

As an application of Remark \ref{rem-noclosed} one can show the following result (see \cite[Proposition 10.4]{FPt1s}): 

\begin{prop}\label{prop.noreebph}
Let $f: M \to M$ be a partially hyperbolic diffeomorphism preserving branching foliations $\cW^{cs}$ and $\cW^{cu}$ and let $\cF_1$ and $\cF_2$ be two closeby approximating foliations. Then, the intersection foliation $\cG = \cF_1 \cap \cF_2$ does not have generalized Reeb surfaces. 
\end{prop}

\begin{proof}
The case of Reeb surfaces is explained in \cite[Proposition 10.4]{FPt1s} and the key idea is the same (using Proposition \ref{prop-closedcurve} instead). Since here we are dealing with a more general case, we will give a proof. Note that there is no loss of generality in assuming that everything is oriented as lifts of partially hyperbolic diffeomorphisms to finite covers are still partially hyperbolic and having a generalized Reeb surface is stable under finite lifts. 

Assume then that the blown up branching foliations give rise to transverse foliations $\cF_1$ and $\cF_2$ whose intersected foliation $\cG$ contains a generalized Reeb surface, say in $\cF_1$. Note that the blowing up given by Theorem \ref{teo-BI2} allows one to pull back the strong stable and strong unstable foliations to leaves of $\cF_1$ and $\cF_2$ respectively since the maps $h_n$ are $C^1$-immersions when restricted to leaves. This gives rise to a foliation in $\cF_1$ which is transverse to $\cG$ and has no closed leaves contradicting Proposition \ref{prop-closedcurve} and completing the proof. 
\end{proof}

\subsection{Collapsed Anosov flows}\label{ss.CAF}
 
Recall that a topological Anosov flow $\varphi_t: M \to M$ is an expansive flow with $C^1$ orbits, and preserving two topological foliations which are topologically transverse (see \cite[\S 5]{BFP} or \cite{BarthelmeMann} for other definitions and context). When $\varphi_t$ is transitive (i.e. it has a dense orbit), it has been proved in \cite{Shannon} that such a flow is orbit equivalent to an Anosov flow (i.e. a $C^\infty$ flow whose time one map is partially hyperbolic). This is still unknown in the non-transitive case.

Here, by orbit equivalence between flows $\varphi_t, \psi_t:M \to M$ we mean a homeomorphism $h: M \to M$ which maps orbits of $\varphi_t$ to orbits of $\psi_t$, and preserves orientation along orbits. A \emph{self orbit equivalence} of a (topological) Anosov flow $\varphi_t$ is an orbit equivalence of $\varphi_t$ with itself. We say that the self orbit equivalence is trivial if it is of the form $\beta(x)  = \varphi_{\tau(x)}(x)$ for some $\tau: M \to \RR$. Two self orbit equivalences are \emph{equivalent} if they differ by a trivial self orbit equivalence.

\begin{definition} \label{def-caf}
We say that a partially hyperbolic diffeomorphism $f: M \to M$ is a \emph{collapsed Anosov flow} if there is a (topological) Anosov flow $\varphi_t: M \to M$, a self orbit equivalence $\beta: M \to M$ of $\varphi_t$ and a continuous and surjective map $h: M \to M$ such that: 

\begin{itemize}
\item one has that $f \circ h = h \circ \beta$, and, 
\item $h$ is $C^1$ along orbits of $\varphi_t$ and $\frac{\partial}{\partial t}|_{t=0} h(\varphi_t(x)) \in E^c(h(x)) \setminus \{0\}$ for every $x \in M$. 
\end{itemize}
\end{definition}

It is shown in \cite[Theorem A]{BFP} that the set of collapsed Anosov flows associated to a given Anosov flow and an equivalence class of self orbit equivalences forms an open and closed set of partially hyperbolic diffeomorphisms. 

A stronger version of collapsed Anosov flows \cite{BFP} is the following:

\begin{definition} \label{def-scaf}
A \emph{strong collapsed Anosov flow} is a collapsed Anosov flow $f$ satisfying the
following additional condition:  it preserves branching foliations $\cW^{cs}$ and $\cW^{cu}$, so that the map $h$ in the definition of collapsed Anosov flow additionally verifies that: 
\begin{itemize}
\item  $h$ maps weak stable and weak unstable leaves of $\varphi_t$ into $\cW^{cs}$ and $\cW^{cu}$-leaves. 
\end{itemize}
\end{definition}

Recently, extending \cite{BGP,BGHP} it was shown that under orientability conditions, every equivalence class of self orbit equivalence of an Anosov flow can be realized by a (strong) collapsed Anosov flow, see  \cite[Appendix]{BowdenMassoni}.

There are two other forms of collapsed Anosov flow that will be used in this paper and which are interrelated. For a partially hyperbolic diffeomorphism $f$ preserving branching foliations $\cW^{cs}$ and $\cW^{cu}$ we call center leaf to any curve which is the projection to $M$ of a connected component of the intersection of a leaf $L \in \wcW^{cs}$ with a leaf $E \in \wcW^{cu}$. We denote by $\cW^c$ to the center (branching) foliation consisting of all center leaves. It is possible to give a topology to the leaf space of $\wcW^{c}$ (i.e. the space of center leaves in $\mt$) which makes it into a two dimensional, possibly non-Hausdorff, simply connected manifold. See \cite[\S 3]{BFP} for more details on this. 

\begin{definition} \label{def-qgph}
A partially hyperbolic diffeomorphism $f$ is a \emph{quasigeodesic partially hyperbolic
diffeomorphism} if $f$ preserves branching foliations
$\cW^{cs}, \cW^{cu}$ with Gromov hyperbolic leaves,
and the center leaves of $\wcW^c$ are uniform quasigeodesics in
leaves of $\wcW^{cs}$ and $\wcW^{cu}$.
\end{definition}

The vast majority of the work in this article is to
prove the quasigeodesic property under certain conditions. The following is \cite[Theorem D]{BFP}
  
 \begin{teo}\label{teo-QGimpliesSCAF}
 Let $f: M \to M$ be an oriented partially hyperbolic diffeomorphism preserving branching foliations $\cW^{cs}$ and $\cW^{cu}$ with Gromov
hyperbolic leaves. Then $f$ is a quasigeodesic partially hyperbolic
diffeomorphism if and only if $f$ 
a strong collapsed Anosov flow. 
  \end{teo}

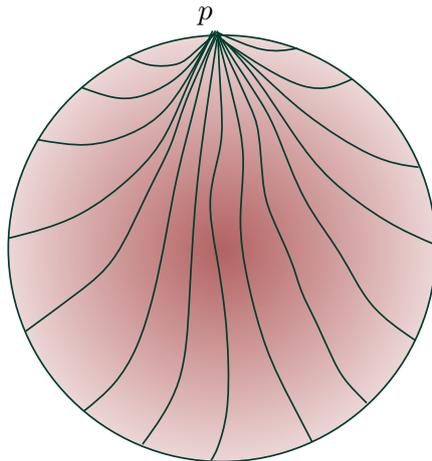
\begin{figure}[htbp]
\begin{center}
\includegraphics[scale=0.56]{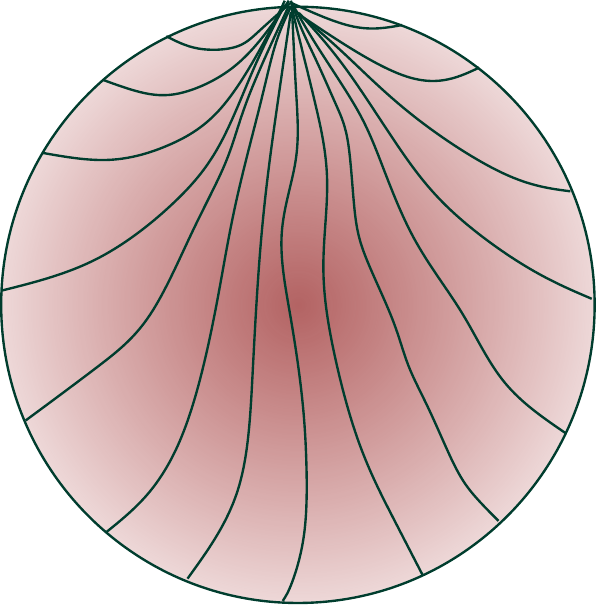}
\begin{picture}(0,0)
\put(-93,167){$p$}
\end{picture}
\end{center}
\vspace{-0.5cm}
\caption{{\small When every curve of $\wcG$ is quasigeodesic in its corresponding leaf, and the foliation comes from a partially hyperbolic diffeomorphism, one can show that in every leaf one sees a \emph{quasigeodesic fan} foliation as depicted in the figure (which is also the figure one sees in leaves of topological Anosov flows). Every quasigeodesic line shares one of the endpoints, depicted by $p$.}\label{f.qgfan}}
\end{figure}

We also need the following:

\begin{definition} \label{def-lscaf}
A partially hyperbolic diffeomorphism $f$ is a \emph{leaf space
collapsed Anosov flow} if it preserves branching foliations
$\cW^{cs}, \cW^{cu}$ and there is a topological
Anosov flow $\varphi_t$, and a $\pi_1(M)$ equivariant
homeomorphim between the orbit space of $\tilde \varphi_t$ and the
center leaf space of $\wcW^c$. 
\end{definition}

In particular, this implies that the leaf space of $\wcW^c$ is homeomorphic to the plane (as is always the case for Anosov flows). See \cite[\S 2.5]{BFP} for more details. 

Upgrading this leaf space condition to the strong collapsed Anosov flow property (which is needed in Section \ref{s.uniquelyintegrable}) necessitated some orientability conditions in \cite{BFP} that we can remove here thanks to the construction of the mentioned topological Anosov flow in $M$. 

\subsection{A result about topological Anosov flows}

In this subsection we show that in the definition of topological Anosov flows, the requirement that orbits
are $C^1$ is not essential and can be eliminated
up to topological equivalence.  We state this separately as it is useful as a stand alone
result.

We say that a one dimensional foliation $\cG$ in
a closed $3$-manifold $M$ is expansive if leaves which
are always very close are the same: specifically
there is $\eps > 0$, so that 
if $\wcG$ is the lifted foliation to the
universal cover $\mt$, then if two
leaves $\ell_1, \ell_2$ have Hausdorff distance 
$< \eps$ from each other, then $\ell_1 = \ell_2$.

\begin{prop} \label{lem-weaktaf}
Suppose that $\cG$ is an expansive one dimensional foliation
in a closed $3$-manifold $M$, 
which is contained in a pair of two dimensional $C^0$-foliations\footnote{In particular, we do not ask leaves to be smooth.} $\cF_1$ and $\cF_2$ which are topologically transverse. Then $\cG$ is topologically equivalent to the flow
foliation of a topological Anosov flow $\varphi_t$.
\end{prop}

\begin{proof}
Assume first that $\cG$ is orientable. Since $\cG$ is expansive, it follows from work of 
 \cite{Pa,IM}, that $\cG$ is topologically equivalent to the
flow foliation of a pseudo-Anosov flow $\varphi_t$. In fact, as in \cite[Theorem 5.9]{BFP} one can show that $\cF_1$ and $\cF_2$ correspond to the weak stable and weak unstable foliations. So, it is enough to show that $\varphi_t$ is orbit equivalent to a flow with $C^1$-orbits.

In \cite{Cal1} it is shown that up to changing the smooth structure on $M$ we can assume that leaves of the weak stable foliation (which we call $\cF_1$) are smooth.
This foliation does not admit holonomy invariant
transverse measures, so it follows from Candel's theorem \cite{Can} that there is a metric on $M$ so that every leaf of $\wcF_1$ is isometric to the hyperbolic plane. 

Let $\cF_2$ the weak unstable foliation of $\varphi_t$ which is (topologically) transverse to $\cF_1$.
The orbit space of $\varphi_t$ is 
is Hausdorff because it is pseudo-Anosov (see \cite[\S 1.3]{BarthelmeMann} or \cite{FenleyMosher}).

In \cite{FP-hsdff} it is shown that under these assumptions, the orbits of $\widetilde \varphi_t$ are quasigeodesics in the leaves of $\wcF_1$ (see in particular \cite[Remark 5.6]{FP-hsdff} which explains the regularity assumptions). 

Given a leaf $L$ of $\wcF_1$, it follows that all flow lines of $\widetilde \varphi_t$
inside $L$ share a common ideal point denoted by
$p$ because $\cF_1$ is the weak stable foliation of $\varphi_t$.
Finally no two flow lines in $L$ share both ideal points in $S^1(L)$:
otherwise the quasigeodesic property implies that 
these two flow lines would be a bounded distance from
each other in $L$ and would contradict expansivity. It follows that for every point $q$ in $S^1(L) \setminus \{ p \}$ there is a unique flow line with ideal points $p, q$.

Then as in \cite[\S 5.4]{BFP} we produce a new flow $\phi_t$
in $M$ so that flow lines are contained in leaves
of $\cF_1$, and in the universal cover,
in a leaf $L$ as above we make the orbits of $\phi_t$ to be the geodesic representatives of each flow line of $\varphi_t$
in $L$. In \cite[Proposition 5.22]{BFP} it is
proved that $\varphi_t$ and $\phi_t$ are orbitally  equivalent (and thus, also orbit equivalent to $\varphi_t$),
and obviously $\phi_t$ has flow lines which are $C^1$ completing the proof. 

If $\cG$ is not orientable, a double
cover of $\cG$ is orientable, and the above applies. But 
the direction where orbits converge in a leaf of $\cF_1$
is invariant, which shows that the original $\cG$ is orientable. This finishes the proof.
\end{proof}


\section{Precise statement of results}\label{s.precise}

The main technical result of this paper is the following result which is a slightly weaker version of Theorem B:

\begin{teo}\label{teo.main}
Let $\cF_1$ and $\cF_2$ be two transverse foliations by Gromov hyperbolic leaves in a closed 3-manifold $M$ with non-virtually solvable fundamental group. Assume that $M$ is hyperbolic or Seifert, or that the foliations are minimal and let $\cG= \cF_1 \cap \cF_2$. Then, either there is a generalized Reeb surface in  one of the foliations or $\cG$ is leafwise quasigeodesic in both $\cF_1$ and $\cF_2$. 
\end{teo}

In \S\ref{s.teoB} we will expand on the strategy and it will be clear when it is possible to ensure the existence of Reeb surfaces (and not just generalized Reeb surfaces) in both. We will also show the stronger result that either $\cG$ is leafwise quasigeodesic in $\cF_1$, or $\cF_1$ has a generalized Reeb surface. We also deduce in Theorem \ref{teo.periodic} that two transverse minimal foliations in a manifold with non-solvable fundamental group, their intersection contains a closed leaf. This can be compared with \cite[Question 10.8]{FP-hsdff}. Note that if $\pi_1(M)$ is (virtually) solvable, then the result above does not hold as stated (see \cite[\S 6]{FP-hsdff}). 

As a consequence of Theorem \ref{teo.main}, we will obtain in \S \ref{ss.teo1} the following direct consequence:  

\begin{teo}\label{teo.mainph}
Let $f: M \to M$ be an oriented (chain)-transitive partially hyperbolic diffeomorphism in a manifold with exponential growth of fundamental group. Then, $f$ is a strong collapsed Anosov flow. 
\end{teo}

In  \S \ref{s.uniquelyintegrable} we will show the following result which allows to remove the orientability assumptions in some cases. Note also that a direct consequence is that under these assumptions all definitions of collapsed Anosov flows coincide (cf \S~\ref{ss.CAF}).

\begin{teo}\label{teo.mainorient}
Let $f: M \to M$ be a partially hyperbolic diffeomorphism and let $\hat f$ be a lift of an iterate of $f$ to a finite cover $\hat M$ which is an oriented partially hyperbolic diffeomorphism and so that every $\hat f$-invariant branching foliation is by Gromov hyperbolic leaves and intersections are leafwise quasgeodesic. Then $f$ preserves unique branching foliations $\cW^{cs}$ and $\cW^{cu}$. Moreover $f$ is a strong collapsed Anosov flow with respect to these branching foliations. 
\end{teo}

\section{Proof of Theorem A}\label{s.teoA}

\subsection{Proof of Theorem \ref{teo.mainph} assuming Theorem \ref{teo.main}}\label{ss.teo1} Let us first explain how to prove Theorem \ref{teo.mainph} using Theorem \ref{teo.main}.  

\begin{proof}[Proof of Theorem \ref{teo.mainph}]
By assumption, we know that leaves of the branching foliation are by Gromov hyperbolic leaves (see Theorem \ref{teo.chaintransitive}).
First suppose that $M$ is virtually solvable. 
Since $\pi_1(M)$ has exponential growth, then it cannot be nilpotent.

In addition since leaves are Gromov hyperbolic leaves the branching foliations do not have torus leaves, and it follows from \cite{HP-Solv} that the pair of branching foliations are uniformly equivalent to the stable and unstable foliations of a suspension Anosov flow in $M$. This proves that the intersection is leafwise quasigeodesic,
and Theorem \ref{teo-QGimpliesSCAF} shows that $f$ is
a strong collapsed Anosov flow.

If $M$ is not virtually solvable then using Proposition \ref{prop.noreebph} we know that the approximating foliations to the branching foliations intersect in a foliation which has no generalized Reeb surfaces. Thus Theorem \ref{teo.main} implies it is leafwise quasigeodesic in both foliations. 
Theorem \ref{teo-QGimpliesSCAF} again implies that 
$f$ is a strong collapsed Anosov flow.
\end{proof}

\subsection{Proof of Theorem A assuming Theorem \ref{teo.mainorient}}\label{ss.teoa}

Consider first $\hat f$ a lift of an iterate of $f$ to a finite (at most four fold) cover of $M$ on which every bundle is orientable and $\hat f$ preserves the orientation of the bundles. Thus, $\hat f$ is an oriented partially hyperbolic diffeomorphism. Proposition \ref{prop-liftcover} implies that $\hat f$ is also chain-transitive, and thus every $\hat f$-invariant branching foliation must be by Gromov hyperbolic leaves due to Theorem \ref{teo.chaintransitive}. Then we can apply Theorem \ref{teo.mainorient} and get that
$f$ is a strong collapsed Anosov flow, hence
 a collapsed Anosov flow, which is the conclusion of Theorem A. 

\begin{remark}
There are other ways to ensure that every branching foliation is by Gromov hyperbolic leaves: for instance, if $f$ is homotopic through partially hyperbolic diffeomorphisms to a (chain)-transitive one and $M$ has exponential growth of fundamental group (see for example \cite[Remark 4.9]{BFP}). Uniqueness of branching foliations does not hold in full generality as can be seen with examples made from the examples of \cite{HHU3} (see also \cite{Martinchich,HaPo-tori}). Still, even when the fundamental group of $M$ is not exponential, the chain recurrence hypothesis provides uniqueness of branching foliations (see \cite{HP-Solv}).
In fact, for such manifolds, in \cite{HP-Solv} it is proved that there are unique branching foliations
unless there are $cs$ or $cu$ tori (tori tangent to the $E^{cs}$ or
$E^{cu}$ bundles). If there are $cs$ or $cu$ torus then there is an
$f$ periodic one, and so (if it is a $cu$) torus, it is attracting and $f$ cannot be chain transitive.
Notice that $cs, cu$ tori can only occur in solvable manifolds, see \cite{CRHRHU}.
\end{remark}

\subsection{Periodic orbits}\label{ss.periodic}

Here we get the following consequence of Theorem \ref{teo.main}. The proof assumes familiarity with \cite[\S 6]{BFP}. 

\begin{teo}\label{teo.periodic}
Let $\cF_1, \cF_2$ transverse minimal foliations in a manifold with non-solvable fundamental group. Then, $\cG= \cF_1 \cap \cF_2$ has a closed leaf.  
\end{teo}

\begin{proof}
If $\cF_1 \cap \cF_2$ has a generalized Reeb surface, then, at least one boundary is closed, so by Theorem \ref{teo.main} we can assume without loss of generality that the flow is leafwise quasigeodesic. 

As in \cite[Proposition 6.9]{BFP} we can deduce that there is a closed lamination $\Lambda$ saturated by $\cF_1$ leaves consisting of weak quasigeodesic fans. 
A \emph{weak quasigeodesic fan}\footnote{See Figure \ref{f.qgfan} for what we call \emph{quasigeodesic fan}, weak quasi geodesic fans allows pairs of leaves to share \emph{both} endpoints.} is a leaf $L$ of $\wcF_1$ where all 
leaves of $\cG_L$ have a common ideal point, call it $p$. Then for any 
other $q$ in $S^1(L)$ there is at least one leaf of $\cG_L$ with ideal
points $q, p$. These leaves must have stabilizers which are either trivial or cyclic, and cannot all be trivial. Consider one such $L$ with non trivial
stabilizer generated by $\gamma$. Then $\gamma$ has to fix the 
common ideal point $p$ of all leaves in $\cG_L$. In addition since 
$\gamma$ is a hyperbolic isometry it fixes another point $q$ in $S^1(L)$.
Let $I$ be the compact interval (possibly degenerate) of leaves with ideal points $q, p$.
The end leaves of this interval are fixed by $\gamma^2$ leading to a closed leaf for $\cG$.
\end{proof}

 
 \section{Structure of Theorem B and reduction to other results}\label{s.teoB}
 
Theorem B splits into two statements (Theorems \ref{teo-injective} and \ref{teo-noninjective}) that together imply the result, but which give more information than the main result alone. 
 
 The starting point is \cite[Theorem 1.1]{FP-hsdff}:
 
 \begin{teo}\label{teo-Hsdff}
 Let $M$ be a closed 3-manifold whose fundamental group is not virtually solvable. Let $\cF_1$ and $\cF_2$ be transverse foliations by Gromov hyperbolic leaves in $M$, and let $\cG= \cF_1 \cap \cF_2$ be the
intersection foliation. Assume that restricted to every leaf $L \in \wcF_1$ and to every leaf $L$ $\in \wcF_2$ we have that the restriction $\cG_L$ of $\wcG$ to $L$ has Hausdorff leaf space. Then, $\cG$ is leafwise quasigeodesic in both $\cF_1$ and $\cF_2$. 
 \end{teo}

Hence the goal of this section is to show that
generalized Reeb surfaces are the source of every leafwise non-Hausdorff behaviour of $\cG$. Motivated by the examples in \cite{MatsumotoTsuboi} and in \cite{BBP},  which are of different nature from
each other, we separated the study in two cases. The first one explains the phenomena in \cite{BBP}, while the second case requires a different approach and produces a structure seen in the examples from \cite{MatsumotoTsuboi} (see also \cite[\S 7]{FPt1s}). 

\begin{teo}\label{teo-injective}
 Let $M$ be a closed 3-manifold. Let $\cF_1, \cF_2$ transverse foliations by Gromov hyperbolic leaves in $M$ and assume that for every $L \in \wcF_1$ and $E \in \wcF_2$ the intersection $L \cap E$ is connected (it could be empty). Then:
 \begin{enumerate}
 \item either there is a generalized Reeb surface in a leaf of $\cF_1$ or there
is a generalized Reeb surface in a leaf of $\cF_2$ (a priori both
can occur), 
 \item or the foliation $\cG = \cF_1 \cap \cF_2$ is leafwise quasigeodesic in both $\cF_1$ and $\cF_2$, 
 \item or there are bubble laminations $\Lambda_1 \subset \cF_1$ and $\Lambda_2 \subset \cF_2$. 
 \end{enumerate}
Moreover, if $M$ is hyperbolic or Seifert or if both $\cF_1$ and $\cF_2$ are minimal and $M$ does not have solvable fundamental group, then, the last option cannot happen. 
\end{teo}

\begin{remark} In \cite{BBP} many examples are introduced
with the features of Theorem \ref{teo-injective}.
In many examples of \cite{BBP} one starts with a closed 
hyperbolic $3$-manifold $M$ that fibers over the circle,
and a well chosen pseudo-Anosov suspension flow.
One blows up the stable and unstable
foliations along the singular orbits to produce laminations, and then
one fills the holes (solid tori) with appropriate foliations,
producing a pair of transverse foliations $\cF_1, \cF_2$.
The examples satisfy the hypothesis of \ref{teo-injective} that
intersections of pairs of leaves one in $\wcF_1$ and one in $\wcF_2$ are connected.
In \cite[\S 5]{BBP} examples are constructed which
are entirely symmetric with respect to the stable and unstable
foliations, and hence there are Reeb annuli of $\cG$ in leaves
of both $\cF_1$ and $\cF_2$. In \cite[\S 2 and 3]{BBP} the filling
is more involved and not symmetric: in this case there are
Reeb annuli in leaves of (say) $\cF_2$, but not in leaves
of $\cF_1$. In addition the leaf space $\cG_L$ for any $L$ in
$\wcF_1$ is Hausdorff, and one can prove the leaves of $\cG_L$
in $L$ are uniform quasigeodesics in $L$.
Details will appear in a forthcoming work.
It was explained to us afterwards by the authors of the paper \cite{BuTa}, that arguing as there provides a way to produce minimal examples with injective developing map and minimal foliations (consider a non-$\RR$-covered Anosov flow and its tilted foliation) to get examples with Reeb annuli. We believe (joint with Chi-Cheuk Tsang) that there are examples of transverse foliations with Reeb crowns.  
\end{remark}

The other case to consider is the one where leaves may intersect in several connected components, which forbids the leafwise quasigeodesic behavior of the intersected foliation, so we get:

\begin{teo}\label{teo-noninjective}
Let $\cF_1,\cF_2$ be two transverse foliations in a closed 3-manifold $M$ and assume that there is a leaf $L \in \wcF_1$ and a leaf $E \in \wcF_2$ so that the intersection $L \cap E$ has more than one connected component. Then, $\cG$ has a Reeb surface in both $\cF_1$ and $\cF_2$. 
\end{teo}

Note that we do not assume in this result that the leaves are Gromov hyperbolic nor anything on the topology of $M$. In fact, not even Reebless needs to be assumed but details of this last fact are left to the reader\footnote{If $\cF_1$ or $\cF_2$ have a Reeb component, the only way to put a foliation transverse to the Reeb component produces a Reeb annulus in the boundary torus.}. 

It is easy to see that Theorem \ref{teo-injective} and \ref{teo-noninjective} together imply Theorem \ref{teo.main}. Sections \ref{s.B1}, \ref{s.lemas} and \ref{s.B2} will be devoted to their proof. 
 
\section{Injectivity of the developing map}\label{s.B1}
The goal of this section is to prove Theorem \ref{teo-injective}. We will assume that $\cF_1, \cF_2$ are transverse foliations by Gromov hyperbolic leaves intersecting in a foliation $\cG= \cF_1 \cap \cF_2$ and the main assumption throughout this section is that given $\ell \in \wcG$ we have that if $L \in \wcF_1$ and $E \in \wcF_2$ are the leaves which contain $\ell$, then $\ell = L \cap E$. Equivalently, for every $L \in \wcF_1$ and $E \in \wcF_2$ the intersection $E\cap F$ has at most one connected component. We will say in this case that $\cG$ has \emph{injective developing map} (motivated by \cite{BarbotFP}). 

%
%
%
%
%
%

The proof starts by showing that if two non-separated leaves in a leaf of 
one of the foliations in $\mt$ have a deck transformation that sends them close by, then, they must return in the same local sheet of the other foliation. 

To produce such pairs of leaves and returns, we use the fact that some of the results of \cite{FP-hsdff}, particularly \cite[Proposition 4.1]{FP-hsdff}, can be extended to the case where intersections between leaves in the universal cover are connected. The hypothesis in \cite{FP-hsdff} is stronger, namely
that the leaf space of $\cG_L$ is Hausdorff for every leaf $L$ in
either $\wcF_1$ or $\wcF_2$. 
Here it is exactly the Hausdorff property that we are trying to
obtain.
Then, part of the general strategy devised in \cite[\S 1.1]{FP-hsdff} will be developed in the case of this section until we achieve \emph{landing} and \emph{small visual measure}:  two geometric properties that require the Gromov hyperbolicity of leaves. Note that small visual measure may fail, giving rise to \emph{bubble laminations}, but under certain assumptions we will be able to exclude them. 

Next, using these tools we will exploit the non-Hausdorfness to produce generalized Reeb surfaces. This second part, which is subtler, diverges entirely from the 
arguments in \cite{FP-hsdff}. In a forthcoming paper \cite{FHP} that supersedes the present one, we shall give a somewhat different proof with some intersection with this one but which avoids completely the use of Gromov hyperbolicity of leaves, landing and small visual measure.

An important consequence of injectivity of the developing map is the following:

\begin{lema} \label{lem.same}
Let $L$ be a leaf of $\wcF_1$ and $\ell_1, \ell_2$ two distinct
leaves of $\cG_L$ which are non separated from each other.
Let $E_i \in \wcF_2$ the leaf containing $\ell_i$. Then $E_1, E_2$ are distinct
and non separated from each other.
\end{lema}

\begin{proof} 
There are $e_n$ in $\cG_L$ converging to $\ell_1 \cup \ell_2$ 
(and maybe other leaves of $\wcG$ as well). Therefore $E_1, E_2$
are non separated from each other. The important point is that
under the running hypothesis in this section, $\ell_i = E_i \cap L$,
hence $E_1 \not = E_2$.
\end{proof}

\subsection{Returns in same local sheets}\label{ss.returns}

Here we consider $L \in \wcF_1$ and two non separated rays $r_1, r_2$ of $\cG_L$. One can consider a curve $\alpha: [0,1] \to L$ whose endpoints are in $r_1,r_2$ respectively, it is transverse to $\cG_L$ at the endpoints and such that the interior does not intersect the leaves of $r_1$ and $r_2$. 

Let $\eps>0$ be small so that $\alpha|_{[0,\eps)}$ is transverse to $\cG_L$. It follows from the non-separation, that if $\ell_t$ is the leaf of $\cG_L$ intersecting $\alpha(t)$ for $t \in (0,\eps)$ then it defines an open arc $c_t$ whose closure is a compact segment of $\ell_t$ whose endpoints are $\alpha(t)$ and $\alpha(1-t)$ (this is achieved by a reparametrization of $\alpha$ if necessary). One then defines: 

\begin{equation}\label{eq:Ut}
U_{r_1,r_2} = \bigcup_{t \in (0,\eps_0)} c_t , 
\end{equation}

\noindent which is an open topological disk in $L$ whose boundary consists of $c_{\eps_0}$, two subarcs of $\alpha$, the non-separated rays of $r_1,r_2$ starting at $\alpha(0),\alpha(1)$ and possibly other (at most countably many) leaves $\{\ell_i\}_i$ of $\cG_L$ which are non separated from both of the rays $r_1,r_2$.

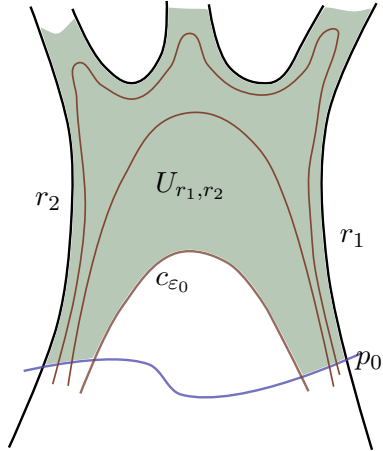
\begin{figure}[htbp]
\begin{center}
\includegraphics[scale=0.56]{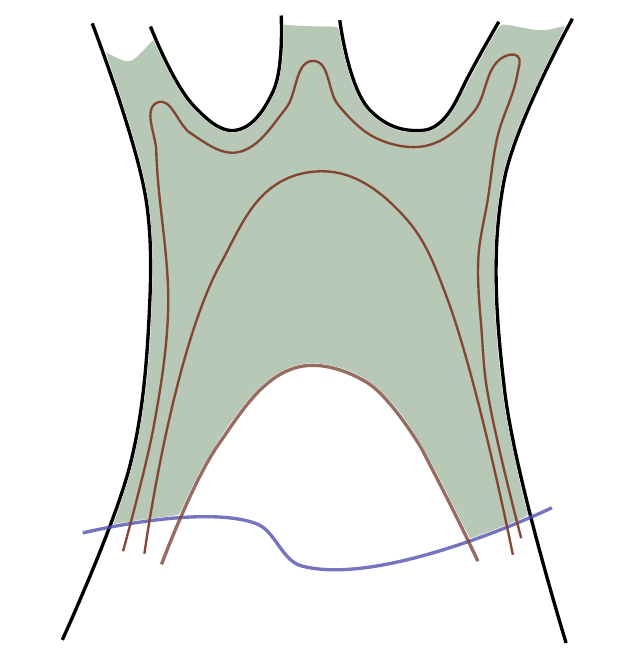}
\begin{picture}(0,0)
\put(-36,86){$r_1$}
\put(-150,100){$r_2$}
\put(-30,40){$p_0$}
\put(-105,105){$U_{r_1,r_2}$}
\put(-105,70){$c_{\eps_0}$}
\end{picture}
\end{center}
\vspace{-0.5cm}
\caption{{\small The shaded region depicts $U_{r_1,r_2}$.}\label{f.ur1r2}}
\end{figure}

We denote by $NS(r_1,r_2)$ the set of complete leaves $\{\ell_i\}_i$ of $\cG_L$ which are in $\partial U_{r_1,r_2}$ (so we are not
including $r_1, r_2$). Note that this set could be empty. 

We want to show that if there are points which have a close return under some deck transformation, then the returns should belong to the same local sheet of the other foliation. 
A return of $z \in r_1$ near $p$ means $\gamma z$ is near $p$ for some
$\gamma$ a deck transformation (an element of $\pi_1(M)$).
Let us formalize this in the following statement. We let $\eps_0$ to be the size of local product structure, see \S~\ref{ss.transverse}. 

\begin{prop}\label{prop-returnlocalsheet}
Let $\cF_1, \cF_2$ be transverse foliations with the injective developing map property. Let $L \in \wcF_1$ so that $\cG_L$ has rays $r_1, r_2$ which are non separated from each other and $\eps<\eps_0$. Assume there is $\gamma \in \pi_1(M)$ and points $p,p' \in  r_1, q \in r_2$ and $q' \in r \in \{r_2\} \cup NS(r_1,r_2)$ so that
\[ d_{\mt}(\gamma p', p) < \eps \text{ and  } d_{\mt}(\gamma q', q)< \eps. \]
Then, if $E_1,E_2 \in \wcF_2$ are such that $r_1 \subset E_1 \cap L$ and $r_2 \subset E_2 \cap L$ then one has  that $\gamma p' \in E_1$ and $\gamma q' \in E_2$. 
\end{prop}

All we will need to use later in the article is that $\gamma p' \in E_1$, but we need to know that the point $q'$ returns close to $q$ to ensure this. This proposition does not use the geometry of leaves, but to produce 
these simultaneous returns is where Gromov hyperbolicity of the leaves of $\cF_1$ will be used crucially. 

\begin{proof} 
Consider small neighborhoods $B_p, B_{p'}, B_q, B_{q'}$ in $L$ around the points $p,q,p',q'$ which are trivially foliated by $\cG_L$ and such that $r_1$ separates $B_p, B_{p'}$ in exactly two connected components and $r_2, r$ separate $B_q, B_{q'}$ respectively in exactly two connected components.

 Denote by $\alpha$ a curve in $L$ joining $p$ and $q$ and intersecting $r_1, r_2$ only at the extremes (and not intersecting other non-separated leaves from $r_1,r_2$). Similarly, one can consider $\alpha'$ joining $p', q'$ and intersecting $r_1, r_2$ only at the extremes. We denote by $B_p^+, B_{p'}^+, B_q^+, B_{q'}^+$ the connected components of $B_p, B_{p'}, B_q, B_{q'}$ minus the corresponding rays $r_1,r_2, r$  which intersect $\alpha$ or $\alpha'$. By the non-separation property of the rays $r_1$ and $r_2$, we can assume without loss of generality that every leaf of $\cG_L$ intersecting $B_p^+$ intersects also $B_q^+$ and similarly for $B_{p'}^+$ and $B_{q'}^+$.

Let now $E_1$ and $E_2$ be the leaves of $\wcF_2$ containing respectively $r_1$ and $r_2$ and let $E_3$ the leaf through $q'$ (it could be equal to $E_2$ if $q' \in r_2$). We want to show that $\gamma E_1 = E_1$ and $\gamma E_3 = E_2$. Let $\ell_1=L\cap E_1$ and $\ell_2 = L \cap E_2$ 
be the leaves of $\cG_L$ containing $r_1$ and $r_2$ respectively,
and let $\ell_3 = L \cap E_3$.
Note here that the injectivity property ensures that $E_1 \neq E_2$, $E_1 \neq E_3$, so $\gamma E_1 \neq \gamma E_3$.  Moreover, $\gamma E_1$ contains $\gamma p'$ which is near
$p$, hence $\gamma E_1$ intersects $L$. 

We proceed by contradiction assuming that $\gamma E_1 \neq E_1$. Suppose first that $\hat \ell_1=\gamma E_1\cap L$ intersects $B^+_p$, and hence it  intersects $\alpha$ in $B_p^+$ in a point $\hat p$ very close to $p$ (recall the orientation hypothesis). This point $\hat p$ is connected to $\gamma p'$ by a curve $\tau_1:[0,\eps] \to \mt$ either in $L$ or transverse to $\wcF_1$. By the choice of $B_p^+$ we know that $\hat \ell_1$ also intersects $B_q^+$. If $\hat \ell_2 =\gamma E_3 \cap L$, then
$\gamma q' \in \gamma E_3$ is very near $q$, hence $\hat \ell_2$ intersects
$B_q$, and  we have that both $\hat \ell_2$ and $\hat \ell_1$ intersect $B_q$.
It follows $\gamma E_1, \gamma E_3$ intersect a common transversal to $\wcF_2$.
and therefore
so do $E_1, E_3$. This is a contradiction to Lemma \ref{lem.same}.

Suppose now that $\hat \ell_1$ does not intersect $B^+_p$. 
In this case we use $\gamma^{-1}$ and do the following:
let $\hat \ell_1' = \gamma^{-1}(E_1) \cap L$. Since $\hat \ell_1$ does
not intersect $B^+_p$, it follows that $\hat \ell_1'$ intersects
$B^+_{p'}$. Hence it also intersects $B^+_{q'}$. 
In the same way as in the previous argument $\gamma^{-1}(E_2)
\cap L$ intersects $B_{q'}$. So $\gamma^{-1}(E_1), \gamma^{-1}(E_2)$
both intersect $B_{q'}$ and hence a common transversal,
again a contradiction. Therefore $\gamma E_1 = E_1$.

Assuming that $\gamma E_3 \neq E_2$ one can produce a similar contradiction. 
\end{proof}

\subsection{Pushing} 

Here we use that the arguments of \cite[Proposition 4.1]{FP-hsdff} are valid in the
context of this section. The proof is exactly the same as in \cite{FP-hsdff} and this can be seen by looking at \cite[Figure 4]{FP-hsdff} which shows how if one pushes a leaf of $\wcG$ transversely to say $\wcF_1$ along a leaf
of $\wcF_2$, and this \emph{splits} into two leaves of $\wcG$, then this forces a leaf of $\wcF_2$ to intersect one leaf of $\wcF_1$ in at least two distinct connected components.

\begin{prop}[Pushing]\label{prop-pushing}
Let $\cF_1, \cF_2$ be two transverse foliations so that $\cG=\cF_1 \cap \cF_2$ has injective developing map. Consider a leaf $\ell_0$ the foliation $\wcG = \wcF_1 \cap \wcF_2$ where $\ell_0 = L_0 \cap E_0$ with $L_0  \in \wcF_1$ and $E_0 \in \wcF_2$. Let $\tau_1, \tau_2 : [0,t_0] \to \mt$ be transversals to $\wcF_1$ satisfying the following conditions:
 \begin{enumerate}
\item $\tau_1(t)$ and $\tau_2(t)$ belong to the same leaf $L_t$ of $\wcF_1$ for all $t$, \item $\tau_1(0)$ and $\tau_2(0)$ belong to $\ell_0$, and, \item $\tau_1(t), \tau_2(t)$ belong to $E_0$ for all $t$. 
\end{enumerate} 
Then, for every $t$ we have that $\tau_1(t)$ and $\tau_2(t)$ are the endpoints of a compact arc $c_t$ contained in $\ell_t = L_t \cap E_0$.  
\end{prop}

\begin{proof}
The statement of the conclusion is in fact what is proven in \cite{FP-hsdff}: if the conclusion is not true, then there is a pair of leaves $L \in \wcF_1, E \in \wcF_2$ so that
$L, E$ intersect in more than one component\footnote{ In \cite[Lemma 2.5]{FP-hsdff} this fact is then
used to get that the leaf space of the intersection cannot be Hausdorff. }, which is
what this proposition claims.
\end{proof}

\subsection{Landing} 

Since leaves of $\cF_1$ are by Gromov hyperbolic leaves, we can compactify each leaf with a circle at infinity. Given $L \in \wcF_1$ we consider $\overline L= L \cup S^1(L)$ its compactification. See for instance \cite[\S 2.3]{FP-hsdff} for more details. Thus, given a properly embedded ray $r: [0, \infty) \to L$ we can consider its boundary $\partial r$ to be the set of accumulation points of $r(t)$ as $t \to \infty$ inside $S^1(L)$ (since $\wcG$ restricted to $L$ is a foliation in $L$, then the infinite arc $r(t)$ is proper, and it does not accumulate in $L$). We say that the ray $r$ \emph{lands} if $\partial r$ is a unique point. See \cite[\S 2.5]{FP-hsdff} for more on this. 

Applying the fact that in \cite{FP-hsdff} only the property\footnote{In \cite[Remark 5.5]{FP-hsdff} it is explained how if one assumes Hausdorff leaf space inside the leaves of $\wcF_1$ it is easier to get a contradiction.} of Proposition \ref{prop-pushing} is used to get landing, we get: 

\begin{teo}\label{teo.landing}
Let $\cF_1,\cF_2$ be foliations by Gromov hyperbolic leaves and assume that $\cG=\cF_1 \cap \cF_2$ has injective developing map. Then, every ray of $\wcG$ lands in its corresponding leaf of $\wcF_1$ (resp. $\wcF_2$). 
\end{teo}

\begin{proof}
The proof of \cite{FP-hsdff} applies verbatim as it only uses the pushing property given by Proposition \ref{prop-pushing}. (See \cite[Remark 5.5]{FP-hsdff} and the comments after \cite[Theorem 5.1]{FP-hsdff}.) 

Let us briefly explain the strategy for completeness. 

The first two steps, which are \cite[Lemma 5.2]{FP-hsdff} and  \cite[Lemma 5.3]{FP-hsdff}, are general properties of foliations with Gromov hyperbolic leaves and one dimensional subfoliations of them. Together they state the
following: suppose that a ray $r$ of $\wcG$ in a leaf $L'$ of $\wcF_1$ does
not land in $S^1(L')$ and limits to a non degenerate interval $J$
in $S^1(L')$. Then ``zoom in" in $L'$ towards an interior point
of $J$, and in the limit produce 
a minimal sublamination of $\cF_1$, so that
$\wcF_1$ contains a leaf $L$ of the lifted sublamination with $L$ having  non-trivial stabilizer, and so that this lamination is related with
the failure of landing of the initial ray. Then consider a geodesic in $L$
projecting to a closed geodesic in a leaf of $\cF_1$.
This is accumulated more and more by arcs of deck translates of $\eta$ which
are limiting to big intervals in $S^1(L'')$ for $L''$ in $\wcF_1$ near $L$.

Then, Proposition \ref{prop-pushing} is used to get a result (see \cite[Lemma 5.4]{FP-hsdff}) obtaining enough arcs (by pushing to $L$ the translates of the ray constructed above). This is where the pushing property is used
in our context. These ``pushed" arcs in $L$ cannot intersect each other
transversely,  because they belong to leaves of the foliation $\wcG$ restricted to $L$. Finally, this is shown to be a contradiction using only the structure that had been obtained in \cite[Lemma 5.4]{FP-hsdff} and applying the elements of the stabilizer of $L$. 
This completes the proof.
\end{proof}

\subsection{Small visual measure and bubble laminations}\label{ss.svm}

Here, we will define the small visual measure property and show that in the case of injective developing map, 
$\cG$ has the small visual measure property unless there
are bubble laminations.

We say that a one dimensional subfoliation $\cG$ of a foliation $\cF$ by Gromov hyperbolic leaves has the \emph{small visual measure property} if for every $\eps>0$ there is a uniform $R>0$ so that if $L \in \wcF$, $x \in L$ and $c$ is a segment (it could be a ray or full leaf) of a leaf of $\widetilde{\cG}$ contained in $L$ and totally outside the $R$-neighborhood of $x$ in $L$, then, the shadow of $c$ in $S^1(L)$ has visual measure less than $\eps$ when measured from $x$. We refer the reader to \cite[\S 7]{FP-hsdff} for a precise definition of visual measure in Gromov hyperbolic leaves. For concreteness, if one uses the Candel metric (which makes every leaf isometric to a hyperbolic plane), this means that if $I$ is the interval of vectors tangent to $x$ for which the geodesic ray starting at $x$ with this tangent vector intersects $c$ (i.e. the shadow of $c$), then the interval $I$  has angle less than $\eps$. Note that the small visual measure property implies the landing of rays, but landing is used to establish it. 

An important property of small visual measure is the following (see \cite[Proposition 7.2]{FP-hsdff}): 

\begin{lema}\label{lem.svm} 
Assume that $\cG$ has the small visual measure property on $\cF$ then, there is a uniform constant $a_0>0$ so that for every arc or ray $c$ of a leaf $\ell \in \widetilde{\cG}$ contained in a leaf $L \in \wcF$ we have that the $a_0$-neighborhood of $c$ contains the geodesic segment (or ray) joining the endpoints of $c$.  
\end{lema}

Adapting some arguments from \cite[\S 8]{FP-hsdff} we show that:

\begin{prop}\label{prop.svm}
Let $\cF_1, \cF_2$ be transverse foliations by Gromov hyperbolic leaves so that $\cG= \cF_1 \cap \cF_2$ has injective developing map. If $\cG$ does not have the small visual measure property on the foliation $\cF_1$, then $\cF_1$ contains a sublamination $\Lambda_1$ so that for every $L \in \widetilde \Lambda_1$ there is a point $\xi(L)$ so that every ray of a leaf of $\cG_L$
verifies that it lands at $\xi(L)$, i.e. $\Lambda_1$ is a bubble lamination. Moreover, the point $\xi(L)$ varies continuously with $L$. 
\end{prop}

\begin{proof}

This is done in \cite[\S 8.1]{FP-hsdff}. First, in \cite[Lemma 8.2]{FP-hsdff} it is shown that failure of small visual measure provides points $y_n$ and arcs $c_n$ of $\cG_{L_n}$ (with $y_n \in L_n$) so that seen from $y_n$ the arc $c_n$ projects to almost all of $S^1(L_n)$. Taking deck transformations one can assume that $y_n$ converges to $y_\infty$ and using the pushing property (Proposition \ref{prop-pushing}) one can see that if $L_\infty \in \wcF_1$ is the leaf through $y_\infty$, then, it contains arcs whose landing points belong to an arbitrarily small interval in $S^1(L_\infty)$ (see \cite[Lemma 8.3]{FP-hsdff}). This shows that $L_\infty$ has the desired property, and since the set of leaves with this property is a closed subset of $\mt$, and it is $\pi_1(M)$ invariant, one gets a minimal lamination $\Lambda_1$ (see \cite[Lemma 8.4]{FP-hsdff}). This lamination verifies that every leaf has cyclic fundamental group (see \cite[Corollary 8.6]{FP-hsdff}) and one can prove a continuity property of the point $\xi(L)$ for leaves of the lamination by using an appropriate topology in the union of circles at infinity (see \cite[Lemma 8.7]{FP-hsdff}). 

This completes the proof of this proposition. 
\end{proof}

\subsection{Non-separated leaves of $\wcG$} 

Here we assume that there is small visual measure of $\cG$ in $\cF_1$. Suppose that there is a leaf $L \in \wcF_1$ and two non-separated rays $r_1, r_2 \in \cG_L$.  We will use the setup considered in  \S \ref{ss.returns}, in particular, we consider a curve $\alpha: [0,1] \to L$ joining points $p_0 \in r_1$ and $q_0 \in r_2$ and we will consider the leaves $NS(r_1,r_2)$ that are in the boundary of the set $U_{r_1,r_2}$ defined in equation \eqref{eq:Ut}. For $t \in (0,\eps)$ we have the arcs $c_t$ which approach as $t \to 0$ all leaves of $NS(r_1,r_2) \cup \{r_1,r_2\}$.

We will show the following useful property that allows to apply Proposition \ref{prop-returnlocalsheet}. 

\begin{prop}\label{prop-escapingpoints}
There is a constant $a_1>0$ so that, for every $n>0$ there is a point  $p_n \in r_1$ with $d_L(p_0,p_n)>n$, and a point $q_n \in r_n \in \{r_2\} \cup NS(r_1,r_2)$ so that $d_L(p_n,q_n)\leq a_1.$
\end{prop}

\begin{proof}
This is an application of the small visual measure property: Let $\xi$ in $S^1(L)$ be the landing point of $r_1$. 
Suppose first that there is a leaf $r'$ in $\{ r_2 \} 
\cup NS(r_1,r_2)$ so that $r'$ has a ray (still denoted by $r'$)
with landing point $\xi$. By Lemma \ref{lem.svm}, 
for every point in a geodesic ray in $L$ with ideal point $x$,
then both $r'$ and $r_1$ have points a fixed bounded distance
from such a point. Hence we can find $p_n$ and $q_n \in r'$ in this case.

Suppose from now on that no ray in $\{ r_2 \} \cup NS(r_1,r_2)$
has landing point $\xi$. Parametrize $r_1$ as $\{ v_s, s \in \RR_{\geq 0} \}$.
If the proposition is not true then 

$$d_s \ := \ d_L(v_s, \{ r_2 \} \cup NS(r_1,r_2)) \to \infty$$

We have assumed that no element in $\{ r_2 \} \cup NS(r_1, r_2)$ has
landing point $x$. Let $I$ be the interval of $S^1(L)$ with
endpoints $\xi$, and
the ideal point $y$ of $r_2$, and so that the arcs $c_t$ of leaves of $U_{r_1,r_2}$ as above 
have points converging to every point in this interval.

First we assume that for every non trivial subinterval of $I$,
there is an element $\ell$ in $NS(r_1, r_2)$ with landing point in this
interval.

There is a lower bound of the visual measure
of $I$ from any point $v_s$ because $I$ is non
trivial and $\xi$ is the landing point of $r_1$, call this bound $2a_1>0$. Given $v_s$ choose $z$ in $I$ so that
the visual measure as seen from $v_s$ of the subinterval $I'$ of $I$ from $y$ to $z$
has visual measure $> a_1$. Now choose $\ell$ in $NS(r_1,r_2)$
with an ideal point not in $I'$.
Then for big enough $t$, \  $c_t$ has a subinterval $c'_t$ with 
an endpoint very close to $q_0$ in $r_2$ and another very
close (in the compactification $L \cup S^1(L)$) to $z$.
The visual measure  of $c'_t$ from $v_s$ can be chosen
to be $> a_1$. Our assumption was that
distance from $v_s$ to $c'_t$ could be chosen to increase
without bound as $s \to \infty$. This contradicts the small visual
measure property.

Finally, suppose that some non degenerate subinterval $I'$ of $I$ with one
ideal point $\xi$ does not have any ideal points of $\ell$ with $\ell$
in $NS(r_1,r_2)$. Then $c_t$ has subintervals very close to $I'$ in
$L \cup S^1(L)$, and they can also be chosen to have distance from $v_s$
in $L$ going to infinity. This also contradicts the small visual
measure property.

This finishes the proof.
\end{proof}

\begin{convention}\label{convention-rays} 
Up to changing $p_0,q_0$ for $p_k,q_k$ and taking subsequences of the points $p_n, q_n$ from the previous proposition we can assume that for every $n$ the points $p_n$ project into an $\eps_1/10$-neighborhood of the projection of $p_0$ and the same for $q_n$ and $q_0$. Here $\eps_1$ is much smaller than the local product structure constant, and verifies that an arc of length $10a_1$ in a leaf of $\wcF_1$ or $\wcF_2$ can be pushed to nearby leaves having some point at distance less than $\eps_1$.  In particular, we can assume that for every $n >0$ there is a deck transformation $\gamma_n \in \pi_1(M)$ so that $\gamma_n p_n$ belongs to an 
$\eps_1/10$-neighborhood of $p_0$, 
and $\gamma_n q_n$ in an $\eps_1/10$-neighborhood of $q_0$. Note that by Proposition \ref{prop-returnlocalsheet} we know that $\gamma_n p_n$ belongs to $E_1$ the leaf of $\wcF_2$ containing $r_1$ and also $\gamma_n q_n \in E_2$ where $E_2 \in \wcF_2$ is the leaf through $r_2$. 

We can consider transversals $\tau_1, \tau_2 : [-1,1] \to \mt$ to $\wcF_1$ so that $\tau_1(0)=p_0$ and $\tau_2(0)=q_0$, so that $\tau_1(t)$ and $\tau_2(t)$ belong to the same leaf of $\wcF_1$ for every $t\in [-1,1]$, these arcs have length smaller than $\eps_1/10$, each arc is contained in a single leaf of $\wcF_2$ and curves of length $\leq a_1$ starting close to $p_0$ can be pushed to the leaf through $\tau_1(t)$ for every $t \in [-1,1]$ and produce a curve which is $\eps_0$ close to the original one.
\end{convention}

\subsection{Trivial foliation}

We first want to produce a transversal to the returns of the ray $r_1$ in the same plaque of $\wcF_2$, and show that when looked in the foliation $\wcF_2$ all leaves of $\wcG$ intersecting this transversal form a trivially foliated band.

As in Convention \ref{convention-rays}, denote by $E_1$ the leaf of $\wcF_2$ containing $r_1$ and $E_2 \in \wcF_2$ the one containing $r_2$. Also consider $\tau_1, \tau_2$ as in Convention \ref{convention-rays}. 

Let $L_t$ be the $\wcF_1$ leaf of $\tau_1(t)$.

For every $n>0$ we know that $\gamma_n p_n$ is in an $\eps_1/10$-neighborhood of $p_0$ inside $E_1$ and we can assume without loss of generality that $\gamma_n p_n \in \tau_1([-1,1])$ and similarly\footnote{We will not actually use the returns $\gamma_n q_n$ too much, but note that their existence is crucial to ensure that we can apply Proposition \ref{prop-returnlocalsheet} and thus get that $\gamma_n p_n \in E_1$.} $\gamma_n q_n \in \tau_2([-1,1])$. 

For every $n>0$, let $t_n \in [-1,1]$ so that $\gamma_n p_n = \tau_1(t_n)$.
If $t_n = 0$ for some $n$  then the projection $\pi(r_1)$ is a closed curve
in a leaf of $\cF_1$ and by Proposition \ref{prop-genRS} the curve $\pi(r_1)$ is the boundary of a generalized Reeb surface. 

If on the other hand $t_n \not = 0$ for all $n$ (without loss of generality we assume $t_n>0$ for all $n$) we consider $u_n$ to be the closure of the infinite
ray contained in $\gamma_n r_1 \setminus \{ \gamma_n p_n \}$, and let $B_{1,n}$ be the closure of the connected component of $E_1 \setminus (\tau_1([0,t_n]) \cup r_1 \cup u_n)$ not containing $\tau_1(1)$. A key property we will use is:

\begin{lema}\label{lem-trivialsides}
The set $B_{1,n}$ is trivially foliated by infinite rays of $\cG_{E_1}$, and the only leaves of $\cG_{E_1}$ that intersect $B_{1,n}$ are contained $L_t \cap E_1$ with $t \in [0,t_n]$.
\end{lema}
\begin{proof}

We assume by way of contradiction that there are non-separated leaves $\ell_1, \ell_2$ in $B_{1,n}$. Without loss of generality, we can assume that one of the leaves, say $\ell_1$, intersects $\tau_1([0,t_n])$ because  otherwise, the union of the rays of leaves through $\tau_1([0,t_n])$ would be open and closed in $B_{1,n}$ and thus the set would be trivially foliated. 

\begin{figure}[htbp]
\begin{center}
\includegraphics[scale=0.69]{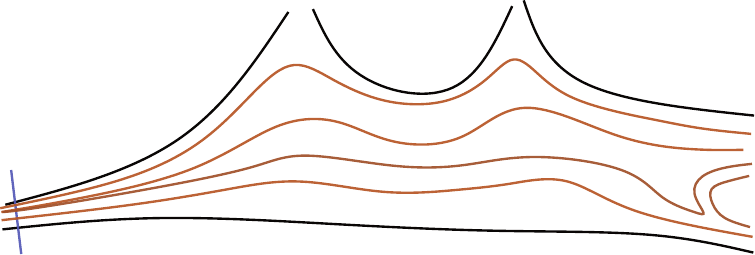}
\begin{picture}(0,0)
\put(-130,66){$\ell_2$}
\put(-190,60){$\ell_1$}
\put(-130,0){$r_1$}
\put(-299,24){{\small $\tau_1([-1,1])$}}
\end{picture}
\end{center}
\vspace{-0.5cm}
\caption{{\small If the region $B_{1,n}$ is not trivially foliated, one gets some non separated leaves.}\label{f.trivialfol}}
\end{figure}
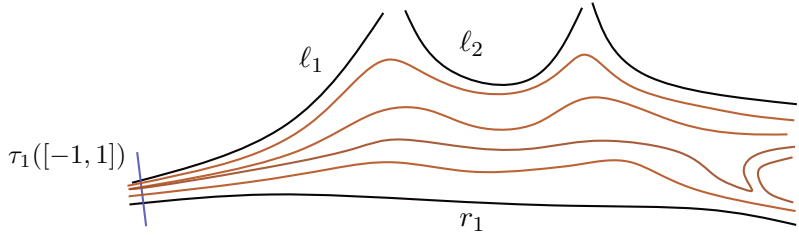

Let $H_i$ in $\wcF_1$ so that $\ell_i = H_i \cap E_1$. Notice that
by Lemma \ref{lem.same} one has that $H_1 \neq H_2$.
Since $\ell_1$ intersects $\tau_1([0,t_n])$ it follows that $H_1= L_t$ for some $t \in [0,t_n]$. Let $L_s$ with $s \in (0,t_n)$ very close to $t$ and so that $\ell_s= L_s \cap E_1$ has points $z_1$ and $z_2$ which are $\eps_1/10$ close to  $\ell_1$ and $\ell_2$ respectively, with $\eps_1$ as in Convention \ref{convention-rays}.

Since $s \in [0,t_n]$, and by the choice of $\eps_1$,
 it follows that $L_s$ intersects $E_2$.
Now, consider $\ell_s'  = L_s \cap E_2$.
Notice that $\ell_s'$ and $\ell_s$ are contained in $L_s$.
We claim that they are distinct, and that they are non separated in 
$\cG_{L_s}$.
The first part is obvious, since they are contained respectively
in $E_1, E_2$ which are distinct.
For the second part,  by the choice of $\eps_1$, we can consider a continuation of the arc $\alpha$ ($\alpha$ is the arc from $p_0$ to $q_0$ in $L$)
to an arc $\alpha_s$ in $L_s$ from $\tau_1(s)$ to $\tau_2(s)$.
We may assume that $\alpha_s$ is transverse to $\cG_{L_s}$ near
the endpoints. Recall that $r_1, r_2$ are non separated in
$\cG_L$, so there is a sequence of arcs
$b_n$ converging to $r_1 \cup r_2$ and maybe 
other leaves as well. 

Let $w_n, y_n$  points in $b_n$ converging
to $\tau_1(0), \tau_2(0)$.
The points $w_n, y_n$ push up along their common leaf of $\wcF_2$ to
points in $\alpha_s$. By Proposition
\ref{prop-pushing} the entire arc $b_n$ pushes to arcs $b_{n,s}$ in $L_s$, 
so that as $n \to \infty$, then $b_{n,s}$ converges to leaves
containing $\tau_1(s), \tau_2(s)$ respectively.
This shows that these leaves (which are $\ell_s$ and $\ell_s'$ respectively) of $\cG_{L_s}$ are distinct
and non separated from each other in $\cG_{L_s}$ as we wanted to show.

Now applying Proposition \ref{prop-pushing} in the other direction, 
we see that for any $n$, the full arc $b_{n,s}$ must push along $E_n$ (the leaf of $\wcF_2$ so that $L \cap E_n$ contains $b_n$ and thus $L_s \cap E_n$ contains $b_{n,s}$) 
to arcs in $H_1$ (which is equal to $L_t$).
On the other hand $\ell_s$ contains $z_2$, so there
are points in $b_{n,s}$ (for $n$ large enough) which are arbitrarily close to $z_2$. Recall that $z_2$ is $\eps_1/10$ close to $H_2$.
Pick $n$ sufficiently big and $v_n$ in $b_{n,s}$ which is sufficiently close
to $z_2$, so $v_n$ pushes along $E_n$ to a point in $H_2$ near a 
point in $\ell_2$.
This contradicts the fact that the entire arc $b_{n,s}$ pushes along 
$E_n$ to an arc in $H_1$, and $H_1$ is non separated from $H_2$.
This finishes the proof.
\end{proof}

\subsection{The limit generalized Reeb surface} 

Throughout this subsection we are considering rays $r_1, r_2$ of $\cG_L$
which are distinct, but non separated from each other.
In this section we will show that the ray $r_i$, when projected to $M$, either
is a boundary component of a Reeb surface or a Reeb crown in
the foliation $\cF_1$, 
or is asymptotic to a closed leaf of $\cG$. 
We keep the notation from the previous sections, in particular, from Convention \ref{convention-rays}. 

For this, the main property we will show is that: 

\begin{prop}\label{lem-converge}
Let $i = 1$ or $i = 2$. Then $\pi(r_i)$ is either a closed 
curve in $M$ or it is asymptotic to a closed curve in $M$.
\end{prop}

We will work with $i=1$ as the statement is symmetric and we assume throughout that $\pi(r_1)$ is not closed (else, the statement is immediate). 

Fix some large $n$ and recall that $\gamma_n L = \gamma_n L_0 = L_{t_n}$ for some $t_n \in [-1,1]$. It follows that there is a curve $\eta_n: [0,1] \to E_1$ transverse to $\cG_{E_1}$ contained in $B_{1,n}$ joining $p_n$ with $\gamma_n p_n$. Iterating by $\gamma_n$ we obtain a transversal

$$\hat \eta_n \ \ := \ \ \bigcup_{k \in \ZZ} \gamma_n^k \eta_n$$ 

\noindent
to $\cG_{E_1}$. This  produces an embedding of $\RR$ into the leaf space of $\cG_{E_1}$ which is invariant under $\gamma_n$ and on which $\gamma_n$ acts without fixed points. Denote by $\Theta_n$ this subset in the leaf space of $\cG_{E_1}$. 

%
%
Recall that by Lemma \ref{lem-trivialsides} the set   $B_{1,n}$ is trivially foliated, and
so is the closure $A_n$ of the non compact component of $B_{1,n} \setminus
\{ \eta_n \}$. Define 

$$D_n = \bigcup_{i \in \ZZ} \gamma^i_n A_n.$$ 

Since the foliation is group invariant, we deduce that $\gamma_n \Theta_n = \Theta_n$ and $\gamma_n D_n = D_n$ as well. We remark that $\gamma_n$ acts freely on $\Theta_n$.

Note that for $m$ different from $n$ we can also construct $\Theta_m$ and we have that $\Theta_m$ and $\Theta_n$ must coincide in a neighborhood of $L \cap E_1$. We will use here that $\pi(E_1)$ is not a torus (because $\cF_2$ is by Gromov hyperbolic leaves). 

\begin{lema}\label{lema-thetanm}
The sets $\Theta_n$ and $\Theta_m$ coincide. 
\end{lema}

\begin{proof}
The sets $\Theta_n,
\Theta_m$ are open intervals in the leaf space which intersect.
The intersection $\cI$ is a common subinterval of each one, which we want to
prove that it is the full set $\Theta_n = \Theta_m$. 
Suppose this is not the case.

There are two possibilities for a given finite end of $\cI$, either it is also an end of $\Theta_n$ or $\Theta_m$, or it limits to points $x \in \Theta_n$ and $y \in \Theta_m$ with $x \neq y$. 

Suppose first that
$\cI$ has an end that limits to leaves $x$ in $\Theta_n$ and $y$ in
$\Theta_m$, so that $x \not = y$. In other words, $x, y$ are
non separated from each other in the leaf space of $\cG_{E_1}$,
and they are limits of the interval $\cI$. That means that
$\Theta_n, \Theta_m$ branch away from each other at this point: there is a sequence of leaves $\ell_k$ in $\cI$ so that $\ell_k$ converges to both $x \in \Theta_n$ and $y \in \Theta_m$ (and maybe other leaves).

Now we use that $\cG_{E_1}$ is oriented. Without loss of
generality assume that $y$ is non separated from $x$ in the
positive direction from $x$, this means the following: let $u \in x$ and $v \in y$
respectively and $\ell_k$ as defined above. Let $u_k, v_k$ in
$\ell_k$ converging to $u, v$ respectively. Then $y$ is in the
positive direction from $x$ means that for sufficiently big
$k$, then $v_k$ is in the positive side of $u_k$ in $\ell_k$.

Up to changing the boundary leaves of $A_n$ slightly we can assume that there is $i \in \ZZ$ so that the leaf $x$ intersects the interior of $\gamma^i_n A_n$ and we can assume that the sequence $\ell_k$ defined above is also in $\gamma^i_n A_n$. Since $A_n$ is closed we get that $y$ is contained in $\gamma^i_n A_n$ contradicting that it is trivially foliated by Lemma \ref{lem-trivialsides}.  


We now analize the case where we have a point $x$ in say
$\Theta_n$ which is an ideal point of $\Theta_m$ (that is, $x$ is a endpoint of $\cI$ but also of $\Theta_m$). The symmetric case exchanging $n$ and $m$ is identical. There are several subccases. Since $\gamma_m$
preserves $\Theta_m$ it follows that $\gamma_m x$ is either mapped to $x$ or to a leaf which is non
separated from $x$. Suppose first that $\gamma_m x$ is distinct
from $x$. Then either $\gamma_m x$ or $\gamma_m^{-1} x$ is
non separated from $x$ and in the positive side of $x$, since $x \in \Theta_n$ this gives the same contradiction as above. 

Suppose now\footnote{Note that from the point of view of obtaining a Reeb crown, this concludes, because $x$ will be the boundary of a Reeb crown in $\cF_1$. However, we continue the argument because it will give us the stronger statement that every non-separated ray accumulates in a Reeb crown.} that $\gamma_m x = x$. Then this projects
to a closed curve $c$ in the projection $\pi(E_1)$ of $E_1$ to $M$. 
We will show that this also leads to a contradiction. Note first that since $\gamma_m x =x$ but $\gamma_n x \neq x$ one has that $\gamma_n$ and $\gamma_m$ cannot belong to the same cyclic group.

Let $\hat c$ be a lift of $c$ to $E_1$ invariant under $\gamma_m$ and intersecting
$\hat \eta_n$ in a unique point $z$. Let $c_0 = \hat c \cap D_n$ a ray of $\hat c$ starting at $z$. By Lemma \ref{lem-trivialsides} 
the region in $D_n$ between $c_0$ and
$\gamma_n(c_0)$ is product foliated by $\wcG_{E_1}$.
But as $k \to \pm \infty$ then $\gamma_n^k(c_0)$ escapes
to the two fixed points of $\gamma_n$, so it follows
that the union of these regions (which is $D_n$)
is a hemisphere of $E_1$
bounded by $\hat \eta_n$ and verifies that every curve of $\wcG_{E_1}$ intersecting $D_n$ must intersect $\hat \eta_n$.

Let $w^-,w^+$ be the ideal points of $\hat \eta_n$ with 
$w^+$ the attracting point of $\gamma_n$.
Let $t^-, t^+$ be the ideal points of $\hat c$ with $t^+$ the
attracting point of $\gamma_m^{-1}$ (in particular
$t^+$ is the ideal point of the ray $c_0$). Since $\gamma_n, \gamma_m$ do not belong to the same cyclic group, the points $w^-,w^+,t^-,t^+$ are pairwise distinct.

Consider $\gamma_m^{-k} \gamma_n(\hat c)$ with $k > 0$.
First $\gamma_n(t^-), \gamma_n(t^+)$ are in
the interval of $S^1(E_1) \setminus \{ t^-, t^+ \}$
which contains $w^+$. They are contained in the
attracting basin of of $t^+$ under $\gamma_m^{-1}$,
in particular for $k$ big enough
$\gamma_m^{-k} \gamma_n(t^-),
\gamma_m^{-k} \gamma_n(t^+)$ are contained in the
interval of $S^1(E_1) \setminus \{ w^-, w^+ \}$ which
contains $t^+$.
These are the ideal points of $\gamma_m^{-k} \gamma_n(\hat c)$.
Since $\hat c$ is fixed and it is a quasigeodesic in $E_1$ it
follows that for $k$ big enough, the whole curve 
$\gamma_m^{-k} \gamma_n(\hat c)$ is contained in 
$\mathring{D_n}$.

This is a contradiction because we proved above that every leaf
intersecting $D_n$ has to intersect $\hat \eta_n$.
This contradiction finally
shows that the assumption that $\Theta_n \not = \Theta_m$ is
impossible and this finishes the proof of the lemma.
\end{proof}

Now, pick $n, m$ large enough. The arguments in the previous proposition allows to show that: 

\begin{lema}\label{lem-cyclicgroupinj} 
The elements $\gamma_n$ and $\gamma_m$ belong to the same cyclic group. 
\end{lema}
\begin{proof}
Let $G = \langle \gamma_n, \gamma_m \rangle$ be the subgroup generated by
$\gamma_n, \gamma_m$. 
Since $\Theta_n = \Theta_m$ it follows that $G$ acts on $\Theta_n$
which is homeomorphic to $\R$.
We first claim that the action is free. Otherwise there is
a non trivial element $\alpha$ of $G$ with a fixed point 
in $\Theta_n$. This implies that $\gamma_n$ and any power
of $\alpha$ cannot belong to the same cyclic group.
Then the last case of the proof of the previous proposition
applies verbatim to produce a contradiction
(apply it to $\gamma_n$ and $\alpha^k$ for $|k|$ sufficiently
big).

By H\"{o}lder's theorem this implies that $G$ is
abelian. If it is not cyclic then it has a $\ZZ^2$ subgroup,
which implies that $\pi(E_1)$ 
is a torus, hence compact.
This is disallowed by hypothesis.
This proves that $G$ is cyclic.
\end{proof}

We can now give the proof of Proposition \ref{lem-converge}.

\begin{proof}[Proof of Proposition \ref{lem-converge} ]
Now we use that $\pi(E_1)$ is a surface, so elements 
are contained in unique maximal cyclic groups.
Let $G$ be a maximal cyclic group containing all the $\{ \gamma_i \}$.
Let $E_G = E_1/_{G}$, an annulus.
Since all $\gamma_n$ are in $G$ and 
$\gamma_n L$ has accumulation points in the interval $J$
of leaves $\{L_t\}_{t \in [-1,1]}$ (with $L=L_0$) it follows that $\gamma_n (L \cap E_1)$ 
converges to a leaf $\ell$ (and maybe other leaves as well).
Since $\gamma_n$ preserves $\Theta_n$ it follows that if 
$\gamma_n \ell \neq \ell$ it must 
be non separated from $\ell$ in $\cG_{E_1}$.

Let $\beta_n$ be segments of bounded
length from $p_n$ to $q_n$ in $L$ so that are transverse to $\cG_L$ in neighborhoods of $p_n,q_n$. We may assume up to 
a further subsequence, that $\gamma_n(\beta_n)$ converges
as $n \to \infty$ 
to a segment $\beta$ in a leaf $L'$ of $\wcF_1$ in the limit
of $\gamma_n(L)$ (as $n \to \infty$) 
which contains the limit of $\gamma_n(p_n)$.
Therefore the non separation between $\gamma_n(r_1)$ and $\gamma_n(r_2)$
in $\gamma_n(L)$ may be pushed to non separation in $L'$.
As we have proved previously, this implies that the product picture in $E_1$ extends beyond
$\ell$. 

This implies that $\gamma_n(\ell) = \ell$, and in particular
this proves that $\pi(r_1)$ is asymptotic to $\pi(\ell)$ which is a closed curve as we wanted to show.
\end{proof}

We can now show: 

\begin{teo}\label{teo.noninjsvm}
Let $\cF_1, \cF_2$ transverse foliations by Gromov hyperbolic leaves and injective developing map. Assume that $\cG= \cF_1 \cap \cF_2$ has the small visual measure property in the foliation $\cF_1$. Then, either it is leafwise quasgeodesic on $\cF_1$ or it has a generalized Reeb surface. 
\end{teo}

\begin{proof}
If $\cG_L$ is Hausdorff for every $L \in \wcF_1$, then, one can apply \cite[Proposition 7.4]{FP-hsdff} to deduce that $\cG$ is leafwise quasigeodesic in $\cF_1$. 

Else consider $L' \in \wcF_1$ and $r_1, r_2$ rays of $\cG_{L'}$
which are in distinct $\wcG$ leaves, but non separated from
each other in $\cG_{L'}$. Hence we are in the conditions of Convention \ref{convention-rays}. 
By Proposition \ref{lem-converge} we know that $\pi(r_1)$ is either closed or
asymptotic to a closed leaf.  

If $\pi(r_1)$ is closed then Proposition \ref{prop-genRS} shows
that $\pi(r_1)$ is a closed boundary component
of a generalized Reeb surface of $\cG$ in a leaf of $\cF_1$.

Suppose then that $r_1$ is asymptotic, but not equal to a leaf $\ell$ of $\wcG$, where $\pi(\ell)$ is
closed. We claim that $\ell$ is non separated from another leaf
of $\wcG$ in its $\wcF_1$ leaf which we call $L$. 

To show this, we use the notations introduced earlier in this section,
in particular $NS(r_1, r_2)$. Let $\gamma$ be a generator of
the stabilizer of $\ell$ moving points in the direction
that $r_1$ is asymptotic to. 
Fix $x$ in $\ell$.
Let $p_n$ be points in $r_1$ which are close to $\gamma^n(x)$
when $n \to \infty$, one can choose $p_n$ so that $d(p_n, \gamma^n(x)) \to 0$ since $\pi(r_1)$ acumulates to $\pi(\ell)$ which is closed.
By Proposition \ref{prop-escapingpoints} 
there are points $q_n \in r_2 \cup NS(r_1,r_2)$ which
are $a_1$ distant in $L'$ from $p_n$. Let $n$ big
enough so that a compact arc $\beta_n$
in $L'$ of length $\leq a_1$ connecting
$p_n, q_n$ can be pushed to an arc of
similar length in $L$. 

Let $E$ be the $\wcF_2$ leaf of $\ell$ (note it is also the $\wcF_2$ leaf containing $r_1$), and let
$E_n = \wcF_2(q_n)$.
Then $E_n$ intersects $L$ and
the intersection is a leaf $\ell'$ of $\wcG$ in $L$.
Exactly as done in the proof of Lemma \ref{lem-trivialsides} 
it follows that $\ell$ and $\ell'$ are non separated in $\cG_L$.
This is the property we wanted to obtain and allows to apply Proposition \ref{prop-genRS} to the projection of the region between $\ell$ and $\ell'$ (and the rest of the non-separated leaves to $\ell$) to $M$ showing the existence of the posited generalized Reeb surface.  
\end{proof}

\subsection{The tricothomy} 

Let us prove the trichotomy in Theorem \ref{teo-injective}. Suppose that $\cG$ does not have generalized Reeb surfaces in either $\cF_1$ or $\cF_2$ (i.e. the first case does not hold). 

Assume that there is small visual measure in say, $\cF_1$, then, we know that the leafspace of $\cG$ in leaves of $\wcF_1$ must be Hausdorff due to Theorem \ref{teo.noninjsvm} and moreover. In addition there must be a closed leaf of $\cG$: this is proved in Theorem 9.1 of \cite{FP-hsdff} $-$ note that in the statement it is assumed that $\wcG$ has Hausforff leaf space, but all that is
used is that $\cG_L$ has Hausdorff leaf space for all $L$ in $\wcF_1$.

We claim that this implies that $\cG$ must also have the small visual measure on $\cF_2$. If that is not the case, then, by  Proposition \ref{prop.svm} applied to $\cF_2$ we get a bubble lamination $\Lambda_2$ in $\cF_2$. Note that since there is a closed leaf of $\cG$ it follows that $\Lambda_2 \neq M$.

We can apply \cite[Lemma 8.10]{FP-hsdff} and \cite[Lemma 8.11]{FP-hsdff}  (see \cite[Lemma 8.12]{FP-hsdff}) to extend the bubble lamination $\Lambda_2$ to the complementary regions, a contradiction with the fact that $\Lambda_2 \neq M$. 
This also needs an explanation as in \cite{FP-hsdff} the assumption is
that $\wcG$ has Hausdorff leaf space: all that is used in these results
of \cite{FP-hsdff}, when analyzing the foliation $\cF_2$, is that
$\cG_L$ has Hausdorff leaf space for all leaves $L$ of the 
{\underline {other}} foliation $\cF_1$. But now we have this Hausdorff 
property as $\cG$ is leafwise quasigeodesic in leaves of
$\cF_1$. 

As a conclusion it follows that $\cG$ has the small visual measure property in $\cF_2$ as well. Applying Theorem \ref{teo.noninjsvm} we conclude the tricothomy.

\subsection{The minimal case}
Assuming minimality, we will show that the small visual measure property always holds as long as $M$ is not a torus bundle. 

\begin{theorem} Let $\cF_1, \cF_2$ be transverse minimal foliations by Gromov hyperbolic leaves
in a closed $3$-manifold, and let $\cG = \cF_1 \cap \cF_2$. 
Suppose that $\cF_1, \cF_2$ have injective developing map, then
one of the following mutually exclusive options occur: 
\begin{itemize}
\item $\cG$ is leafwise quasigeodesic,  
\item $\cG$ has a generalized Reeb surface in one of the foliations, or, 
\item $M$ is a torus bundle and both foliations equivalent to the weak stable foliation of a suspension Anosov flow.
\end{itemize} 
\end{theorem}

Note here that by minimality once one shows leafwise quasigeodesic, automatically has that $\cG$ is the orbit foliation of an Anosov flow due to an argument of Calegari (see \cite[\S 6]{BFP}). 

\begin{proof}
We start from the trichotomy which we have already obtained. We need to show that if $\cG$ is not leafwise quasigeodesic in neither of the foliations nor it has a generalized Reeb surface, then, we are in the third option of the statement. By the trichotomy already proven, we get that there are sublamination $\Lambda_1$ and $\Lambda_2$ of $\cF_1$ and $\cF_2$ made of bubble leaves. By minimality we have that $\Lambda_1=\cF_1$ and $\Lambda_2= \cF_2$. 

In particular there is a distinguished ideal point in leaves
of $\wcF_i$ which is equivariant and continuous. By a result
of Calegari, as explained in \cite[\S 6]{BFP} it follows that $\cF_i$
are Anosov foliations. Denote by $\varphi_i$ the Anosov flow associated to $\cF_i$. We orient $\varphi_i$
towards the distinguished ideal point in the leaves
$\wcF_i$. 

The proof is by contradiction, assuming that $\varphi_i$ is not a suspension Anosov flow. We separate the contradiction depending on the class of manifold $M$:

\vskip .1in
\noindent
{\bf {Case 1}} $-$ $M$ is Seifert fibered.

Since the foliations are minimal, then they must be $\RR$-covered. As we are in the injective developing map case the leaf space of $\wcG$ in any leaf $L$ of $\wcF_1$ or $\wcF_2$ is
Hausdorff (see Lemma \ref{lem.same}). 

With this property we can apply directly the main result of \cite{FP-hsdff} to get 
leafwise quasgeodesic property of $\wcG$, a contradiction.  

\vskip .1in
\noindent
{\bf {Case 2}} $-$ $M$ is hyperbolic.


Let $L \in \wt{\cF_1}$ be a leaf invariant under $\gamma \in \pi_1(M) 
\setminus \{ \mathrm{id} \}$ which must then fix $\xi(L)$. Let $c \subset L$ be a geodesic with one endpoint in $\xi(L)$ and the other in the other fixed point of $\gamma$ in $S^1(L)$ (up to a bounded
distortion within $L$, then $c$ corresponds to the lift of a periodic orbit of $\varphi_1$). This geodesic  $c$ projects to a quasigeodesic in $M$ (because it is a closed curve) whose endpoints in $\partial \wt M \cong \partial \mathbb{H}^3$ are $\gamma^-$ and $\gamma^+$. Consider $\gamma^+$ (up to taking the inverse of $\gamma$) to correspond to the ray of $c$ converging to $\xi(L)$. 

Since $L$ is a bubble leaf, it follows that every curve $\ell \in \cG_L$ has both rays converging to $\xi(L)$ in $L \cup S^1(L)$ and since $\cF_1$ has the continuous extension property (see \cite{Fen1})  it follows that both rays of $\ell$ converge to $\gamma^+$ in $\partial \mathbb{H}^3$. 

Consider an arbitrary leaf $F \in \wt{\cF_2}$ intersecting $L$. By assumption, it is a bubble leaf with endpoint $\xi(F) \in S^1(F)$ and $F\cap L$ verifies that both rays converge to $\gamma^+$ in $\mathbb{H}^3$. By the continuous extension property applied to $\cF_2$, it follows that $\xi(F)$ corresponds to $\gamma^+$.

Since periodic orbits of the flow $\varphi_2$ associated to $\cF_2$ are dense, we find several deck transformations $\eta, \eta' \in \pi_1(M)$ not in the same cyclic group and which must therefore have fixed point in $\gamma^+$ a contradiction.

\vskip .1in
\noindent
{\bf {Case 3}} $-$ $M$ has non trivial JSJ decomposition.

Since we have assumed that $\varphi_i$ is not a suspension Anosov flow, it follows
that a torus in the JSJ decomposition cannot cut $M$ into a product. 

Let $T$ be a torus in the
JSJ decomposition, homotope it to be quasi-transverse to $\varphi_1$ (see \cite[\S 6.3]{BarthelmeMann}).
 Since $\cF_1$ is minimal, $\varphi_1$ is transitive and thus there
is a periodic orbit of $\varphi_1$ intersecting $T$ transversely. 

Let $L$ be the $\wcF_1$ leaf of a lift of this periodic orbit
of $\varphi_1$, and let $c$ be the lift of the periodic
orbit in $L$, which is a curve
invariant under $\gamma \in \pi_1(M) \setminus \{ \mathrm{id} \}$
a deck transformation leaving $L$ invariant.

The lifts of $T$ may not be disjoint $-$ this is because of the quasi-transverse property which means $T$ may self intersect tangentially along some 
periodic orbits contained in $T$. But each such lift to $\mt$ 
is a properly embedded
plane, which is quasi-isometrically embedded in $\mt$.
This follows from \cite[Theorem 1.1]{KL} (see also \cite[Section 3.1]{Ng}) because  $M$ is not a torus bundle.  

Let $\wwt_i$ be the lifts of $T$ 
which are successively  intersected by $c$,
and let $t_i := \wwt_i \cap L$.

\begin{claim}
The $t_i$ are quasigeodesics in $L$ (with constants independent of $i$) whose endpoints in $S^1(L)$ are different from $\xi(L)$. Moreover, in $L \cup S^1(L)$ the closure of $t_i$ and converge to $\xi(L)$.
\end{claim}

\begin{proof}
The intersection $\cF_1 \cap T$ is
a one dimensional foliation in the torus $T$.

If this foliation does not have Reeb annuli, then, the foliation is by quasigeodesics when lifted to any lift $\widetilde{T}$ of $T$ and since $\widetilde{T}$ is quasigeodesically embedded in $\wt M$, it follows that every intersection with $L$ is a quasigeodesic (with uniform constants) and we deduce that $\ell_i$ are all quasigeodesics.

Suppose that $\cF_1 \cap T$ has Reeb annuli, and
suppose $c$ intersects a lift of $T$ in one
of these annuli, let $t$ be the intersection of $L$
with this lift. Then $t$ is a curve in $L$ that intersects
$c$ in a single point. 
The rays of $\pi(t)$ in the torus $T$ are asymptotic
to closed curves in $T$, hence the rays of $t$
are quasi-isometrically embedded in $\widetilde T$, hence
also quasi-isometrically embedded in $\widetilde M$,
and in $L$. In particular each such this ray
is a bounded distance from the corresponding
geodesic in $L$. Suppose that one of these rays
has an ideal point $q$ in $S^1(F)$, which is fixed by $\gamma$
(recall that $\gamma(c) = c$).
Then the same happens for all $\gamma^n(t)$ in $L$.
Choose $\gamma$ so that $q$ is the repelling fixed point.
Then $\gamma^n(t)$ for $n \geq 0$, all intersect
a fixed bounded neighborhood of a geodesic in $L$
with ideal points the fixed points of $\gamma$ in $S^1(L)$.
In particular they accumulate in $L$, and hence also
in $\widetilde M$. But these curves $\gamma^n(t)$
are in distinct lifts $\gamma^n(\widetilde T)$ of $T$,
and they cannot accumulate in $\widetilde M$, because
$T$ is compact. This is a contradiction.  

It follows that the ideal points of $t$ is $S^1(L)$ are
distinct from each other and so $t$ is a quasigeodesic
in $L$. The arcs $\gamma^n(t)$ all have pairwise disjoint
sets of ideal points, so they form a neighborhood basis
of $\xi(L)$ which is one of the fixed points of $\gamma$.
This proves the claim.
\end{proof}

The foliation $\wcG$ in $L$ is a bubble foliation, and all leaves
of $\wcG$ in $L$ limit only in $\xi(L)$.
Since $t_i$ are all quasigeodesics with same constants and endpoints are different from $\xi(L)$, It follows that a ray of an arbitrary leaf $\ell$ of $\wcG$ in $L$ 
will intersect
the $t_i$ for $i$ big enough and these rays  escape to $\xi(L)$. 

Now fix an open interval of $\wcF_2$ leaves intersecting $L$, parametrize 
this interval as $F_t, \ 0 < s < 1$. Let $a_s = F_s \cap L$.
For simplicity assume Candel metric on 
the foliation $\cF_2$. Now homotope $T$  to be quasi-transverse
to $\varphi_2$.
Let $T'$ be the resulting torus, then lifts $\wwt_i$ of $T$ correspond to lifts $\wwt'_i$ of $T'$ at bounded distance from $\wwt_i$. Denote by $t'_i = \wwt'_i \cap L$ which are bounded  distance from the corresponding $t_i$.

For each $s$, then in $E_s$  the $\wcG$ leaves $a_s$ intersect every $t'_i$ for $i$ large. In particular the
flow lines of $\widetilde \varphi_2$ in $E_s$ also do that.
Fix a flow line in $E_0$ which intersects all $\wwt'_i$ (remove
some initial $\wwt'_i$ if necessary).

The collection of the flow lines which intersects 
all $\wwt'_i$ transversally is weak 
stable $\wt \varphi_2$ saturated: if an orbit $v$ intersects all $\wwt'_i$ and $v'$ is asymptotic to $v$ then it must eventually be in the same connected component of $\mt \setminus \wwt'_i$ and thus intersect it. (Note that if an orbit of $\widetilde \varphi_2$ intersects a lift of a quasi-transverse torus $T'$ transversally, then, it cannot intersect this same lift of $T'$ again\footnote{The only orbits not intersecting $\wwt'_i$ transversally are those for which the orbit is contained in $\wwt'_i$, moreover, if a segment of orbit intersects $\wwt'_i$ in two points, then, it is contained in $\wwt'_i$.}.)

Suppose that two points not in the same weak stable leaf go through
all tori lifts $\wwt'_i$. The weak stable and weak unstable leaves
of $\widetilde \varphi_2$
separate $\mt$, so it follows that two points on same 
weak unstable leaf intersect all $\wwt'_i$. The set of points in the same
weak unstable leaf that go through two lifts of tori is connected, so we get a non trivial
segment of orbits in the same weak unstable leaf, which 
goes through all lifts $\wwt'_i$.

	But then the forward saturation of this weak unstable strip is dense when projected to $M$. For any
periodic orbit of $\varphi_2$, it is forward asymptotic to an orbit in this projection, because $\varphi_2$ is transitive. Therefore any
periodic orbit of $\varphi_2$ has a lift intersecting all 
$\wwt'_i$ (recall that the sets intersecting all $\wwt'_i$ are 
are stable saturated).
	This forces $T'$ to intersect all orbits of the flow
$\varphi_2$.  This is impossible as the flow $\varphi_2$ is not
a suspension. This completes this case and ends the proof. 
\end{proof}

\subsection{Hyperbolic and Seifert manifolds} 

When $M$ is Seifert, then, being by Gromov hyperbolic leaves implies that every leaf is horizontal. In particular both foliations
are $\RR$-covered: so we obtain either a Reeb annulus
or leafwise quasigeodesic behavior.
Thus, the third option in the new version of Theorem B cannot happen. 

Now we discuss further on the third option which in particular imply that it cannot happen either when $M$ is hyperbolic. 

Let us assume that there are sublaminations $\Lambda_1, \Lambda_2$ of $\cF_1$ and $\cF_2$ of bubble leaves. We consider the sets $\Lambda_i$ to be actually the set of \emph{all} bubble leaves of $\cF_i$, this is closed saturated set, so a sublamination.  Without loss of generality we assume everything is orientable. 

\begin{lemma}
In the laminations $\Lambda_i$ every leaf is a cylinder or a plane, moreover, there is some small $\eps>0$ so that the boundary of $B_\eps(\Lambda_i)$ is a union of tori. 
\end{lemma}

\begin{proof}
Consider $\widetilde{\Lambda_i}$ to be the lift of $\Lambda_i$ to $\mt$. Since for each $L \in \wt{\Lambda_i}$ there is a distinguished point $\xi(L)$ in $S^1(L)$ we deduce that the stabilizer of $L$ is cyclic. 

Now, take an octopus decomposition of $M \setminus \Lambda_i$ and note that boundary leaves must be annuli (else, one can push the bubble leaves along, as in \cite[Lemma 8.10]{FP-hsdff}). In addition there is an ordering around, thus giving rise to a torus boundary. 
\end{proof}

Now, let us show that these tori boundary cannot bound a solid torus: 

\begin{lemma}
Each boundary torus of $B_\eps(\Lambda_i)$ is an incompressible torus.
\end{lemma}

\begin{proof}
Assume by way of contradiction that one such torus is compressible,
and hence bounds a solid torus.

There are two cases, either the complementary region $U$ has two boundary components, in which case we can argue exactly as in  \cite[Lemma 8.12]{FP-hsdff} to show that actually the leaves in the complementary region are also bubble leaves and so belong to $\Lambda_i$.

If $U$ has more boundary leaves, we can consider a lift $\tilde U$ of $U$ and since it is a solid torus, the boundary leaves of $\wt{\Lambda_i}$ will be in one to one correspondence with the boundary leaves of $U$, and are a finite sequence $L_1, \ldots, L_k$ circularly ordered so that $L_i$ and $L_{i+1}$ share a half space corresponding to an \emph{arm} of the octopus decomposition of the complementary region $U$ (see \cite[\S 8.3]{FP-hsdff}). All leaves $L_i$ are fixed by some deck transformation $\gamma$ which corresponds to the deck transformation fixing $\tilde U$. 

Now, consider a leaf $\ell_1 \in \cG_{L_1}$ which is in the axis of $\gamma$ non trivial deck transformation, and let $F \in \wt{\cF_2}$ by so that $\ell = F \cap L_1$. Arguing exactly as in \cite[Lemma 8.11]{FP-hsdff} we can show inductively in $i$ that $L_i \cap F \neq \emptyset$ and call $\ell_i = F\cap L_i$.  By injectivity after turning around we return to the same leaf, so we obtain that $F \cap \tilde U$ is an open disk whose boundary consists of the leaves $\ell_i$. Since there are more than $2$ arms in the octopus decomposition, an index argument makes this incompatible with filling the region $\tilde U$ with a transverse foliation. 
\end{proof}

With this we readily get the result for hyperbolic 3-manifolds since in this case we cannot have incompressible tori and so $\Lambda_i = \cF_i$. After this is done, we can argue as in the previous section to get a contradiction as we have obtained that $\cF_i$ is actually a blow up of an Anosov foliation. This completes the proof of Theorem \ref{teo-injective}.

\section{Some general position constructions}\label{s.lemas}
In this section we compile some useful general results that will be needed in \S \ref{s.B2}. In \S~\ref{ss.nonsepinter} we define a notion of non-separated intersection between leaves of $\wcF_1$ and $\wcF_2$ that will be crucial in the proof. 

\subsection{Curves in general position connecting non
separated leaves} 

Here we show the following useful result that we will use later to produce some relevant regions in three dimensional space (the reader should compare with the construction of the set $U_{r_1,r_2}$, see Figure \ref{f.ur1r2}). 

\begin{lema}\label{lema-curvealpha}
Let $\cG$ be a one dimensional foliation 
of a plane $L$ and let $\ell_1, \ell_2 \in \cG$ be two non-separated leaves. Then, given $x \in \ell_1, y \in \ell_2$ there is a $C^1$-arc $\alpha:[0,1] \to L$ so that $\alpha(0)=x$, $\alpha(1)=y$ and so that up to changing orientation of $\cG$ we have: 
\begin{itemize}
\item $\alpha$ is tangent to $\cG$ exactly at $\alpha(1/2)$ where it intersects a leaf $\ell'$ and it is its unique intersection point with $\ell'$,
\item for every $t \in (0,1/2)$ the curve $\alpha$ is positively transverse to $\cG$ at $\alpha(t)$ and intersects each leaf in a unique point in that interval,
\item for every $t \in (1/2, 1)$ the curve $\alpha$ is negatively transverse to $\cG$ at $\alpha(t)$ and intersects each leaf in a unique point in that interval. Moreover, $\alpha(t)$ intersects the same leaf as $\alpha(1-t)$. 
\end{itemize}
Let $U$ be the the union of the compact segments in leaves of $\cG$ which intersect $\alpha$ in two points, and have both endpoints in $\alpha$.  The interior  of $U$ is homeomorphic to a disk, and the closure of $U$ in $L$ consists of $\alpha$, the non-separated rays of $\ell_1$ and $\ell_2$ and possibly some other non-separated leaves of $\cG$ from $\ell_1, \ell_2$. 
\end{lema} 

\begin{proof}
Using that $\ell_1$ and $\ell_2$ are non-separated we have that close to $x,y$ we can consider transversals $\tau_1,\tau_2: [0,1/2] \to L$ to $\cG$ so that $\tau_1(0)=x, \tau_2(0)=y$ and such that $\tau_1(t)$ belongs to the same leaf of $\cG$ as $\tau_2(t)$. Call $c_t$ the arc in a leaf of $\cG$ from $\tau_1(t)$ to $\tau_2(t)$ and $\hat c_t$ the same arc with the endpoints removed. 

The set $\bigcup_{t \in [1/4,1/2]} c_t$ is a foliated neighborhood and one can modify the curve $\tau_1|_{[1/4,1/2]}$ to make it transverse to $\cG$ and to limit as $t \to 1/2$ to the midpoint of $c_{1/2}$. Doing the same with $\tau_2$ and concatenating those curves one obtains the desired $\alpha.$ 

Finally, the fact that $\bigcup_{t \in (0,1/2]} \hat c_t$ is an open disk is direct since it is homeomorphic to the product of open intervals. By construction, the closure contains the image of $\alpha$ as well as the non-separated rays of $\ell_1$ and $\ell_2$. Since the rest of the closure needs to be $\cG$ saturated, it follows by definition that the leaves that belong to these closure must be non-separated from $\ell_1$ and $\ell_2$. 
\end{proof}

\subsection{A family of transverse disks}\label{ss.transversals}

Let $\cF_1, \cF_2$ be two transverse foliations in $M$ and let 
$\cG := \cF_1 \cap \cF_2$ be the intersected foliation. We will consider a fixed Riemannian metric in $M$ once and for all. 

To get some uniform estimates, we will consider a family of disks $\{O_x\}_{x\in M}$ transverse to $\cG$ of uniform size. We can consider these disks to be smooth and to vary continuously 
with $x$ by considering the projection by the exponential map of the orthogonal space (with respect to the chosen Riemannian metric) to the tangent space to the $\cG$ curve through $x$.

We can assume without loss of generality that the intersection of $O_x$ and both foliations $\cF_1$ and $\cF_2$ gives a trivial foliation by curves that intersect in a unique point and have a local product structure (i.e. the disks are in fact open rectangles). This allows one to define quadrants in $O_x$. We will consider the size of the sets $O_x$ as $\eps_0$ (related to local product structure (cf. \S \ref{ss.transverse}). 

\begin{convention}[Sign convention on quadrants]\label{conv-sign}
If one has a local tranverse orientation of $\cF_1$ and $\cF_2$ this allows to name these quadrants as $(+,+)$, $(+,-)$, $(-,+)$ and $(-,-)$ according to these orientations. 
Specifically the first entry refers to the positive or negative
transverse orientation to $\cF_1$ at $x$ and similarly for
the second entry.
\end{convention}

We will consider, for $p \in \mt$ the disk $\tilde O_p$ which lifts $O_x$ so that $p$ projects to $x$ in $M$. If $p,q \in \ell \in \wcG$ we can define the map $\Psi_{p,q} : A_{p,q} \to \tilde O_q$ where $A_{p,q} \subset \tilde O_p$ is the domain of the $\wcG$-holonomy from $p$ to $q$.

\subsection{Doubly non-separated intersection and disks in good position}\label{ss.nonsepinter} 

Let $\cF_1$ and $\cF_2$ be two transverse foliations and $\cG= \cF_1 \cap \cF_2$. Let $L \in \wcF_1$ and $E \in \wcF_2$ be leaves intersecting in at least two\footnote{By transversality the intersection is at most countably many leaves.} leaves $\{\ell_i\}_{i \in I}$ so that $\ell_1, \ell_2$ are non separated in \emph{both} $\cG_L$ and in $\cG_E$ (recall that $\cG_L$ and $\cG_E$ are the respective restrictions of $\wcG$ to $L$ and $E$). 
We stress that in general if $L \cap E$ contains two components $\ell_1, \ell_2$, we know that the leaf spaces of $\cG_L$ and $\cG_E$ cannot be Hausdorff, but it could be that $\ell_1$ and $\ell_2$ are separated in both, or non-separated in only one of those leaf spaces. 

When such a configuration holds, we say that $L$ and $E$ have a \emph{doubly non-separated intersection}.

\begin{lema}\label{lem-nonseprays}
Let $\ell_1, \ell_2$ be doubly non-separated leaves in $L \cap E$. Then, there are rays $r_1, r_2$ of $\ell_1$ and $\ell_2$ so that these are the non-separated rays of $\ell_1$ and $\ell_2$ in both $\cG_L$ and $\cG_E$. (Recall Figure \ref{f.nsr}.)
\end{lema}

\begin{figure}[htbp]
\begin{center}
\includegraphics[scale=0.56]{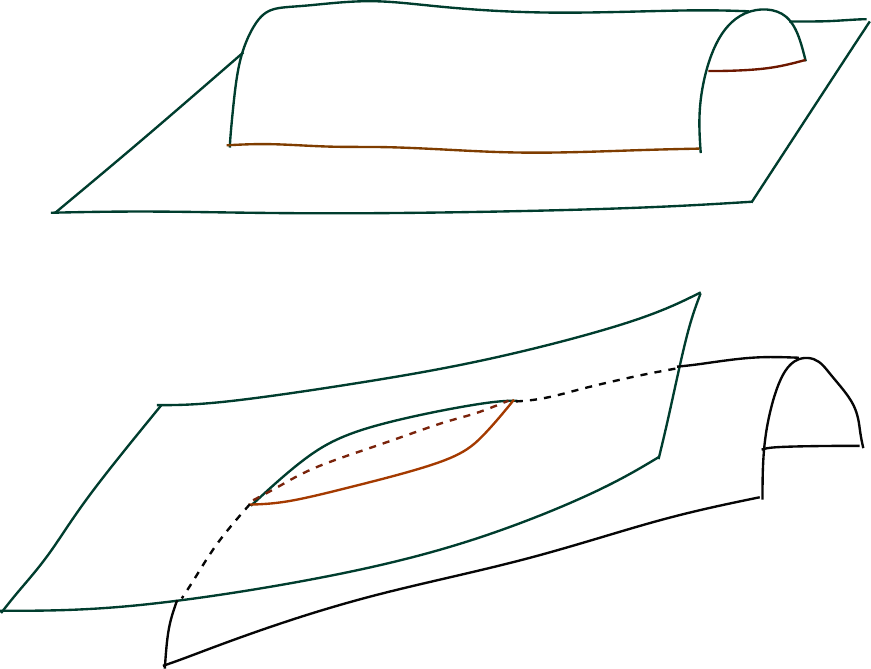}
\begin{picture}(0,0)
\put(-93,167){$E$}
\put(-93,131){$\ell_1$}
\put(-202,129){$L$}
\put(-33,152){$\ell_2$}
\put(-198,48){$L'$}
\put(-45,70){$E'$}
\end{picture}
\end{center}
\vspace{-0.5cm}
\caption{{\small Above two leaves $L \in \wcF_1$ and $E \in \wcF_2$ intersecting in a pair of doubly non-sepatated leaves $\ell_1$ and $\ell_2$. Below is depicted what would happen if closeby leaves $L'$ and $E'$ escape in different directions and produce a circle.}\label{f.escape}}
\end{figure}
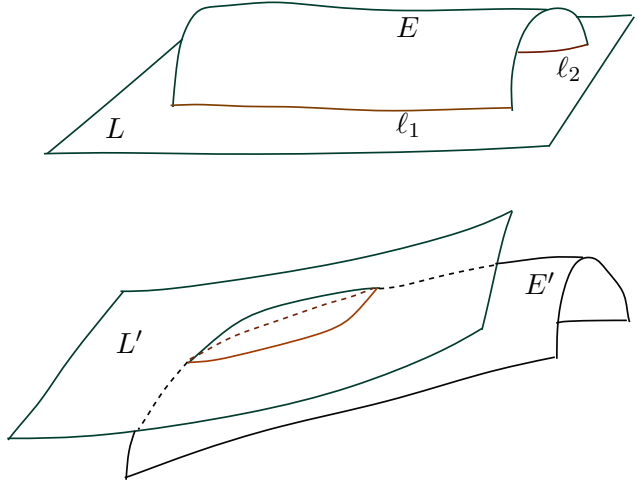

\begin{proof}
Let $r_1$ and $r_2$ be the non-separated rays of $\ell_1$ and $\ell_2$ in $L$. Since $\wcG$ is oriented, and being non-separated rays they must be oriented in oposite direction (one towards the starting point of the ray and the other against), it follows that in $E$ either $r_1$ and $r_2$ are non-separated, or it is the complementary rays $r'_i= \ell_i \setminus r_i$ which are non-separated. Let us see that the latter case cannot happen by contradiction. If it were the case, one would get very close leaves $E'$ and $L'$ to $E$ and $L$ respectively so that the intersection of $E'$ with $L$ has an arc joining points very close to the endpoints of $r_1$ and $r_2$ and oriented in the direction of the non-separation. Similarly with $L'$ and $E$. By continuity, this implies that $E' \cap L'$ contains a closed curve, which is a contradiction since both are planar leaves and $\wcG$ has no singularities.  See Figure \ref{f.escape}.
\end{proof}

Our goal is to pick a curve $\alpha_1$ joining $\ell_1$ and $\ell_2$ in $L$ and a curve $\alpha_2$ joining the same points in $E$ and produce a region between $L$ and $E$ that goes in the direction of the non-separated rays $r_1,r_2$ of $\ell_1$ and $\ell_2$ bounded by at disk whose boundary is the union of $\alpha_1$ and $\alpha_2$. We can give a transverse orientation to $\wcF_1$ and $\wcF_2$ so that for the rays $r_1, r_2$ it holds that the positive direction is the one that gives the non-separation between the rays in both foliations. This defines quadrants for each local transversal $\tilde O_p$ as defined in \S~\ref{ss.transversals}.     

We first define a disk bounded by curves in $L$ and $E$ that will allow us to define in the next subsection a region in between the leaves that will be technically important in the proof of Theorem
\ref{teo-noninjective}. 

\begin{prop}\label{prop-nonsepintersection}
Let $L \in \wcF_1$ and $E \in \wcF_2$ have doubly non-separated intersection. Let $\ell_1$ be a connected component of $L \cap E$ and $r_1$ be a ray which is non-separated in both $\cG_L$ and $\cG_E$ from other components (could
be a single one) of $L \cap E$. For every $p \in r_1$ and $q \in \ell \in L \cap E$ so that the ray $r_1$ is non-separated from $\ell$ in both $\cG_L$ and
$\cG_E$,  there is a disk $D_{p,q}$ whose boundary consists of two curves $\alpha_1 \subset L$ and $\alpha_2 \subset E$ connecting $p$ and $q$ verifying the following properties: 
\begin{enumerate}
\item The curve $\alpha_1$ intersects each leaf of $\cG_L$ in at most two points and is transverse to $\cG_L$ except at one point. Similarly with $\alpha_2$ and $\cG_E$. (Recall Lemma \ref{lema-curvealpha}.)
\item The interior disk $D_{p,q}$ is transverse to both $\wcF_1$ and $\wcF_2$.
\item The intersection of $D_{p,q}$ with a leaf of $\wcF_2$ intersecting its interior is exactly a curve going from $\alpha_1$ to itself. Symmetrically with $D_{p,q}$ and $\wcF_1$ and the curve $\alpha_2$. 
\item With appropriate transverse orientations we have that $D_{p,q} \cap \tilde O_p$ and $D_{p,q} \cap \tilde O_q$ contain the full corresponding $(+,+)$-quadrants.  
\end{enumerate}
\end{prop}

\begin{proof}
We choose transverse orientations to $\wcF_1, \wcF_2$ so that
the double non-separation at $p$ is in the $(+,+)$ quadrant.
We will construct such a disk which intersects $\tilde O_p$ and $\tilde O_q$ in an open set of the quadrant $(+,+)$ and then pushing along the flow of $\wcG$ one can deform the disk so that it verifies the last property. 

We start by choosing curves $\alpha_1$ in $L$, $\alpha_2$ in $E$ both
from $p$ to $q$ and satisfying the properties of Lemma 
\ref{lema-curvealpha}.

To do this, we start with a smooth vector field $X$ transverse to $\wcF_1$ in a neighborhood of $\alpha_1$ and similarly a smooth vector field $Y$ transverse to $\wcF_2$ in a neighborhood of $\alpha_2$. (One can choose the vector fields to be tangent to $\tilde O_p$ and $\tilde O_q$ 
 where it makes sense and also so that $X$ is tangent to $\alpha_2$ and $Y$ to $\alpha_1$ close to $p$ and $q$.) This allows to construct small \emph{walls} $W_1, W_2$ transverse to $\wcF_1$ and $\wcF_2$ respectively along $\alpha_1$ and $\alpha_2$. 

We further require that  $W_i$ is transverse
to both $\wcF_1$ and $\wcF_2$. To do that we adjust $W_i$ by an 
isotopy, we explain this for $W_1$. Let $z$ be the unique point
of $\alpha_1$ where $\alpha_1$ is not transverse to $\wcF_2$.
Consider the leaf $\hat E$ of $\wcF_2$ through $z$. 
Consider a small compact arc $\beta$ in $\wcG(z)$ with $z$ in the 
interior. This arc is contained in $L \cap \hat E$ and has endpoints
$z_1, z_2$. Now push this arc slightly along $\hat E$ in the positive
direction transverse to $\wcF_1$ to an arc $\beta'$ still contained
in $\hat E$ and now with endpoints $z'_1, z'_2$. Now connect the
very near endpoints $z_1, z'_1$ and $z_2, z'_2$ respectively
by very short arcs $\delta_1, \delta_2$ in $\hat E$. The union

$$\beta \cup \beta' \cup \delta_1 \cup \delta_2$$

\noindent
is a closed curve in $\hat E$ which bounds a disk in $\hat E$,
which we denote by $\hat D$. 

Now consider arcs of $\wcG$ in $L$ near $z$ and with both
endpoints in $\alpha_1$. Parametrize these as $\gamma^*_t$ where
$\gamma^*_t$ converges to the point $z$ when $t \to 0$.
Extend these arcs in $L$ crossing $\alpha_1$, so that the
extended arcs, denoted by $\gamma_t$ 
now limit to $\beta$ when $t \to 0$. We can do the same procedure
as above for each $t > 0$: perturb $\gamma_t$ slightly along
their $\wcF_2$ leaves, and eventually produce disks $\hat D_t$
which have boundary being the union of an arc in $L$ (the $\gamma_t$),
and an arc in the $\wcF_2$ leaf of $\gamma_t$. These
disks converge to $\hat D$ when $t \to 0$.

Consider this family of disks $\hat D_t, t \geq 0$. They describe
the foliation $\wcF_2$ near $z$ on the positive side
of $L$, and the negative side of $\hat E$. Using this we can then isotope 
$W_1$ so that it is also transverse to $\wcF_2$
near $z$, except at $z$. 
This is because we described the foliation $\wcF_2$ near $z$
on the positive side of $L$ and the negative side of  $\hat E$
(positive and negative refer to the transverse orientations).
In particular the $\wcF_2$ leaves will intersect the adjusted
$W_1$ in either the point $z$ if the $\wcF_2$ leaf is $\hat E$,
and in an arc contained in $\hat D_t$ which has both endpoints
in $\alpha_1$ and otherwise it is contained in the positive
side of $L$.

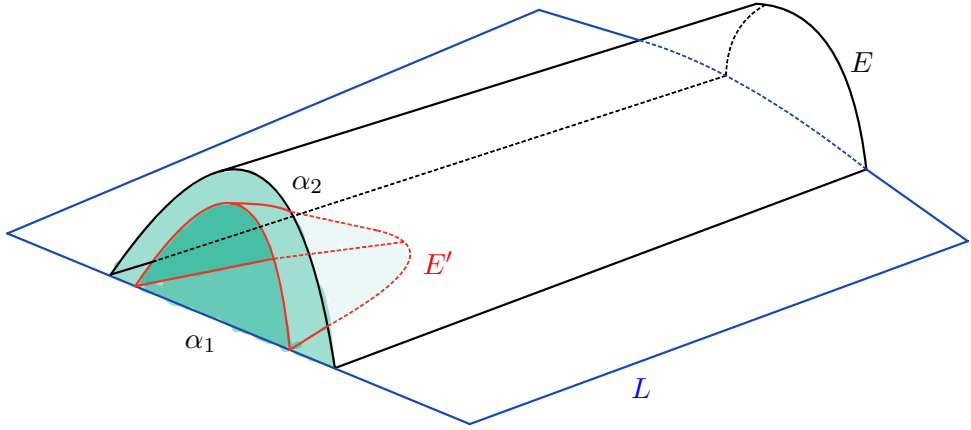
\begin{figure}[htbp]
\begin{center}
\includegraphics[scale=0.69]{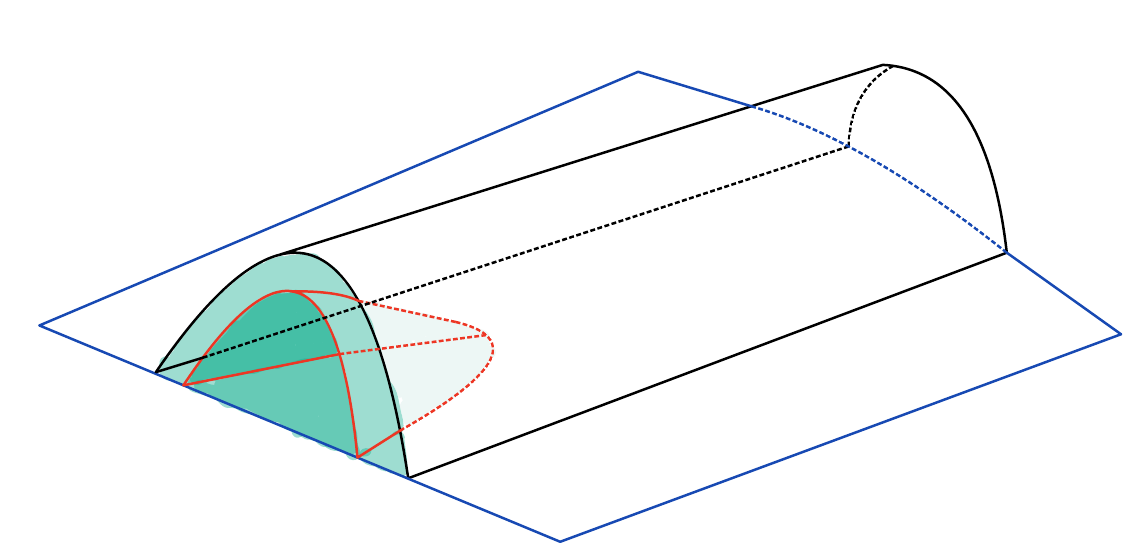}
\begin{picture}(0,0)
\put(-50,136){$E$}
\put(-210,60){{\color{red}$E'$}}
\put(-132,13){{\color{blue}$L$}}
\put(-299,32){$\alpha_1$}
\put(-259,92){$\alpha_2$}
\end{picture}
\end{center}
\vspace{-0.5cm}
\caption{{\small Depiction of some objects in the construction of $D_{p,q}$.}\label{f.Dp}}
\end{figure}

Then we can adjust the rest of $W_1$
so that $W_1$ will be transverse $\wcF_1$, and transverse
to $\wcF_2$ except at two points in the boundary (one of which
is $z$).
After this is done for both $W_1$ and $W_2$ we can then adjust
one of them so that $W_1, W_2$ coincide in small neighborhoods of
$p$ and $q$.

Now we explain how to ``fill in" $W_1 \cup W_2$ to a disk.
Choosing leaves $L'$ and $E'$ above $L$ and $E$ respectively which intersect the walls $W_1$ and $W_2$ in curves $\alpha'_1$ and $\alpha_2'$ close to $\alpha_1$ and $\alpha_2$ (above means in the positive side
defined by the transversal orientation). Putting these together we get a small annular band $B$ which verifies the following properties:

\begin{itemize}
\item the \emph{outside} boundary of $B$ is $\alpha_1  \cup \alpha_2$ and the \emph{inside} boundary is $\alpha'_1 \cup \alpha_2'$. (The upper part of $B$ is depicted in Figure \ref{f.Dp}.) 
\item the interior of $B$ is transverse to both $\wcF_1$ and $\wcF_2$. In addition for every $L'' \in \wcF_1$ between $L$ and $L'$ and $E'' \in \wcF_2$ between $E$ and $E'$ it follows that each of $L'' \cap B$ and $E'' \cap B$ 
is a compact arc (hence connected), and these
two curves intersect in exactly two points.
\item for every $L'' \in \wcF_1$ intersecting $B$ but not in between $L$ and $L'$ one has that the intersection between $L'' \cap B$ and a leaf $E''$ intersecting $B$ can be either empty, one point, or two points, and symmetrically, 
\item for every $E'' \in \wcF_2$ intersecting $B$ but not in between $E$ and $E'$ one has that the intersection between $E'' \cap B$ and a leaf $L''$ intersecting $B$ can be either empty, a single point or two points. 
\end{itemize}

Note that $E' \cap L'$ contains a compact arc $\eta$ which intersects $\partial B$ exactly in $\partial \eta$. Then, we can define disks$D_{L'} \subset L'$ bounded by $\eta \cup \alpha_1'$ and $D_{E'} \subset E'$ bounded by $\eta \cup \alpha_2'$ whose union is a disk $\hat D$. 

Note that $B \cup \hat D$ verifies the desired conditions except at $\hat D$ where it is tangent respectively to $L'$ (along $D_{L'}$),
 and to $E'$ (along $D_{E'}$)  but it still verifies that for every other leaf the third condition is verified. So, it is enough to modify slightly the disk in a neighborhood of $\hat D$ in order to make it transverse to the foliations. For this, one can consider a deformation so that we get the first three properties. To deform $D_{E'}$ consider the part of its boundary that is contained in $B$: it is the arc $\alpha'_2$. Now in the 
part in $B$ very near $D_{E'}$, push it to be becoming
more and more tangent to $\wcF_2$, tagging along $D_{E'}$ so that
at the end it becomes the arc $\eta$. This can be made transverse to
both foliations: transverse to $\wcF_1$ because it is becoming more 
and more tangent to $D_{E'}$ which is a disk in a leaf of $\wcF_2$
which is transverse to $\wcF_1$. We also use that $B$ is transverse
to $\wcF_1$ in its interior. The transversality to $\wcF_2$ is
more immediate. Then do the same to perturb $D_{L'}$. 

The last condition was discussed at the beginning. 
\end{proof}

\begin{figure}[htbp]
\begin{center}
\includegraphics[scale=0.69]{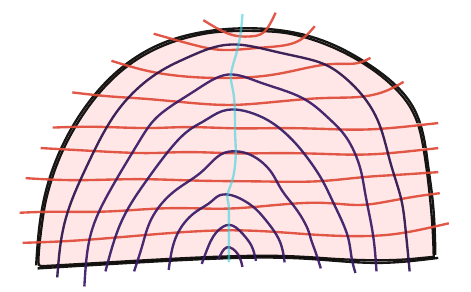}
\begin{picture}(0,0)
\end{picture}
\end{center}
\vspace{-0.5cm}
\caption{{\small There is a line of tangencies between the foliations induced by $\wcF_1$ and $\wcF_2$ in $D_{p,q}$.}\label{f.Disk}}
\end{figure}

\begin{remark} The proposition produces a disk $D_{p,q}$ 
so that its interior is transverse to both
$\wcF_1$, $\wcF_2$, which then induce one dimensional foliations
$\wcF^D_1, \wcF^D_2$ in the interior  of $D$. Notice
that $\wcF^D_1, \wcF^D_2$ are not transverse
to each other in the interior of $D$, even though $\wcF_1, \wcF_2$
are transverse to each other in $\mt$. See Figure \ref{f.Disk}. 
\end{remark}

\subsection{A region between a doubly non-separated intersection} 

Here, we let $L \in \wcF_1$ and $E \in \wcF_2$ be two leaves with doubly non-separated intersection and we consider a point $p$ in a doubly non-separated ray $r_1$ of a component of the intersection between $L$ and $E$. We also consider a point $q$ in another leaf of $\wcG$ contained  in
$L \cap E$ and doubly non separated from $r_1$. We let $D_{p,q}$ be the disk constructed in Proposition \ref{prop-nonsepintersection}. Here, we shall construct an open region $V_{p,q}$ in $\mt$ which is trivially foliated by disks of leaves of $\wcF_1$ and $\wcF_2$ and bounded by $D_{p,q}$ and parts of $L$ and $E$ as well as possibly other leaves of these foliations. This region will be crucial for the constructions in the proof of Theorem \ref{teo-noninjective} but we will try to condense in a statement all the needed properties of this region so that its construction can be read independently of the proof of Theorem \ref{teo-noninjective}.

\begin{prop}\label{prop-tubeconstruction}
Let $L \in \wcF_1$ and $E \in \wcF_2$ be two leaves with doubly non-separated intersection. Let $p,q$ be points in different non-separated leaves of $L \cap E$ in both $\cG_L$ and $\cG_E$ defining non-separated rays $r_1$ and $r_2$. Let $D_{p,q}$ be a disk as defined in Proposition \ref{prop-nonsepintersection}. Then, there is a well defined open set $V_{p,q}$ homeomorphic to a ball so that it verifies: 
\begin{enumerate}
\item If $L' \in \wcF_1$ intersects $V_{p,q}$ then $L' \cap V_{p,q}$ is an open disk with compact closure in $L'$ whose boundary consists of an arc in $E$ and an arc in $D_{p,q}$. Similarly, if $E' \in \wcF_2$ intersects $V_{p,q}$ then $E' \cap V_{p,q}$ is an open disk whose closure is a compact disk whose boundary is an arc in $L$ and an arc in $D_{p,q}$. 
\item The boundary of $V_{p,q}$ consists of the closure of $D_{p,q}$ some pieces of $L$ and $E$ and possibly some pieces of leaves of $\wcF_1$ and $\wcF_2$ non separated from $L$ and $E$ respectively. 
\item In particular, if $L' \in \wcF_1$ and $E' \in \wcF_2$ intersect $V_{p,q}$, then, their intersection in $V_{p,q}$ (i.e. $L' \cap E' \cap V_{p,q}$) is a connected arc whose closure has both endpoints in $D_{p,q}$. 
\end{enumerate}
\end{prop}

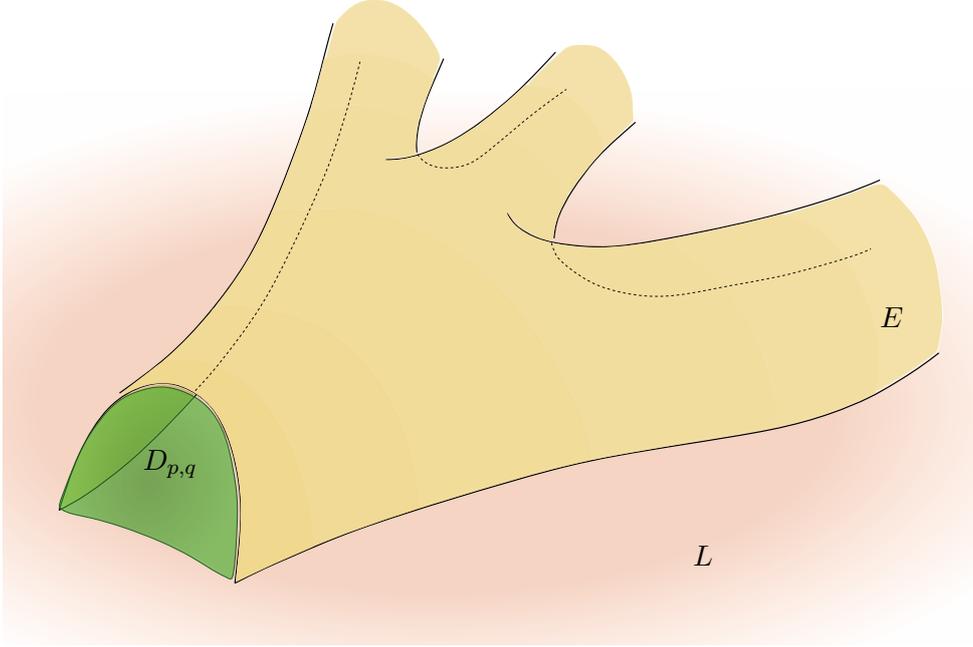
\begin{figure}[htbp]
\begin{center}
\includegraphics[scale=0.59]{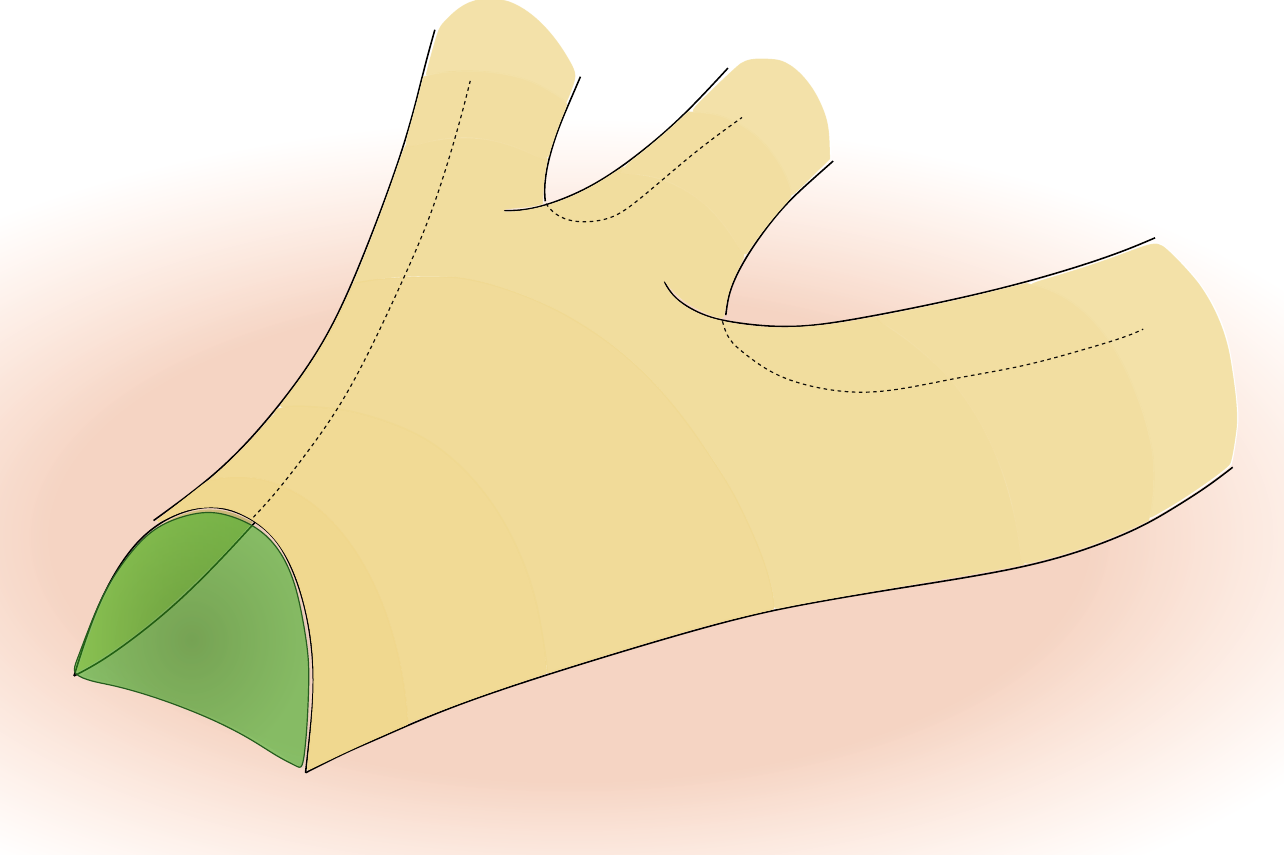}
\begin{picture}(0,0)
\put(-316,66){$D_{p,q}$}
\put(-40,120){$E$}
\put(-110,30){$L$}
\end{picture}
\end{center}
\vspace{-0.5cm}
\caption{{\small The region $V_{p,q}$ between the leaves bounded by the choice of $D_{p,q}$.}\label{f.vpq}}
\end{figure}

By \emph{well defined} we mean that the set $V_{p,q}$ is defined by those properties once the disk $D_{p,q}$ is chosen. See Figure \ref{f.vpq}.

\begin{proof}
This is just a region defined by the choice of the disk $D_{p,q}$ in Proposition \ref{prop-nonsepintersection} so that it sees the non-separated region. To define it properly we consider leaves $L_n \to L$ from above and $E_n \to E$ from above (that is, they all intersect the $(+,+)$ quadrants of $\tilde O_p$ and $\tilde O_q$ in $D_{p,q}$).

Let $\beta'_n = L_n \cap D_{p,q}$ and $\gamma'_n = E_n \cap D_{p,q}$.
These are compact arcs in $D_{p,q}$ with boundary in the interior
of $\alpha_2$ and $\alpha_1$ respectively. 
For $n$ big they intersect in well defined points $p_n, q_n$
so that $p_n \to p, q_n \to q$ when $n \to \infty$.
Let $\beta_n$ be the compact subarc of $\beta'_n$ between
$p_n, q_n$, and similarly define $\gamma_n$.

We claim that $p_n, q_n$ are also the boundary points of uniquely
defined compact arcs $\delta_n$ which are in $L_n \cap E_n$.
To prove this, let $p'_n, q'_n$ be the endpoints of $\gamma_n$,
with $p'_n \to p, q'_n \to q$. Notice that $p'_n, q'_n$
are in the same leaf $E_n$. By hypothesis of non separation of
$r_1, r_2$ in $\cG_L$ it follows that $p'_n, q'_n$ are in
same leaf $g_n$ of $\wcG$ in $L$. Let $c_n$ be the compact arc 
in that leaf between $p'_n$ and $q'_n$. Hence $c_n \to r_1 \cup r_2$
(and maybe converges to other $\wcG$ leaves as well).
Then $c_n \cup \gamma'_n$ is a simple closed curve in $E_n$,
which therefore bounds a unique, well defined disk $D'_n \subset E_n$.
This disk $D'_n$ intersects $L_n$ in a compact arc,
denoted by $\delta_n$ and which has endpoints in $p_n, q_n$.
In particular $\delta_n$ is a subset of $E_n \cap L_n$ and
it is the desired arc.

 We consider the open balls $B_n$ obtained as the region in $\mt$ bounded by the following disks: 

\begin{itemize}
\item The disk $D_{E_n}$ inside $E_n$ bounded by 
$\delta_n \cup \gamma_n$,
\item The disk $D_{L_n}$  inside $L_n$ bounded by $\beta_n \cup \delta_n$,
\item The subdisk $\hat D_n \subset D_{p,q}$ bounded by the union 
$\beta_n \cup \gamma_n$.
\end{itemize}

The union of these disks gives 
a (tamely) embedded sphere in $\mt$ and thus (since
$M$ is irreducible) it bounds an open  ball which we denote
by $B_n$. 

We consider $V_{p,q}$ to be the union in $n$ of all the balls $B_n$. It is not hard to see that this set is well defined independently on the choice of the sequences $L_n$ and $E_n$ $-$ 
because if $L'_n, E'_n$ are other sequences converging to $L$ and $E$
respectively, then 
$L_n, L'_n$ intercalate, and so do $E_n, E'_n$. 
Hence there is an inclusion between balls
$B_n$ and $B'_m$ for $m$ sufficiently big, given $n$.

To verify the properties: first notice that if $F \in \wcF_1$
intersects $V_{p,q}$ then it intersects the compact disk $D_{p,q}$
in the interior. From this, and the way $V_{p,q}$ was expressed
as an increasing union of balls with compact closure, explicitly described
by their boundaries, it is immediate that properties (1)
and (3) follow. 
To prove property (2): clearly the boundary of $V_{p,q}$ contains
the closure of $D_{p,q}$. The rest of the boundary is obtained
as the limit of the sequences $D_{L_n}$ and $D_{E_n}$.
We consider only $D_{L_n}$. Certainly it converges to a subset
in $L$ which contains $r_1$ and $r_2$ and all other leaves
of $\wcG_L$ non separated from these rays and in between them,
and the region bounded by all of these and also $\alpha_1$. 
Any other limit points have to be in leaves non separated from
$L$. 
This proves the result.
\end{proof}


\section{Non-connected intersections and Reeb surfaces}\label{s.B2}

In this section we prove Theorem \ref{teo-noninjective}. The strategy is as follows: we first apply a result from \cite{BarbotFP} to obtain a pair of leaves with doubly non-separated intersection in $\mt$ (as defined in \S~\ref{ss.nonsepinter}). The goal is to show that, when projected to $M$,
these doubly non-separated leaves bound a Reeb surface. 
The proof will be by contradiction assuming that is not the case. 

Given a double non-separation, we consider a $3$-dimensional region
$V_{p,q}$ as defined in the previous section, once the disk 
$D_{p,q}$ is chosen. Let $r_1, r_2$ be the initial non separated
rays. Assume that
$\pi(r_1)$ is not closed. We look at accumulation points of $\pi(r_1)$, that is, the returns of $\pi(r_1)$ to a given foliated box. Using this we will produce what we call a {\em serrated set}: this is obtained from a piece of the region $V_{p,q}$ which projects into a solid torus in a convenient quotient of $\mt$. The returns of $\pi(r_1)$ to a foliated box, correspond to different deck transformations. This will be shown to be impossible unless the deck transformations in question belong to a common cyclic group. This will allow to show that the ray $\pi(r_1)$ is accumulating in some way which contradicts the expansion in holonomy produced by the double non-separation. This contradiction will show that $\pi(r_1)$ is closed, and then since the same argument applies to $\pi(r_2)$, we will obtain a Reeb surface in $M$. 

Note that the model example we have for this section is the example from \cite{MatsumotoTsuboi} that was studied in detail in \cite[\S 7]{FPt1s}. In that case, the doubly non-separated intersection occurs between lifts of cylinder leaves of $\cF_1$ and $\cF_2$, and the boundaries consist of closed curves. This is exactly what we aim to show here that is always the case. 

\subsection{Non injectivity implies there is doubly non separated intersection} 
We will use here the properties shown in \cite[\S 8]{BarbotFP} to get a good configuration from the hypothesis of Theorem \ref{teo-noninjective}. Note that this section considerably simplifies if one assumes that both $\cF_1$ and $\cF_2$ are $\RR$-covered. 

We will precisely state \cite[Theorem 8.1]{BarbotFP}: 

\begin{teo}\label{teo-nonsep}
Let $\cF_1, \cF_2$ be  transverse foliations and $\cG = \cF_1 \cap \cF_2$. Assume there are leaves $L \in \wcF_1$ and $E' \in \wcF_2$ which intersect in at least two leaves $\ell'_1, \ell'_2 \in \wcG$ so that $\ell'_1$ and $\ell'_2$ are non-separated in $\cG_{E'}$. Then there is $E \in \wcF_2$ so that the
intersection $E \cap L$ contains leaves $\ell_1,\ell_2 \in \wcG$ which are non separated in both $\cG_E$
and $\cG_{L}$, that is, the leaves $E$ and $L$ have a doubly non-separated intersection in the sense of \S~\ref{ss.nonsepinter} above.
\end{teo}

\begin{cor}\label{cor-nonsepinter}
Let $\cF_1, \cF_2$ be transverse foliations and assume there are leaves $L' \in \wcF_1$ and $E' \in \wcF_2$ intersecting in more than one connected component. Then, there is a pair of leaves $L \in \wcF_1$ and $E \in \wcF_2$ such that they have a doubly non-separated intersection. 
\end{cor}

\begin{proof}
In order to apply Theorem \ref{teo-nonsep} it is enough to show that if $L' \cap E'$ has more than one connected component, then, there is a leaf $L \in \wcF_1$ so that the intersection $L \cap E'$ has distinct components
(leaves of $\wcG_{E'}$) which are non-separated in $\cG_{E'}$.

To obtain this result we apply \cite[Claim 8.2]{BarbotFP} so let us explain the setup and terminology used there. The leaf 
space of $\wcG_{E'}$ is a simply connected one dimensional manifold, which is
non Hausdorff. Between any two points $x, y$ in this leaf space
there is a minimal ``path" (denoted by $[x,y]$)
 which is the set of leaves that have
to be intersected by any path in $E'$ from leaf $x$ to leaf $y$.
This path is a union $\cup_i [x_i,y_i]$ where $[x_i,y_i]$ are
segments (homeomorphic to $[0,1]$) in the leaf space, and 
$y_i$ is non separated from $x_{i+1}$, see details in \cite{BarbotFP}.
So we consider this path for $x$ a leaf of $\wcG_{E'}$ in $E' \cap L'$ and
$y$ another leaf in $E' \cap L'$. The leaves $x, y$ are distinct so the
path $[x,y]$ is not a singleton. Since they are both contained
in $L'$ then $[x,y]$ cannot be a single segment for then
it would produce a transversal to $\wcG_{E'}$ in $E'$ from $x$ to $y$,
so a transversal in $\mt$ from $L'$ to itself. It follows
that $[x,y]$ has more than one segment. Then \cite[Claim 8.2]{BarbotFP}
shows that there is some $i$ so that $y_i$ and $x_{i+1}$ are
in the same leaf of $\wcF_1$.
In sum, if we let $L$ be the
leaf of $\wcF_1$ containing both $x_{i+1}, y_i$, then 
$x_{i+1}, y_i$ are distinct components of $L \cap E'$ which
are non separated in $\wcG_{E'}$ completing the proof.
\end{proof}

A key property of the doubly non-separated intersection is the following (recall that by Lemma \ref{lem-nonseprays} we have that the non-separated rays are well defined): 

\begin{prop}\label{prop-expandingholonomy}
Let $L \in \wcF_1, E \in \wcF_2$ with doubly non-separated intersection and let $r_1$ be a ray in $L\cap E$ which is non-separated in both $\cG_L$ and $\cG_E$ with some other ray in $L\cap E$. Then, for every $x \in r_1$ and $\delta>0$ one has that there is $R:=R(x,\delta)>0$ so that if $y \in r_1$ verifies that $d_{r_1}(x,y)>R$ then one has that the domain $A_{x,y}$ of the holonomy map $\Psi_{x,y}$ (see \S~\ref{ss.transversals}) is strictly contained in the $\delta$-neighborhood of $x$ in the quadrant $\tilde O_x$ associated with the double non-separation.
\end{prop}
The same happens when projecting by $\pi: \mt \to M$ the universal covering projection.

\begin{proof}
Let $\ell_1$ be the leaf of $\wcG$ containing the ray $r_1$ and $L \in \wcF_1$ the leaf containing $\ell_1$. Since $r_1$ is non separated
from $r_2$ in $L$, it follows that for any
$z$ near $x$ in $L$ in the component of $L  \setminus \ell_1$
containing $r_2$, it follows that eventually the $\wcG$ leaf
of $z$ is more distant in $L$ from $\ell_1$ than the size of
a foliated box of $\wcG$. The same happens in $\wcF_2$, 
and this proves the result.
\end{proof}

Theorem \ref{teo-noninjective} is a direct consequence of the following result whose proof will occupy the rest of this section. 

\begin{teo}\label{teo-doublynsreeb}
Let $\cF_1, \cF_2$ be transverse foliations and assume that $L \in \wcF_1$ and $E \in \wcF_2$ have doubly non-separated intersection. Then
the projections of $L$ and $E$ to $M$ contain Reeb surfaces. 
\end{teo}

The result is unchanged by taking finite lifts, so for simplicity of arguments, we will assume in the remainder of this section that $M$ is orientable and that $\cF_1, \cF_2$ are transversally orientable (therefore also $\cG$ is orientable). 

\subsection{The setup: defining types of returns}\label{ss.setup}

We consider a pair of doubly non-separated rays $r_1, r_2$ 
contained in $L \cap E$ where $L \in \wcF_1$ and $E \in \wcF_2$ have doubly non separated intersection. 
We assume that $r_1$ is non separated from $r_2$ in
its positive side in both $L$ and $E$. This means the following:
given the induced transverse orientation in $\cG_L$ from $\wcF_2$ we have
that $r_2$ is in the positive side component of $L \setminus r_1$,
and similarly in $E$. 

By compactness of $M$, it follows that $\pi(r_1)$
must accumulate somewhere in $M$. We consider $p_n \in r_1$ so that $\pi(p_n) \to x_\infty$. 

Up to changing the sequence, we can assume that for every $n$ we have that all the $\pi(p_n)$ belong to the $\eps/2$-neighborhood of $x_\infty$ and thus to a foliated box for $\cF_1$ and $\cF_2$. In particular, we can assume there are deck transformations $\gamma_n \in \pi_1(M)$ so that $d(\gamma_n p_n , p_0) <\eps$. Here we are considering $\eps \ll \eps_0$ of the local product structure (see \S~\ref{ss.transversals}).

We assume without loss of generality that $p_n$ are ordered in $r_1$ so that if $m>n$ then the point $p_m$ is farther in $r_1$ from $p_0$ than $p_n$. 

Using Proposition \ref{prop-expandingholonomy} and taking a subsequence, we can further assume that the length in $r_1$ from $p_n$ and $p_{n+1}$ in $r_1$ is always larger than $R$ given by that proposition so that the holonomy expands points at distance $\eps \ll \eps_0$ to be at distance larger than $\eps_0$. We will further consider subsequences later, but keep the notation $p_n$ and $\gamma_n$. 

We will consider $p_0$ fixed and fix also a point $q_0$ in another component of $L\cap E$, which is doubly non separated from 
$r_1, r_2$. This component may or may not be the leaf containing $r_2$. Let $D_0$ be a disk given by Proposition \ref{prop-nonsepintersection} (that is $D_0 = D_{p,q}$) associated to these points and similarly $V_0$ ($= V_{p,q}$) the set given by Proposition \ref{prop-tubeconstruction}. 
To obtain $D_0$ we pick compact arcs $\alpha_1$ in $L$ from $p$ to $q$
and $\alpha_2$ in $E$ from $p$ to $q$ with same boundary points. 
Here $\partial D_0 = \alpha_1 \cup \alpha_2$.

We will classify the \emph{returns} $\gamma_n p_n$ depending on the local position with respect to $p_0$ and the intersection between $L$ and $E$. In particular, we will say that $\gamma_n L$ is \emph{above} $L$ if it intersects the connected component of $\mt \setminus L$ containing $V_0$. Similarly, $\gamma_n E$ is above $E$ if it intersects the connected component of $\mt \setminus E$ containing $V_0$. Note that being above is being in the direction where we know that the holonomy expands (according to Proposition \ref{prop-expandingholonomy}).  A priori we do not know whether the other direction has expanding holonomy or not, and as we will see throughout the proof that this complicates the analysis.

We then classify the return $\gamma_n p_n$ as $(+,+)$, $(+,-)$, $(-,+)$ and $(-,-)$ acording to Convention  \ref{convention-rays}  (see Figure \ref{f.sign}):

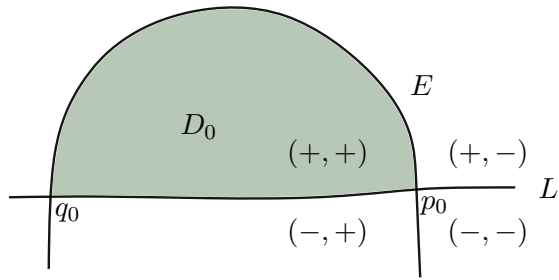
\begin{figure}[htbp]
\begin{center}
\includegraphics[scale=0.79]{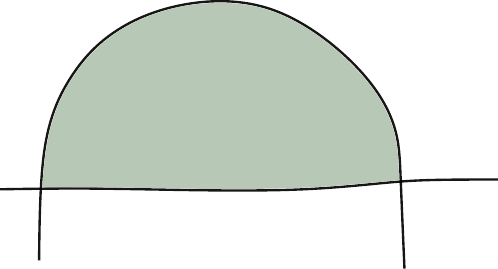}
\begin{picture}(0,0)
\put(5,30){$L$}
\put(-43,69){$E$}
\put(-39,26){$p_0$}
\put(-89,44){$(+,+)$}
\put(-29,44){$(+,-)$}
\put(-29,14){$(-,-)$}
\put(-89,14){$(-,+)$}
\put(-129,54){$D_0$}
\put(-176,23){$q_0$}
\end{picture}
\end{center}
\vspace{-0.5cm}
\caption{{\small The sign convention for returns.}\label{f.sign}}
\end{figure}

We will call the returns of the form $(+,-)$ and $(-,+)$ of \emph{mixed type}. Note that in the above list we are not considering the cases where $\gamma_n L = L$ or $\gamma_n E = E$ which will be considered separately. 

\subsection{Returns in the same leaf} 

This section studies the case where $\gamma_n L = L$ (the case $\gamma_n E = E$ is analogous). This is the easiest case because we can produce a 2-dimensional argument. Some of the ideas later are extensions of this case. In this case we will show:

\begin{prop}\label{prop-goodreturns}
If there is an infinite sequence of $\gamma_n$ such that $\gamma_n L = L$ then the ray $r_1$ projects to a closed curve in $M$.\end{prop}

\begin{proof} 
We will show that $\gamma_n E = E$ for $n$ sufficiently big.
Assume then by contradiction, that up to subsequence $\gamma_n E \not
= E$ for all $n$. 

There are two cases: the first case is that
$\pi(r_1)$ returns outside (negative side), that is, one has that $\gamma_n E$ is strictly below $E$. In this case we will find a closed leaf of $\cG$ in the projection of $L$ and this will lead to a contradiction.
Let $\beta$ be a compact segment in $L$ transverse to $\cG_L$, starting
in $p_0$ and ending in a point $y$ which is  $\delta$ distant from $p_0$.
We can take $\beta$ to be a subsegment of $\alpha_1$.
Let $I$ be the compact interval of leaves of $\wcG_L$ intersected
by $\beta$.
Now choose $n$ big enough so that the holonomy from $p_0$ to $p_n$
expands $\beta$ so that the final result has endpoints more than $\eps_0$ apart. 
Consider the corresponding $\gamma_n$ and the projections to $M$ of the returns.
In other words holonomy pushes $\beta$ to a transverse segment $J$
with one endpoint in $p_n$ and the other more than $\eps_0$ away $-$ 
this is because of expanding holonomy
(Proposition \ref{prop-expandingholonomy}).
Then $\gamma_n$ maps $J$ to a transverse segment 
which can be arranged to contain $I$ in its interior.
This is immediate on the $y$ end because of expanding holonomy.
It also follows in the $p_0$ end because $\gamma_n E$ is smaller
than $E$. Therefore this holonomy in $L/_{< \gamma_n >}$ has fixed points
between $\pi(p_0)$ and $\pi(y)$, meaning $\cG$ has a 
closed curve intersecting $\pi(\beta)$ between $\pi(p_0)$
and $\pi(y)$.
We choose the fixed point of holonomy closest to $\pi(p_0)$, and let $x$ the lift
to $\beta$.
Therefore the leaves of $\wcG$ between the leaf through $x$ and the leaf through $p_0$ 
project to leaves of $\cG$ which spiral towards a closed
leaf. In particular they are not closed leaves.

Now reapply this argument to $m \gg n$ and also to a subsegment
$\beta_1$ of $\beta$ with an endpoint $y'$ which is closer to $p_0$ than
$x$ is. This produces a fixed leaf of $\gamma_m$ intersecting
$\beta_1$ between $y'$ and $p_0$, projecting to a closed
leaf of $\cG$ $-$ but we showed in the previous paragraph
that all of these $\cG$ leaves are not closed, a contradiction.
This finishes the analysis in this case.

Notice that the arguments above also show that there is $\delta > 0$ so that
leaves of $\wcG_L$ not containing $r_1$ but having points within 
$\delta$ of $p_0$ do not project to closed leaves of $\cG$ in $M$.
Indeed, if there is a leaf near $p_0$ which projects
to a closed curve in $M$, the holonomy based
on the curve going in the direction that $\pi(r_1)$ goes, must be expanding (else $r_1$ could not be non-separated from other leaves) and thus contracting in the other direction (and therefore there cannot be leaves projecting to closed curves in between). 

Now consider the second case, that the returns 
are inside, that is, $\gamma_n E$ is strictly above $E$.
 we will use the infinitely many returns to get a contradiction. 
Choose the $p_n$ so that holonomy going back to $p_{n-1}$ sends
any transversal of length $\leq \eps_0$ to a transversal of
length less than $\delta$.

First note that
 up to subsequence, we can assume  that the points $\gamma_n p_n$  are ordered increasingly with respect to some transversal to $\wcF_2$:
otherwise if for some $m > n$ we have that $\gamma_m E$ is below $\gamma_n E$ then we have that we can apply the arguments of the first case 
 to get a leaf of $\wcG$ fixed by $\gamma_n^{-1}\gamma_m$ at distance less than $\delta$ from $p_0$ a contradiction.

Finally, since the returns happen in a neighborhood of $p_0$ we can assume that the sequence $\gamma_n p_n \to z_\infty$. 
Now consider $k$ large enough so that the distance between $\gamma_k p_k$
and all other $\gamma_n p_n$ as well as $z_{\infty}$ is smaller than 
$\delta$. Now consider the holonomy going in the backwards direction
from $p_k$ to $p_0$. 
Using Proposition \ref{prop-expandingholonomy}, we can choose $k$ 
big enough so that the holonomy back from $p_k$ to
$p_0$ sends all points $\gamma^{-1}_k \gamma_n p_n$ (with $n > k$) and $\gamma^{-1}_k z_{\infty}$ to be $\ll \delta$ close to $p_0$.
These points correspond to returns of $\pi(r_1)$ and they 
are all closer to $\pi(p_0)$ than $\pi(\gamma_1 p_1)$. This contradicts what was 
proved in the beginning of the argument in the second case and completes the proof.
\end{proof}

The following will be useful in later arguments:

\begin{lema}\label{lem-monotone} 
If $r_1$ does not project into a closed curve, then, up to subsequence, we can assume that all the returns $\gamma_n p_n$ are of the same type. 
In addition we can assume that $\gamma_n L, \gamma_n E$ are both
strictly monotonic sequences.
\end{lema}

\begin{proof}
First assume there is a subsequence $\gamma_n$ so that $\gamma_n L$ 
is constant. Then replace $L$ by $\gamma_1 L$ and $\gamma_n$ by
$\gamma_n \circ \gamma_1^{-1}$, and apply the previous lemma
to show that $\pi(r_1)$ is closed.
Therefore there is subsequence still denoted by $\gamma_n$ so 
that $\gamma_n L$ are pairwise distinct, and so are $\gamma_n E$.
Now every pairwise distinct sequence in a closed interval has a subsequence which is either increasing of decreasing. Hence we can assume that up to subsequence, the returns in each foliations are all monotonic (either increasing or decreasing) and thus, all returns are of the same type with respect to the previous returns. 
\end{proof}

\subsection{Restrictions on returns} 

From now on we will assume up to subsequence that no returns
are on the same $\wcF_1$ or $\wcF_2$ leaf, that is, for
all $n$, $\gamma_n L \not = L$ and $\gamma_n E \not = E$.
The goal of this subsection is to we show the following: 

\begin{prop}\label{prop-noinside}
Up to subsequence, one can assume that all returns $\gamma_n p_n$ are of the form $(+,-)$ or $(-,+)$ (i.e. of mixed type). 
\end{prop}

In other words, in this subsection we, up to subsequence, 
 will rule out the returns
of type $(-,-)$ or $(+,+)$.
First, we will show that the returns below produce fixed leaves using the expansion of holonomy. This will be used again later: 

\begin{lema}\label{lem-fixedleaf-}
Assume that $\gamma_n L$ is below $L$, then, there is a leaf $L_n$ above $L$ which is fixed by $\gamma_n$. Moreover, this leaf intersects $V_0$. 
\end{lema}

\begin{proof}
This is a direct application of Proposition \ref{prop-expandingholonomy}, as explained in the proof of Proposition \ref{prop-goodreturns}. 
\end{proof}

The following again is very similar to the idea in Proposition \ref{prop-goodreturns}, but now the argument is 3-dimensional. 

\begin{lema}\label{lem-imposible--}
There cannot be a subsequence of returns $\gamma_n p_n$ consecutively of
type $(-,-)$ with respect to the previous one. 
\end{lema}

\begin{proof}
Suppose there are infinitely many such returns.
Using the previous lemma we get a pair of leaves $L_n \in \wcF_1$ and $E_n \in \wcF_2$ fixed by $\gamma_n$ and very close to $L$ and $E$. We choose $L_n$ so that $\gamma_n$ does not fix any leaf between $L_n$ and $L$ and similarly for $E_n$ and $E$. The $\gamma_n$ is obtained by tracking $r_1$, therefore it follows that $\gamma_n$ fixes the $\wcG$ leaf $\ell$ in $L_n \cap E_n$ with
a point very close to $p_0$. Now look at the rectangle of leaves of $\wcG$ with corners in $r_1$ and $\ell$. Then $\gamma_n$ fixes $\ell$ and for any other $\wcG$ leaf $\ell'$ in this rectangle, it follows that $\gamma_n^{-i}(\ell')$ converges to $\ell$ when $i \to \infty$. In particular for any leaf $\ell'$ in the rectangle and not $\ell$, then $\pi(\ell')$ is not a closed leaf.

Now take $m \gg n$. Since the $\gamma_m$ returns are in the $(-,-)$ position with respect to the $\gamma_n$ returns (see Lemma \ref{lem-monotone}), the same argument produces a leaf of $\wcG$ fixed by $\gamma_m \circ \gamma_n^{-1}$ and very close to $\gamma_n p_n$. In fact by letting $m$ as big as we want we can get the leaf as close as we want to $\gamma_n p_n$.
Mapping by $\gamma_n^{-1}$ we obtain a leaf $\ell"$ of $\wcG$ as close as we
want to $p_n$ and in the $(+,+)$ position with respect to $p_n$, 
and $\ell"$ projects to a closed leaf of $\cG$.
Now use holonomy along $r_1$ from $p_n$ to $p_1$ (that is, 
moving backwards along $r_1$ which contracts points): then $\ell"$ has a point 
arbitrarily close to $p_0$ and in the $(+,+)$ quadrant. In particular
$\ell"$ is in the rectangle we defined previously and is not $\ell$.
This contradicts the property of the rectangle.
This finishes the proof of the lemma.
\end{proof}

\begin{lema}\label{lem-imposible++}
There cannot be a subsequence of 
returns $\gamma_n p_n$ consecutively of type $(+,+)$ with
respect to the previous one.
\end{lema}

\begin{proof}
If this happens we will produce another sequence which has 
returns consecutively of type $(-,-)$, which will be a contradiction
by the previous lemma.

We have that $\gamma_n p_n \to z$, all points are  in a small neighborhood of $p_0$, and every point is
in consecutive $(+,+)$ position with respect with the previous one.
In addition we know that for all sufficiently large $n$, we have that the holonomy from $p_0$ to $p_n$ along $r_1$ maps a small neighborhood of $p_0$ in its $(+,+)$ quadrant (containing $\gamma_n p_n$) to a set containing the full $\tilde O_{p_n}$ quadrant. 

We want to apply holonomy backwards. This is better done in $M$:
we are doing holonomy backwards from $\pi(p_n)$ to $\pi(p_0)$.
For every $m \gg n$ we have that holonomy backwards sends
$\pi(p_m)$ to a point $q_n$ in $O_{\pi(p_0)}$, and $q_n$
is in the $(-,-)$ position with respect to $\pi(\gamma_n p_n)$.
For $m$ big these points are in $\pi(r_1)$  and after
$\pi(\gamma_n p_n)$ $-$ this is because
holonomy from $\pi(\gamma_n p_n)$ to $\pi(p_0)$ moves a bounded
amount and $p_m$ is farther and farther in the future from $p_0$
in $r_1$. 
This now produces arbitrarily long returns from $\gamma_n p_n$ 
which are in the $(-,-)$ position with respect to $\gamma_n p_n$.
In addition these returns are consecutively in the $(+,+)$ position
with respect to
to previous ones. Now pick one such return and call it $w_1$.
We can restart with $w_1$: the same argument produces $w_2$
which is in the $(-,-)$ position with respect to $w_1$
and is in the future of $w_1$. In this way we produce a
sequence of consecutive returns starting from $\gamma_n p_n$,
so that they are consecutively in the $(-,-)$ position from
the previous one.
This was shown to be impossible in Lemma \ref{lem-imposible--}.
This finishes the proof of the lemma.
\end{proof}

Putting together Lemmas \ref{lem-imposible++} and \ref{lem-imposible--} with Proposition \ref{prop-goodreturns}, we complete the proof of Proposition \ref{prop-noinside}.

\begin{convention}\label{conv-returns}
Since we assumed  that $r_1$ does not project to a closed curve, we have that there are infinitely many returns, and up to changing the roles of $\cF_1$ and $\cF_2$ we can and we will assume without loss of generality that all returns are of the form $(-,+)$. By Lemma \ref{lem-monotone}, we can take a sequence of returns ordered by their position in $r_1$ so that each return is of the for $(-,+)$ with respect to the previous one. Our goal is to show that there are infinitely many such returns for which the deck transformations $\gamma_n$ belong to the same cyclic group of $\pi_1(M)$. This will allow us to obtain a contradiction with the expansion of holonomy given by Proposition \ref{prop-expandingholonomy}. 
\end{convention}

\subsection{A serrated set}

Here we construct a serrated set $\cS_n$ for some sufficiently large return $\gamma_n$. Recall that $V_0$ is the $3$-dimensional set as constructed
in the previous section, given a choice of disk $D_0$. We assume that $\gamma_n p_n$ is in $(-,+)$ possition with respect to $p_0$ by Convention \ref{conv-returns}. 
The following statement defines a serrated set $\cS_n$ associated to $\gamma_n$ and states the properties we will need. The rest of this subsection is devoted to proving this statement: 

\begin{prop}\label{prop-serratedset}
Let $n$  be sufficiently large so that $\gamma_n p_n$ is in $(-,+)$ position with respect to $p_0$. Then, there exists a closed subset $S_n$ of $\overline{V_0}$ which verifies:
\begin{enumerate}
\item the interior of $S_n$ is disjoint from $\gamma_n S_n$ (and hence disjoint from $\gamma_n^{-1} S_n$ as well), 
\item the interior of the projection of $S_n$ to $\mt/_{\langle \gamma_n \rangle}$ is homeomophic to $\RR^2 \times S^1$, and defines a \emph{serrated set} $\cS_n  = \bigcup_{k \in \ZZ} \gamma_n^k S_n$, which is also a lift of the projection of $S_n$ to $\mt/_{\langle \gamma_n \rangle}$, 
\item the boundary of $S_n$ consists of: $D_0$, a disk $D_n$ whose boundary contains an arc in $L$ and an arc in $E$ and has $p_n$ in a corner, and some subsets of $L$ and $E$ (whose boundaries are either the boundary curves of $D_0$ and $D_n$ or pieces of leaves of $\cG_L, \cG_E$ non-separated from $r_1$), 
\item the intersection of $L_n$, the leaf fixed by $\gamma_n$ intersecting $V_0$ and closest to $L$, and $\cS_n$ is a band $B_n$ whose projection to $\mt/_{\langle \gamma_n \rangle}$ is a compact annulus $A_n$. 
\end{enumerate}
\end{prop}

\begin{figure}[htbp]
\begin{center}
\includegraphics[scale=0.59]{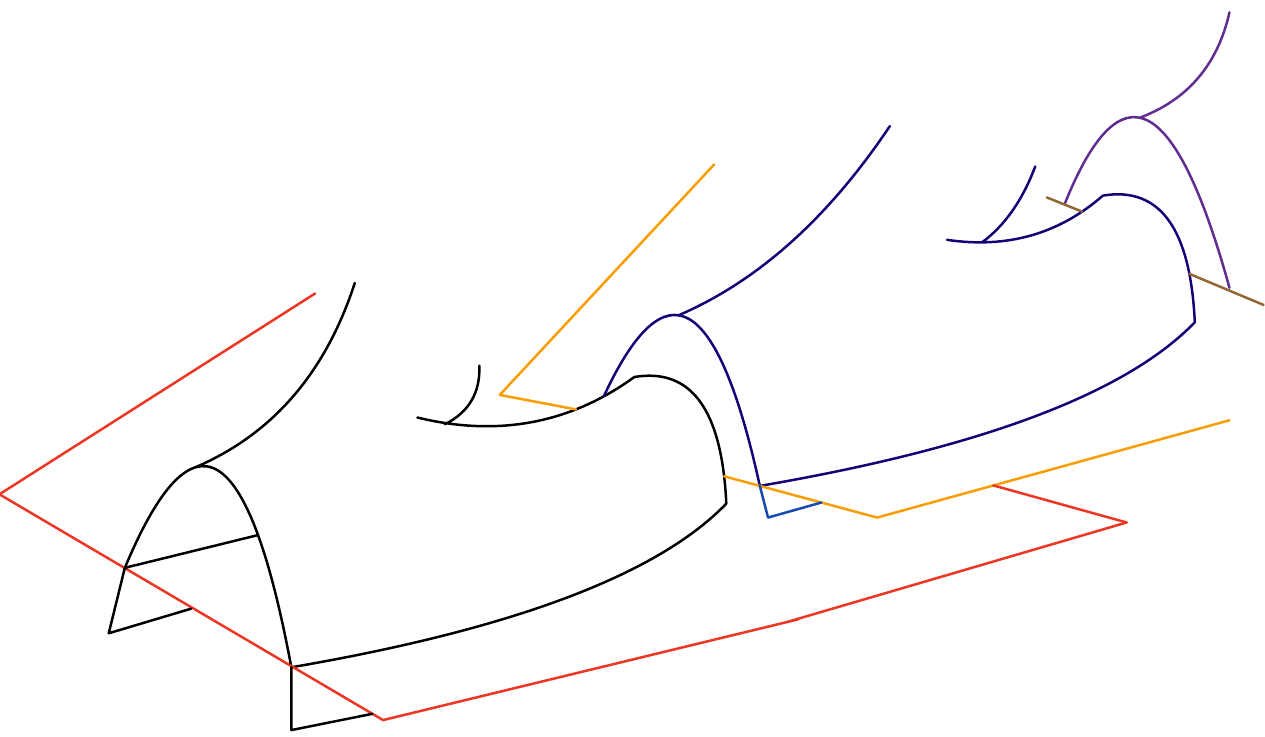}
\begin{picture}(0,0)
\put(-130,96){$\gamma_n^{-1} S_n$}
\put(-220,69){$S_n$}
\put(-340,16){$E$}
\put(-189,12){$L$}
\end{picture}
\end{center}
\vspace{-0.5cm}
\caption{{\small The serrated set $\cS_n$ in $\mt$.}\label{f.serrated}}
\end{figure}

\begin{proof}
We have fixed the disk $D_0$ constructed in \S~\ref{ss.setup} as in Proposition \ref{prop-nonsepintersection} with boundaries $\alpha_1 \subset L$ and $\alpha_2 \subset E$ and corners in $p_0$ and some point $q_0$ in a leaf $\ell_2$ doubly non separated from $\ell_1$, which is  the leaf containing the ray $r_1$. Choose $n$ big enough so that
$d(\gamma_n D_0, D_0) \gg Diam(D_0)$, and also that $\gamma_n p_n$ is
very close to $p_0$. Here $\gamma_n p_n \in O_{p_0}$ is in the $(-,+)$
quadrant of $p_0$, and thus,
the leaf $\gamma_n E \in \wcF_2$ intersects $D_0$ in an arc, which we denote
by $\beta_0$ and has boundary points $p', q'$ in
$\alpha_1$. By switching $p', q'$ if necessary,
assume that $p'$ 
is very near $p_0$ (that is $p' \in O_{p_0}$ is the point near $\gamma_n p_n$).
Then $\beta_0$ together
with a compact subarc $c'_0$ of $\alpha_1$ bounds a subdisk of $D_0$,
which we denote by $D^*_1$. 
Let $c_0 = \gamma_n^{-1}(c'_0)$.

Consider $\gamma_n^{-1}(D^*_1)$.
This is a disk that has a corner $\gamma_n^{-1}p'$
very close to $p_n$. Extend the curve $\gamma_n^{-1}(\beta_0)$
very slightly beyond $\gamma_n^{-1}(p')$
along $E$ until it hits $L$. Notice that $\gamma_n^{-1}(L)$ is very near $p_n$ as
$L$ is very near $\gamma_n p_n$. Denote the point 
where it hits $L$ by $p"$ (one could even choose $p"$ to be $p_n$ as it belongs to $r_1$ and is very close to $p_n$).  Let the segment from $\gamma_n^{-1}p'$
to $p"$ be denoted by $c_1$. 

Next we want to extend $\gamma_n^{-1}(\beta_0)$ along $E$,
in the other
end until it hits a $\wcG$ leaf non separated from 
$r_1$.
Do this in such a way that the extension stays far from
$D_0$, ending in $q"$.  This new arc is denoted by $c_2$ (note that this arc may not be so short a priori).
Now connect the endpoints
$p", q"$ by a compact arc $c_3$ in $L$ still avoiding
$D_0$.
The union 

$$c_0 \cup c_1 \cup c_2 \cup c_3$$

\noindent
is a simple closed curve in $\mt$ which bounds a disk in
$\mt$. We can choose such a disk $D_n$ so that its interior
is contained in $V_0$ and verifies properties like the ones in Proposition \ref{prop-nonsepintersection} for $p'',q''$. 
We also choose $D_n$ so that it contains $\gamma_n^{-1}(D^*_1)$.

We let $S_n$ be the region so that the interior is
the subset of $V_0$ contained between $D_0$ and $D_n$ (that is, the closure of $V_0 \setminus V_n$ where $V_n$ is the set as in Proposition \ref{prop-tubeconstruction} associated to $D_n$). See Figure \ref{f.serrated}. 

Property (2) follows from the fact that the interior of $S_n$ is $V_0 \setminus V_n$ so homeomorphic to a ball and the quotient is then homeomophic to $\RR^2 \times S^1$. Property (3) is just by construction. 

Let us now prove (1): 
Notice that the orientation in $r_1$ induces a transverse orientation
in $D_0$. Recall that $\cG$ is orientable and $\gamma_n$
preserves orientation.

Suppose by way of contradiction that there is $z$ in the interior
of $S_n$ so that $\gamma_n^{-1}(z)$ is also in $S_n$.
Now consider $w$ in $\alpha_1$ sufficiently near $p_0$.
Let $C_w$ be the intersection of $E_w$, the $\wcF_2$ leaf containing $w$, with $V_0$.
By the construction of $V_0$ this is a compact disk.
If $w$ was chosen near enough $p_0$ then $C_w$ intersects
$D_n$ in a single compact arc $\beta_2$. This arc
together with an arc in $D_n \cap L$ bounds a disk $D_w$
which is a subdisk of $D_n$.
In addition a subdisk of $C_w$ union a subdisk of $D_0$ and a disk in $L$
bound a ball $B_w$ which is contained in $V_0$.

Notice that if $w$ is near enough $p_0$ then $\gamma_n^{-1}(z)$
is contained in $B_w$, and in the component $Y'$ of $B_w \setminus D_w$
which contains a subdisk in $D_0$.
The disk $D_w$ separates $B_w$ into one part contained
in $V_0$ and another part not contained in $V_0$.

Now apply $\gamma_n$ to $B_w$. Notice $z$ is in $\gamma_n(Y')$
(as $z = \gamma_n (\gamma_n^{-1}(z))$).
This set $\gamma_n(Y')$ has part of the boundary contained
in $D_0$ (this part is contained in $\gamma_n(D_w)$.
This part of the boundary of $\gamma_n(Y')$ separates 
the part of $\gamma_n(Y')$
contained in $V_0$ from the part not contained in $V_0$. 

By the preservation of orientations by $\gamma_n$, it follows
that $z = \gamma_n (\gamma_n^{-1}(z))$ 
belongs to the part not in $V_0$. This contradicts
the assumption and proves
property (1).

\begin{figure}[htbp]
\begin{center}
\includegraphics[scale=0.89]{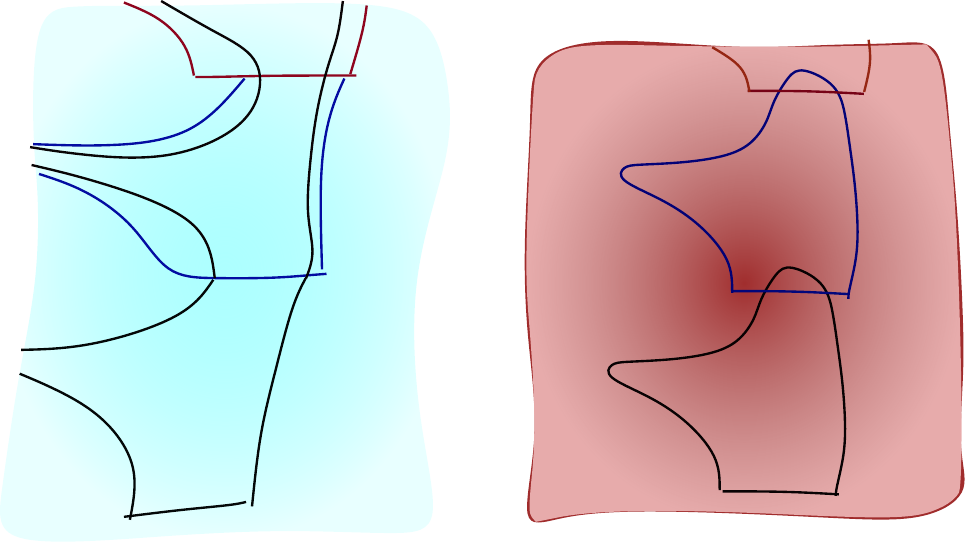}
\begin{picture}(0,0)
\put(105,156){$\gamma_n^{-1} K_n$}
\put(120,69){$K_n$}
\put(120,24){$d_0$}
\put(160,64){$d_1$}
\put(125,114){{\color{blue}$d_2$}}
\put(180,30){$L_n$}
\put(-189,32){$L$}
\put(-90,47){$r_1$}
\put(-120,17){$\alpha_1$}
\put(-190,67){$r_2$}
\put(-180,187){{\color{blue}$\gamma_n^{-1}E \cap L$}}
\end{picture}
\end{center}
\vspace{-0.5cm}
\caption{{\small In the left, the intersection of the serrated set $\cS_n$ with $L$ and in the right the intersection with $L_n$ the lowermost fixed leaf of $\gamma_n$.}\label{f.serratedann}}
\end{figure}
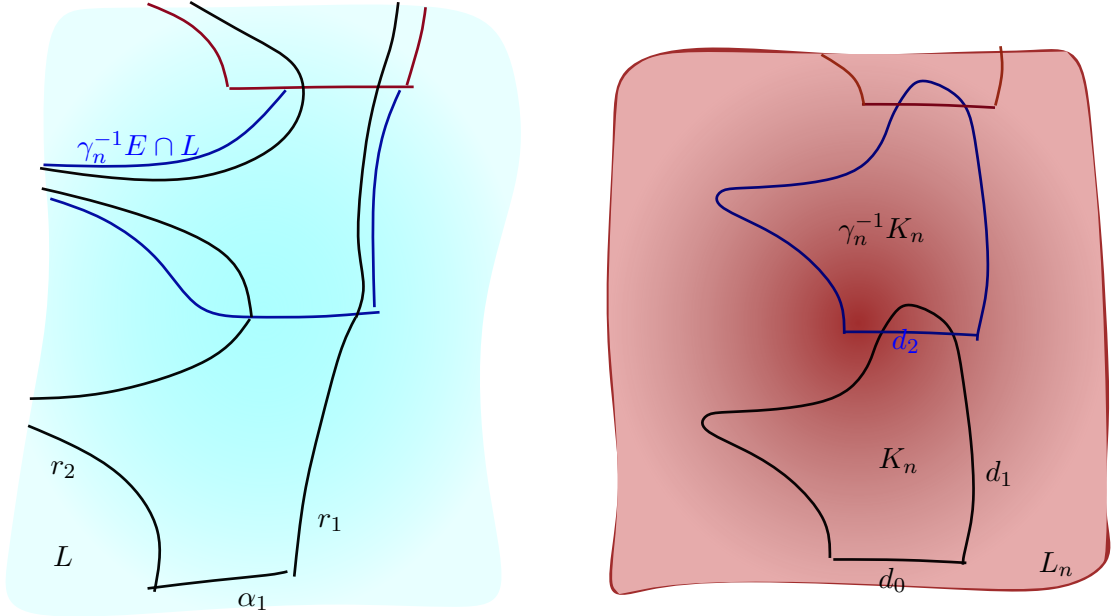

Finally, to get property (4) we note first that since $\gamma_n p_n$ is in the $(-,+)$ quadrant, the expansion of holonomy forces the existence of a compact interval of leaves of $\wcF_1$ where $\gamma_n$ acts in a repelling manner. Let $L_n$ be the lowermost leaf fixed by $\gamma$ in this interval. 
Then $L_n$ intersected with $V_0$ is a compact disk $K_n$ (see Figure \ref{f.serratedann}). 

The boundary of $K_n$ is the union of an arc $d_0 \subset D_0$ with endpoints in $\alpha_2$ and an arc $d_1 \subset L_n \cap E$. The image of $K_n$ by $\gamma_n^{-1}$ is another disk which intersects $K_n$ in a subdisk cutting the arc $d_1$ by some arc $d_2$ ($= \gamma_n^{-1}(d_0)$),
so that its intersection with $K_n$ is a compact arc 
 joining two points in $d_1 \cap \gamma_n^{-1}D_0 \cap K_n$. The union of the $\gamma_n$ orbit of $K_n$ then produces a closed band $B_n$ in $L_n$ which is $\gamma_n$ invariant and the projection $A_n$ of this band to $\mt/_{\langle \gamma_n \rangle}$ is an annulus which coincides with the projection of $K_n$ to $\mt/_{\langle \gamma_n \rangle}$ and therefore is compact. 
\end{proof}

\subsection{Returns in the same cyclic group}
In this section we will use the serrated set constructed in the previous section to obtain a contradiction with the fact that $r_1$ does not project into a closed curve. 

Let us fix $p_n$ and $\gamma_n$ so that $\gamma_n p_n$ is sufficiently close to $p_0$ in the $(-,+)$ position,  so  that the disk $D_n$ can be defined and $\cS_n$ a serrated set as in Proposition \ref{prop-serratedset}.  
Hence we are fixing $n$ and we will consider $m > n$.
The goal of this subsection is to prove the following:

\begin{prop} \label{prop-nomixedreturns} 
The situation of infinitely many $(-,+)$ returns
is impossible.
\end{prop}

The proof will be divided into two parts: first we show
that there is a subsequence $m_i$ so that
all $\gamma_{m_i} \circ \gamma_{m_0}^{-1}$ are contained
in a cyclic subgroup of $\pi_1(M)$. Then we show that this
is impossible given the expansion of holonomy generated by double
non separation.

We use the objects provided by Proposition \ref{prop-serratedset}.
First we consider the leaf $L_n$ of $\wcF_1$ fixed by $\gamma_n$,
intersecting the interior of $V_0$ and the lowest with this
property, that is, $\gamma_n$ does not fix any leaf of $\wcF_1$ between
$L$ and $L_n$. 

We consider the manifold $M_n := \mt /_{\langle \gamma_n \rangle}$ and $\pi_n: \mt \to M_n$ the covering projection.
Now consider the subset $B_n$  inside $\cS_n$ as in Proposition \ref{prop-serratedset} which projects to a compact annulus $A_n$ in $M_n$. The disk $K_n$ (definded in the proof of Proposition \ref{prop-serratedset}) contains a subset $U_n$ which is a fundamental domain for $B_n$ and has boundaries two arcs of $E \cap L_n$, an arc inside $D_0$ and an arc inside $D_n$ (and inside $\gamma_n^{-1}(D_0)$). See Figure \ref{f.serratedann}, it is the set bounded by the arcs $d_0,d_1,d_2$.

We now consider any $m > n$ and analyze how it interacts with
the objects associated with $n$. We show:  

\begin{lema}\label{lem-intersects}
One has that $\gamma_m D_0 \cap  B_n  \neq \emptyset$.
\end{lema}
\begin{proof}
As in Proposition \ref{prop-serratedset}, one can construct a disk $D_m \subset V_0$ so that
$\gamma_m D_m$ has a subdisk contained in $D_0$, and this subdisk 
has a corner very close to $\gamma_m p_m$ and thus to $p_0$
as well. In fact $\gamma_m D_m \cap D_0$ is contained
in $\cS_n$. We will consider the
intersection of $\gamma_m V_0$ with the fixed leaf $L_n$.
The curve $\gamma_m V_0 \cap D_0$ separates 
$\gamma_m V_0 \cap L_n$ into two compact components,
 we consider the component
that contains $\gamma_m D_0 \cap L_n$.

We denote this component
of $\gamma_m V_0 \cap L_n \setminus (\gamma_m V_0 \cap D_0)$
 by $H_m$, and the other component by $C_m$. For
$C_m$ see Figure \ref{f.hm}.

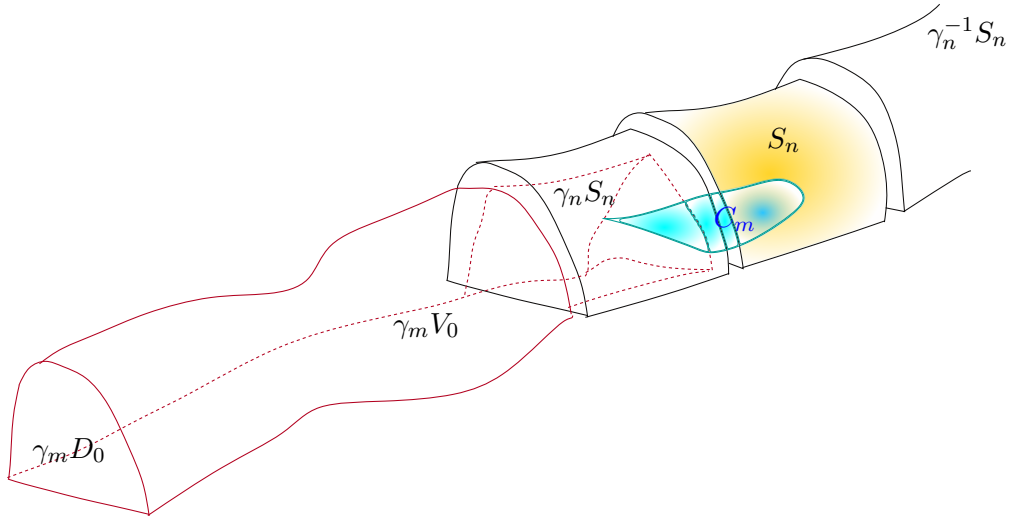
\begin{figure}[htbp]
\begin{center}
\includegraphics[scale=0.69]{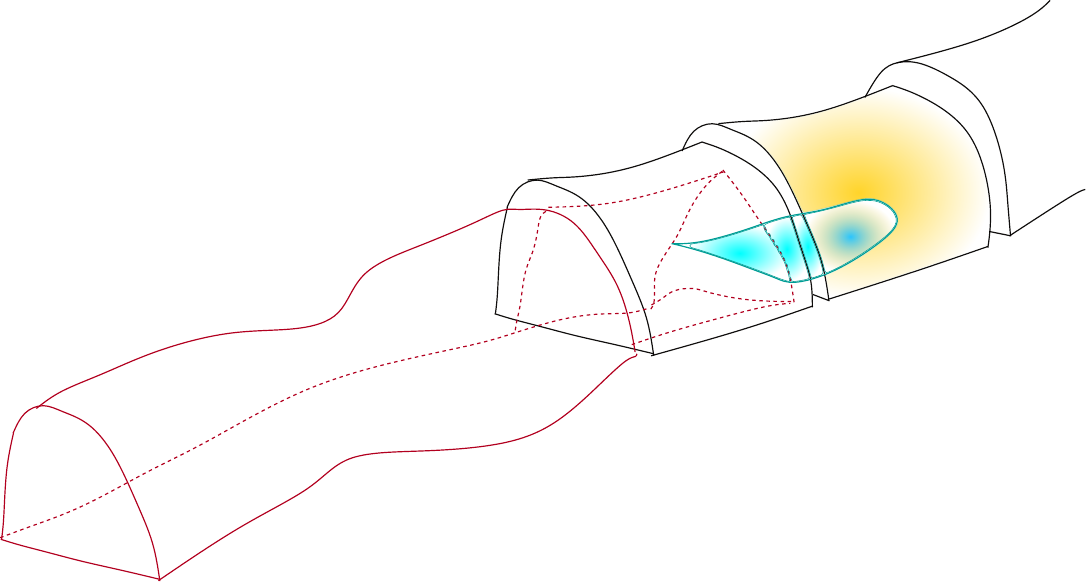}
\begin{picture}(0,0)
\put(-80,139){$S_n$}
\put(-20,179){$\gamma_n^{-1} S_n$}
\put(-160,119){$\gamma_nS_n$}
\put(-100,109){{\color{blue}$C_m$}}
\put(-220,69){$\gamma_m V_0$}
\put(-355,22){$\gamma_m D_0$}
\end{picture}
\end{center}
\vspace{-0.5cm}
\caption{{\small Returns by $\gamma_n$ and $\gamma_m$.}\label{f.hm}}
\end{figure}

First, as with $n$, we know that $\gamma_m V_0 \cap L_n$ is
a compact disk, and its boundary is the union of an arc
in $\gamma_m D_0$ and an arc in $\gamma_m V_0 \cap L_n$
whose interior is disjoint from $\gamma_m D_0$. The disk
$D_0$ cuts this disk so now the boundary of $H_m$ is made up as the
union of four compact arcs (see Figure \ref{f.hm2}): 
\begin{itemize}
\item the intersection $a_1$ of $\gamma_m D_0$
with $L_n$, \item two compact arcs $a_2$ and $a_4$ in $\gamma_m E \cap L_n$, and \item the arc $a_3=\gamma_m E \cap D_0 \cap L_n$ (which is one of the boundaries of $C_m$ defined above). \end{itemize}

We start in  the arc $a_3$ and move along both boundary arcs $a_2$ and $a_4$, that is, we are moving outside of $V_0$. We are interested in how the disk $H_m$ keeps intersecting
the serrated two dimensional set $B_n$. 
Notice that both sets are contained in the leaf $L_n$.
We also remark that the interior of the arc
$a_3$ is entirely contained in the interior of $B_n = \cS_n \cap L_n$ (see Proposition \ref{prop-serratedset}).

As we follow the two curves $a_2$ and $a_4$, notice that both are 
contained in the same leaf of $\wcG$, and they connect to each
other in the region disjoint from $\gamma_m S_m$, so the
curves in $a_2$ and $a_4$ never connect to each other. As long
as they are contained in the serrated set $B_n$,
they will obviously intersect this serrated set.

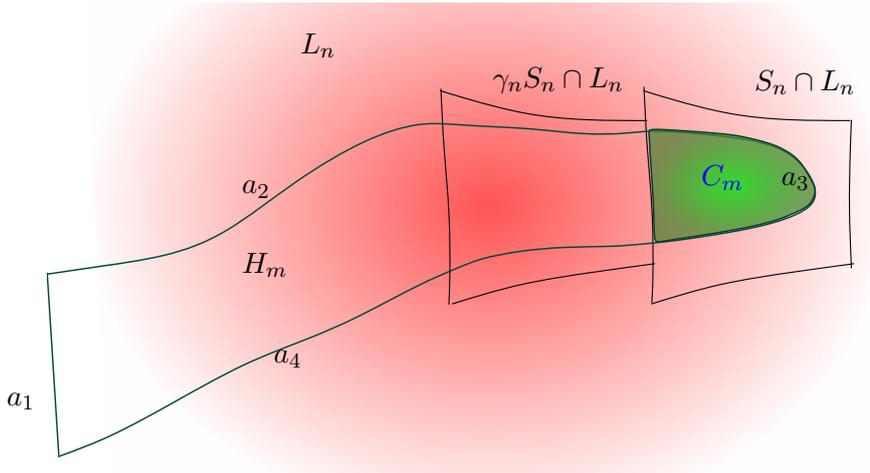
\begin{figure}[htbp]
\begin{center}
\includegraphics[scale=0.49]{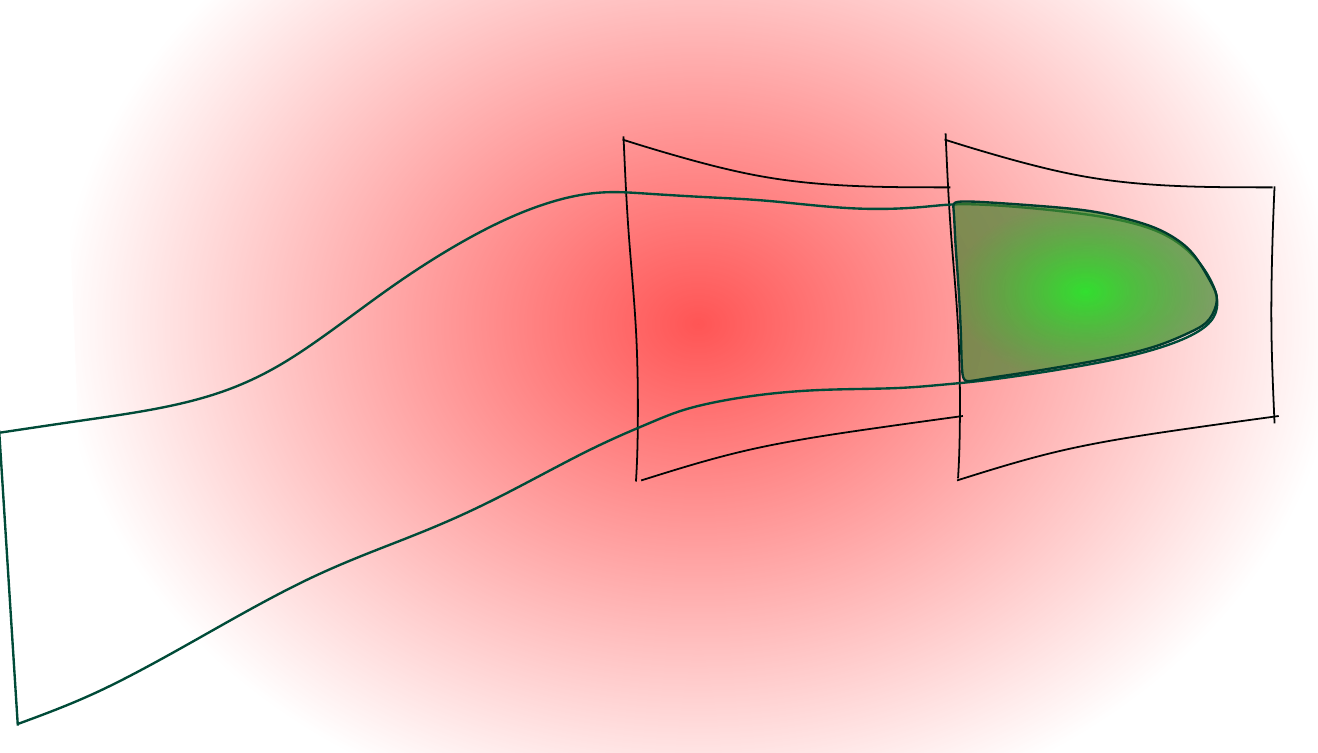}
\begin{picture}(0,0)
\put(-50,145){$S_n \cap L_n$}
\put(-220,159){$L_n$}
\put(-148,146){$\gamma_nS_n \cap L_n$}
\put(-70,109){{\color{blue}$C_m$}}
\put(-40,109){$a_3$}
\put(-330,26){$a_1$}
\put(-242,106){$a_2$}
\put(-242,76){$H_m$}
\put(-230,42){$a_4$}
\end{picture}
\end{center}
\vspace{-0.5cm}
\caption{{\small An image of the leaf $L_n$ and the boundaries of the disk $H_m$.}\label{f.hm2}}
\end{figure}

The first option is that at least one of the arcs $a_2$ and $a_4$ never escape $B_n$,  in which case 
$\gamma_m D_0$ clearly intersects $B_n$.
The second option is that eventually both the curves $a_2$ and $a_4$ escape 
from $B_n$.
They can only do this in the part of the boundary of $B_n$ 
which is transverse to $\wcG$ restricted to
$L_n$ (as the other parts are contained in a leaves 
of $\wcG$ and $a_2,a_4$ are pieces of leaves of $\wcG$). Since each of these arcs in $\partial B_n$ is 
transverse to $\wcG$ in $L_n$,
it follows that the two curves $a_2$ and $a_4$ must escape through different components
of $\partial B_n$, because $a_2, a_4$ are contained in the same leaf
of $\wcG$. Here notice that $B_n$ is an infinite
band and has two boundary components, each of which is a
bi-infinite, properly embedded curve.
This property of $a_2, a_4$ exiting $B_n$ through different
components of $\partial B_n$,  
implies that the endpoints belong to different connected components of $L_n \setminus B_n$ and therefore the disk $H_m$ still keeps intersecting
$B_n$ even after the curves $a_2$ and $a_4$ leave
the serrated set. 

The final conclusion is that the only way the endpoints
of the two curves $a_2$ and $a_4$ can be connected by the arc $a_1$ is
if the arc $a_1$ intersects $B_n$. This completes the proof of the
lemma. 
\end{proof}

We have just proved that for any $m > n$ then 
the disk $R_m := \pi_{n}(\gamma_m D_0)$ intersects the set
$A_n$ in $M_n$ which is a compact annulus (note that $M_n$ is not compact, so compactness of $A_n$ is important). 

This allows us to obtain: 

\begin{lema}\label{lema-cyclicsubgroup}
Up to a subsequence of $\{ \gamma_k, k > n \}$, one has that for all $m >n$ the elements $\gamma_m \circ \gamma_{n+1}^{-1}$ belong to the (cyclic) group generated by $\gamma_n$. 
\end{lema}
\begin{proof}
Since the diameters of $R_m$ are all equal to each other, we may assume up to subsequence of
the $\gamma_m$ with $m > n$
that for any $m > n$ then $R_m$ are all very close to each other.

But these are all compact sets, which are projections
of deck translates (in $\mt$) of the same fixed compact set.
It follows that up to subsequence we may assume that they
are all equal. For simplicity of notation we still denote this
subsequence by $\gamma_m, m > n$.
So we have that for any $m > n$, then

$$\pi_{n}(\gamma_{m} D_0) \ \ = \ \ \pi_{n}(\gamma_{n + 1} D_0).$$

Since $\pi_1(M_{n}) = < \gamma_n >$, 
this implies that there are integers $j_m$ with 
$\gamma_m = \gamma_n^{j_m} \circ \gamma_{n+1}$.

In other words $\gamma_m \circ \gamma_{n+1}^{-1}$ are all in the
subgroup generated by $\gamma_n$ as desired. 
\end{proof}

We can now complete the proof of Proposition \ref{prop-nomixedreturns}.

\begin{proof}[Proof of Proposition  \ref{prop-nomixedreturns}.]
Let $m>n$ in the subsequence of Lemma \ref{lema-cyclicsubgroup}. 
The proof analyzes $V_m$, and in particular the projection
of $V_m$ to the fixed $M_n$. 
Lemma \ref{lem-intersects} shows that $\pi_n(V_m)$ 
never stops intersecting 
$A_n$. Recall that $\pi_1(M_n) = <\gamma_n>$, Lemma
\ref{lem-intersects} shows that moving around in 
$\pi_n(V_m)$, has to keep intersecting $A_n$. Then in terms of homotopy, $\pi_n(V_m)$ has to be going around 
$\gamma_n$ iterates, either forwards or backwards.
We will look at the direction that the leaves of
$\pi_n(\wcF_1)$ and $\pi_n(\wcF_2)$
in the
boundary of $\pi_n(V_m)$ move as one goes around iterates of $\gamma_n$, and
obtain that they are nested, and this will lead to a 
contradiction with expanding holonomy along $r_1$. (See Figure \ref{f.returns}.)

\begin{figure}[htbp]
\begin{center}
\includegraphics[scale=0.76]{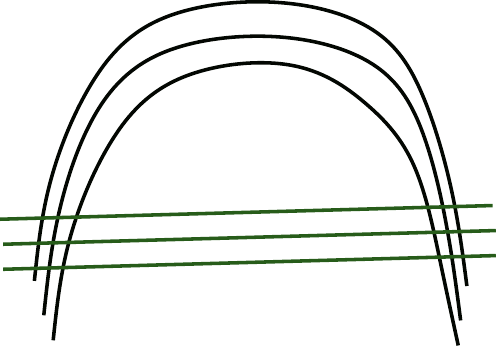}
\begin{picture}(0,0)
\put(-18,57){$p_0$}
\put(-52,22){$\gamma_m p_m$}
\put(-193,46){$L$}
\put(-83,128){$\partial D_0 \subset E \cup L$}
\put(-109,89){$\eta_m \subset \gamma_mE$}
\put(-203,36){$\gamma_n L$}
\put(-208,26){$\gamma_m L$}
\end{picture}
\end{center}
\vspace{-0.5cm}
\caption{{\small Since returns are in the same cyclic group, one should see expansion from $\gamma_m L$ to $\gamma_n L$ if $m \gg n$ but this is not seen.}\label{f.returns}}
\end{figure}
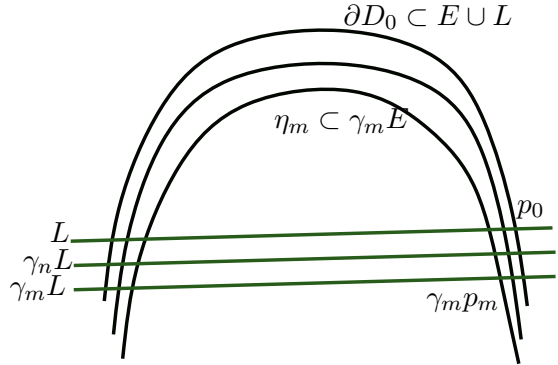

We abuse notation and denote the foliations
in $M_{n}$ by $\cF_1, \cF_2, \cG$ even if $M_n$ is a covering of $M$ (and quotient of $\mt$).

For any $m>n$, consider $\pi_n(\gamma_m D_m \cap D_0)$. This is a subdisk
of $\pi_n( D_0)$ with boundary the union of an arc $\eta_m$ in 
$\pi_n (\gamma_m E)$ and an arc in $\partial \pi_n( D_0)$.

The boundary of $\eta_m$ consist of two points which also bound
an arc $\eta'_m$ in $\pi_n(L \cap \gamma_m(E))$. Then $\eta_m \cup \eta'_m$
bounds a disk $W_1$ in its $\cF_2$ leaf. Now we move along
the boundary of $\pi_n(\gamma_m V_m)$ contained in $\pi_n(\gamma_m E)$,
starting in $\eta_m$ and
moving backwards (negative direction in $\pi_n(\gamma_m r_1)$).
Go once around following the direction of $\gamma_n$. 
The disk $W_1$ continues in this direction and returns
to $\pi_n(D_0)$ to a curve intersecting it 
between $\eta_m$ and the $\cF_2$ boundary of $\pi_n(D_0)$.

This means than going along $\gamma_n$ in the negative direction then 
the $\cF_2$ leaves in the boundary of
$\pi_n (\gamma_m V_m)$  return to the minus side transversely
of the original one (recall that the plus side points
into the disks $D_0$ and related disks).
Since the returns are on the negative side after one
$\gamma_n$ iteration, then they
are always on the negative side.
This means than when going forward (positive direction
along $\pi_n(r_1)$ then these returns are on the
positive $\cF_2$ side).

This takes care of the $\cF_2$ returns. Now consider the
$\cF_1$ returns. 
Here we take the opposite approach and move forward 
(with respect to $\pi_n(r_1)$) in the 
$\gamma_n$ direction as opposed to moving backwards as we did
in the previous analysis.
Since we are in the cases of no returns in
the same leaves they either return on the positive or
negative side (transverse to $\cF_1$). We want to show that
the forward returns are on the negative side
.
Start at $\pi_n(\gamma_m D_0)$ and move forward in the $\gamma_n$ direction.
The $\cF_2$ leaf returns on the positive side so intersects
the interior of the disk $\pi_n(\gamma_m D_0)$. Suppose first
that the return of
the $\cF_1$ leaf in the boundary is also on the positive
side. Then this
return it also intersects the interior of the disk 
$\pi_n(\gamma_m D_0)$. So this produces a compact 
subdisk $W_2$ in the interior
of $\pi_n(\gamma_m D_0)$. 
Since the return is inside this disk, recall that the intersection
of a leaf of $\wcF_1$ with the interior
of $V_0$ is a compact disk. So the
corner orbits intersect this in a compact set going forward
from the first return. But since $\pi_n(\gamma_m D_0)$ is the same
for all $m > n$, it follows that the sets $\pi_n(\gamma_m V_m)$
do more and more turns along the $\gamma_n$ direction.
This is impossible.

Therefore the $\gamma_n$ returns of the $\cF_1$ leaf are on the
negative side. Therefore the ray $\pi_n(\gamma_m r_1)$ 
in a leaf of $\cG$ always returns in the $(-,+)$ position
as we go around $\gamma_n$. As $m$ increases there are more and
more turns around $\gamma_n$ because $\gamma_m 
= \gamma_n^{j_m} \circ \gamma_{n+1}$ so in fact going around 
$\pi_n(\gamma_m V_m)$ is tracking $\gamma_n$ returns.

What we showed is that the returns associated
with the holonomy of positive direction of $\pi_n(r_1)$
are in the $(-,+)$ position with respect to the previous one.
The same happens for $\pi_n(\gamma_m p_m)$. Therefore if
one assumes they are all within $\ll \delta$ from the
limit $z$ of $\pi_n(\gamma_m p_m)$, the following happens:
consider a small transversal segment to $\cF_2$ from
the leaf containing $\pi_n(\gamma_{n+1} p_{n+1}))$ to the leaf containing
$\pi_n(\gamma_{n+2} p_{n+2}))$.
The arguments show that going forward along the corresponding
$\cG$ leaves, they are asymptotic to $z$.
Therefore the distance between the $\cF_2$ leaves
goes to zero. But this is in the positive transverse direction
to $\cF_2$ from doubly non separated rays. 
This contradicts the property of doubly non separated rays that have
to separate $> \eps_0$ from some time on going forward
for any other leaf in the same $\cF_1$ leaf (see Proposition \ref{prop-expandingholonomy}).
This contradiction shows that this situation cannot happen and finishes the proof of 
Proposition \ref{prop-nomixedreturns}.
\end{proof}

\subsection{Finding the Reeb surface}\label{ss.contradiction}

Using Propositions \ref{prop-noinside} and \ref{prop-nomixedreturns} we deduce that $r_1$ needs to project in a closed circle.  This also applies 
for $r_2$. Now we can apply Proposition \ref{prop-genRS} to deduce that the region between the rays project to Reeb surfaces in both the projection of $E$ and $L$. This completes the proof of Theorem \ref{teo-doublynsreeb} (which implies Theorem \ref{teo-noninjective}).

\section{Unique integrability, Proof of Theorem A}\label{s.uniquelyintegrable}
In this section we obtain uniqueness 
 of branching foliations under the assumption that these have Gromov hyperbolic leaves, which will allow us to prove Theorem \ref{teo.mainorient}.

\subsection{Unique integrability}

Note that in Theorem \ref{teo.mainph} we showed that if $f$ is a chain-transitive orientable partially hyperbolic diffeomorphism in a manifold with exponential growth of fundamental group, then, if $\cW^{cs}$ and $\cW^{cu}$ are branching foliations, then, after blow up, the flow foliation of a (topological) Anosov flow.  To generalize the setting, we will say that branching foliations $\cW^{cs}$ and $\cW^{cu}$ of a general orientable partially hyperbolic diffeomorphism intersect in \emph{good position} if their intersection is, after blow up, the orbit foliation of a topological Anosov flow.  The goal of this subsection is to show: 

\begin{teo}\label{teo.unique} Let $f: M \to M$ be an oriented partially hyperbolic diffeomorphism (cf. \S \ref{ss.branchingfol}) and let $\cW^{cs}_1$ and $\cW^{cs}_2$ be two $f$-invariant branching foliations tangent to $E^{cs}$ and $\cW^{cu}$ an $f$-invariant  branching foliation so that both $\cW^{cs}_1$ and $\cW^{cu}_1$ are in good position with respect to $\cW^{cu}$. Then, $\cW^{cs}_1 = \cW^{cs}_2$. 
\end{teo}

Let us consider $\cG_i$ to be the 1-dimensional branching foliations $\cW^{cs}_i \cap \cW^{cu}$ which, by assumption must be leafwise quasigeodesic
in both $\cW^{cs}_i$ and $\cW^{cu}$ (which must be Gromov hyperbolic).

The arguments that follow have some resemblance with some arguments in \cite[\S 10]{FP-dt} but add some new aspects that allow us to deal with the general case. 

We will use the following fact (see \cite[\S 5]{BFP}, \cite[\S 2.8]{FP-hsdff}) which is a consequence of $f$ being a strong collapsed Anosov flow with respect to both pairs $(\cW^{cs}_1,\cW^{cu})$ and $(\cW^{cs}_2,\cW^{cu})$: if $L \in \wcW^{cu}$ is a leaf and $S^1(L)$ denotes its circle at infinity, then, it follows that there is a unique point $p_i \in S^1(L)$ called the \emph{non-marker point} which verifies that every leaf of $\wcG_i$ is a quasigeodesic so that one of its endpoints is $p_i$. 
The non marker point depends only on the foliation $\wcW^{cu}$ and not
on the subfoliation $\wcG_i$, because it is the unique ideal point
of $L$ so that nearby leaves of $\wcW^{cu}$ are not asymptotic to
$L$ in that direction. It follows that $p_1 = p_2$, and we denote it
by $p$.
Moreover, for every $q \in S^1(L)$ there is a unique leaf of $\wcG_1$ (resp. $\wcG_2$) contained in $L$ and such that it limits in $p$ and $q$ (see Figure \ref{f.qgfan}). Since the curves of $\wcG_i$ are quasigeodesics and $f$ maps leaves of $\cG_i$ to leaves of $\cG_i$ it follows that if $\hat f$ is a lift of an iterate of $f$ to $\mt$, then it induces an action on $S^1(L)$ for every $L \in \wcW^{cu}$.

\begin{lema}\label{lema-notuniquethen}
Assuming that $\cW^{cs}_1 \neq \cW^{cs}_2$ then, there is a leaf $L \in \wcW^{cu}$ for which there are center curves $c_1 \in \wcG_1$ and $c_2 \in \wcG_2$ contained in $L$, with the same endpoints in the circle at infinity $S^1(L)$ of $L$ and such that there is a strong unstable leaf in between. 
\end{lema}

\begin{proof}
First we remark that if $\cW^{cs}_1 \not = \cW^{cs}_2$, then 
there is a leaf $L_0$ of $\wcW^{cu}$ for which there is a 
leaf of $\wcG_1$ in $L_0$ which is not a leaf of $\wcG_2$ in $L_0$.
The proof is exactly the same as in \cite[Proposition 10.7]{FP-dt}.

Let $\ell_1, \ell_2$ leaves of $\wcG_1, \wcG_2$ respectively
in $L_0$, so that $\ell_1 \neq \ell_2$ but have
the same ideal points in $S^1(L_0)$.
Suppose first that they intersect. Then there is 
$x$ in the intersection so that it is an endpoint
of compact segments in $\ell_1, \ell_2$ which intersect
only in $x$. Take a small unstable segment $\zeta$ near $x$,
but not through $x$ and intersecting both $\ell_1, \ell_2$,
which we shorten if necessary so that $\zeta$ intersects
them at the distinct endpoints. Let $t$ be an interior point of $\zeta$.
Project to $M$ and iterate by positive powers of $f$
and take a subsequence so that $f^{n_j}(t)$ converges in $M$.
Lifting to $\mt$ we get a sequence of lifts $g_j$ of iterates of $f$
 so that $g_j(t)$ converges to a point $y$. Leaves of $\wcG_i$ in the
same leaf of $\wcW^{cu}$ which have same ideal points
are at a globally  bounded Hausdorff distance from each other
in their leaf of $\wcW^{cu}$ because they are uniform
quasigeodesics. Since the leaves
$g_j(\ell_1), g_j(\ell_2)$ have the same ideal points, we
can find further subsequences of these two sequences that converge to a pair
of leaves $\tau_1, \tau_2$ of $\wcG_1, \wcG_2$
respectively in a leaf $L$ of $\wcW^{cu}$.
In addition they have the same pair of ideal points in $S^1(L)$.
The $g_j$ iterates of the segment $\zeta$ converge to at least the full 
leaf $u$ of $\wcW^u$ through $y$, which cannot intersect $\tau_1, 
\tau_2$. If $u$ separates $\tau_1$ from $\tau_2$ in $L$ 
these are the required leaves in the statement.
Otherwise both ends of $u$ converge to the same ideal
point of $S^1(L)$. In this case we zoom in to this ideal
point and apply deck transformations to bring back to a 
compact part of $\mt$, in the limit we obtain the desired
set of leaves.

On the other hand if $\ell_1, \ell_2$ do not intersect,
we do the following: if there is
 an unstable segment which intersects both of them, we
do the same procedure as above. Otherwise there is an unstable
leaf between them, and we apply the last part of the
proof.
\end{proof}

We now orient $E^c$ so that each curve of $\wcG_1$ points towards the 
non marker point in its $\wcW^{cu}$-leaf. We next prove
the following:

\begin{lema}\label{lem-opoori}
Each curve of $\wcG_2$, with the orientation of $E^c$, points away from
the non marker point in its $\wcW^{cu}$-leaves.
\end{lema}

\begin{proof}
We will use the objects given by the previous lemma.
Let $B$ be the region in $L$ which is
bounded by $c_1 \cup c_2$.
The set $B$ has only two accumulation points in $S^1(L)$: the non-marker point 
$p$, and another point which we denote by $q$.
Let $I, J$ be the connected components of $S^1(L) \setminus
\{ p, q \}$, where $I$ is contained in the accumulation set of the component
$L \setminus 
\{ c_1 \}$ which does not contain $c_2$, and same for $J$ (i.e. $I$ is ``closest'' to $c_1$ and $J$ ``closest'' to $c_2$).

Let $x \in c_i$. Since there is $u_0$ an unstable leaf inside $B$ joining $p$ and $q$, it separates $c_1$ from $c_2$, so the ray $r_x$ of  the unstable leaf $u_x$ through $x$ entering $B$ must converge to either $p$ or $q$. Assume without loss of generality that $x \in c_1$ and assume that $r_x$ converges to $q$, so, the ray of $c_1$ from $x$ to $q$ together with $r_x$ bound an open disk $D \subset B$. If $z \in D$ it follows that its $\wcG_1$ leaf $\ell$ has an ideal point in $J$ and an ideal point in $p$, which forces it to intersect $r_x$ twice, \footnote{It could also intersect $c_1$.} a contradiction. This contradiction shows that $r_x$ limits to $p$ for all $x \in c_1 \cup c_2$.

Let now $u_i$ be the unique unstable leaf contained in $B$,
separating $c_1$ from $c_2$ and so that $u_i$ is
closest to $c_i$. Consider $x$ in $c_1$ and $r_x$ as in the previous
paragraph and let $c$ a leaf of $\wcG_1$ intersecting the disk bounded by the ray of $c_1$ from $x$ to $p$ and $r_x$. We know that $c$ limits in $p$ and in some point in the interior of $J$, in particular, it must intersect and cross $c_2$ (it is no problem to cross $c_2$ as $c, c_2$ are not
in the same foliation!).  It follows that $c$ has to intersect both $u_1$ and
$u_2$ (which could coincide).  Let $z$ be the first point (starting
from the ray to $p$) of intersection between $c$ and $c_2$.

Let $u_z$ be the
unstable leaf of $z$. Let also $r_z$ be the ray in $u_z$ pointing inside $B$ which we proved above must limit to $p$. Note that $c$ can intersect $r_z$ only at $z$. By assumption the orientation of $E^c$ along $c_1$ points towards the non marker point $p$.
Since $c$ intersects $u_1, u_2$, and only intersects
$r_z$ in $z$ this forces the following property: $c$ intersects $c_2$ at $z$ as one goes away from $p$ in the direction of $c_2$ pointing towards $p$. In particular, the orientation of $E^c$ in $z$ in $c_2$ has to be against the direction of $p$.

Finally along leaves of $\wcG_2$,
seen inside leaves of $\wcW^{cu}$,  either the orientation
of $E^c$ always points towards the non marker point
in the $\wcW^{cu}$ leaf, or always points away from the
non marker in such a leaf. Since in $z \in L$ we showed
that the orientation of the leaf $c_2$ of $\wcG_2$ 
points away from the non marker point $p$, this proves
the lemma.
\end{proof}

Consider $\ell_1 \in \wcG_1$ which is fixed by some lift $\hat f$ of an iterate $f^k$ of $f$. Such a curve exists by partial hyperbolicity (fix a recurrent point and use that transverse to the center direction one has hyperbolic behaviour, so an index argument suffices). The leaf $\ell_1$ is in a leaf $L \in \wcW^{cu}$ (unrelated to the leaf obtained in Lemma \ref{lema-notuniquethen}). Associated to $\ell_1$ there is a unique leaf $\ell_2$ of $\wcG_2$ contained in $L$
with the same endpoints in $S^1(L)$ as $\ell_1$.
Since $\hat f$ preserves the foliations $\wcG_i$, it follows that
it acts on $S^1(L)$. In particular, by uniqueness of $\ell_2$ it follows
that $\hat f$ also leaves $\ell_2$ invariant. 

Because the orientation of $E^c$ points towards the non marker point
of $L$ on $\ell_1$ and away from the non marker point on  $\ell_2$ (due to Lemma \ref{lem-opoori}) we deduce that these two cannot coincide. In fact, we will produce a contradiction from the possible ways these curves intersect:   

\begin{lema}\label{lem-notdisj}
The leaves $\ell_1$ and $\ell_2$ cannot be disjoint. 
\end{lema}

\begin{proof}
Assume that $\ell_1 \cap \ell_2 =\emptyset$. Let $B$ be the region in $L$ between
$\ell_1$ and $\ell_2$, with ideal points $p$, the non marker
point in $S^1(L)$, and $q$. As in Lemma \ref{lem-opoori} let $I, J$ be the components
of $S^1(L) \setminus \{ p, q \}$ with $I$ ``closest'' to $\ell_1$, and $J$ ``closest'' to $\ell_2$.

Note first that unstable leaves cannot intersect both $\ell_1$ and $\ell_2$ because otherwise, the $E^c$ orientation of both would point towards the non-marker point (or against) contradicting Lemma \ref{lem-opoori}. 

Then, as in the proof of Lemma \ref{lem-opoori}, for any
unstable $u$ intersecting either $\ell_1$ or $\ell_2$, the
ray of $u$ entering $B$ has to limit to the non marker point $p$.
Then consider rays $u_{x_n}$ for $x_n$ a sequence in $\ell_1$ so that  $x_n \to q$. The limit contains a unique unstable leaf which separates $\ell_1$ from $\ell_2$ in $L$.

As in the proof of Lemma \ref{lem-opoori}, let 
$u_1, u_2$ the unstable leaves in $B$, with distinct ideal points $p, q$, fixed by $\hat f$, $u_1$ closest to $\ell_1$ and
$u_2$ closest to $\ell_2$.  

Consider $u_2$: since $\hat f$ fixes $u_2$, then $\hat f$ fixes a unique point
in $u_2$ call it $z_2$, and $\hat f$  is globally expanding 
in $u_2$.

Consider the set of leaves of $\wcG_1$ intersecting $u_2$ in the open ray of $u_2 \setminus \{ z_2 \}$  limiting to $q$. This is an open
interval of leaves of $\wcG_1$ in $L$, invariant by $\hat f$. 

Since these leaves intersect $u_2$, these leaves all have
one ideal point in $J$. Fix $c$ be the leaf of $\wcG_1$ through the point $z_2$ so that its limit point distinct from $p$ is closest to $q$ in $J$ and let $J' \subset J$ be the open interval between $q$ and $t$, the limit point of $c$, not equal to $p$. It follows that $\hat f$ acts without fixed points in $J'$ and for every $t' \in J'$ one has that $\hat f^n (t') \to q$ as $n \to + \infty$ and $\hat f^n(t') \to t$ as $n \to -\infty$. 

Notice that the rays of $c$ and $\ell_2$ with ideal points not equal to $p$
have to intersect in a compact set. The endpoints
of this intersection 
are fixed by $\hat f$. Let $w$ be the endpoint ``closest" to $t$ in $c$.
Let $r_w$ be ray of the unstable of $w$ in the outside of $B$.  Then $r_w$ is fixed by $\hat f$ and has ideal point in the closure of $J'$ so it must be either $t$ or $q$. 

We consider each case separately. Suppose first that $r_w$ limits on $t$. Then consider leaves of $\wcG_1$ intersecting
$r_w$ near $w$. They will have an ideal point in the interior of $J'$ and hence will have to intersect $r_w$  again, which is
impossible. 

The other option is that $r_w$ has ideal point $q$. Now, instead of $\wcG_1$ we consider leaves of $\wcG_2$ which intersect $r_w$ near $w$ which must have an ideal point in $J'$ and hence will have to intersect $r_w$
a second time, again a contradiction.
This finishes the proof of the lemma.
\end{proof}

Finally, we elliminate the other case and complete the proof of Theorem \ref{teo.unique}: 

\begin{lema}\label{lem-notintersectl1l2}
The leaves $\ell_1$ and $\ell_2$ cannot intersect. 
\end{lema}

\begin{proof}
Note that Lemma \ref{lem-opoori} states that the positive directions of $E^c$ on $\ell_1$ point towards $p$ and in $\ell_2$ point away from $p$ where $p$ is the non-marker point in $L$. Denote by $q$ the other limit point of both $\ell_1$ and $\ell_2$. Because of the orientation, we know that $\ell_1$ and $\ell_2$ cannot share a ray. Assume first that there is a sequence $z_n \in \ell_1 \cap \ell_2$ escaping to infinity, and assume that it is ordered so that $z_n$ is monotonic in $\ell_1$, in particular, it either converges to $p$ or to $q$ and thus, up to subsequence, we can also assume that $z_n$ is monotonic in $\ell_2$ (but if it is increasing in $\ell_1$ it is decreasing in $\ell_2$). Now, it follows from the orientation that the arc of $\ell_1$ from $z_n$ to $z_{n+1}$ and the arc of $\ell_2$ from $z_{n+1}$ to $z_n$ is a smooth closed curve whose self-intersections are tangent, therefore it contains a smooth simple closed curve. This simple closed curve is transverse to $E^u$ and bounds a disk, which contradicts the fact that $E^u$ has no singularities.

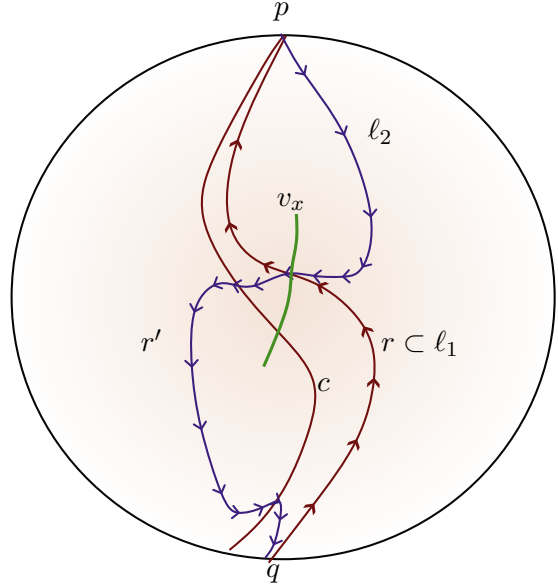
\begin{figure}[htbp]
\begin{center}
\includegraphics[scale=0.68]{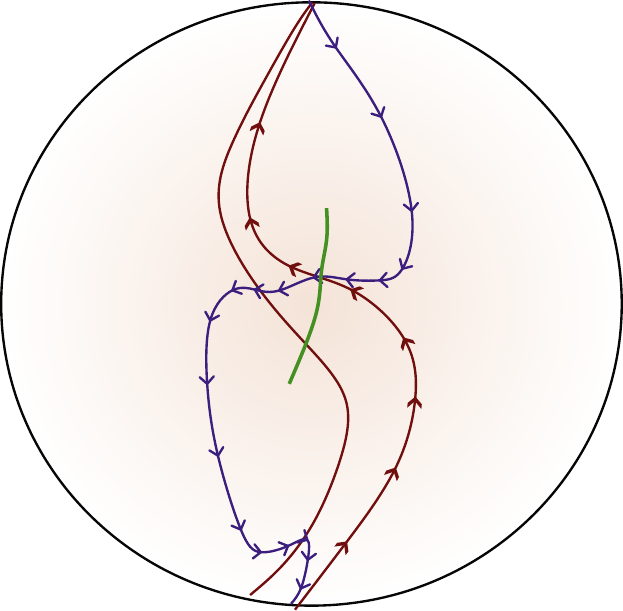}
\begin{picture}(0,0)
\put(-109,207){$p$}
\put(-112,-5){$q$}
\put(-74,160){$\ell_2$}
\put(-69,80){$r \subset \ell_1$}
\put(-159,80){$r'$}
\put(-108,135){$v_x$}
\put(-93,65){$c$}
\end{picture}
\end{center}
\vspace{-0.5cm}
\caption{{\small A depiction of the final part of the argument of Lemma \ref{lem-notintersectl1l2} .}\label{f.orient}}
\end{figure}

We can then assume that $\ell_1 \cap \ell_2$ is compact. Let $r$ be the ray of $\ell_1 \setminus \ell_2$ which limits to $q$. The ray $r$ has starting point $x$ which also belongs to $\ell_2$. Let $r'$ be the ray in $\ell_2$ from $x$ to $q$. 
The union $r \cup r'$ bounds an open disk $D$ in $L$ which
limits only on $q$ in $S^1(L)$. 
In addition again by Lemma \ref{lem-opoori} it follows that the
boundary of $D$ is a smooth curve which implies that $D$ is cut into two by the ray $v_x$ of the unstable of $x$ entering $D$. 
Notice that $v_x$ limits on $q$ as it is contained in $D$.

Let $Y$ be the component of $L \setminus \ell_1$ disjoint
from $D$, and $I$ the component of $S^1(L) \setminus \{ p, q \}$
which is contained in the limit of $Y$.  Let $z$ be a point in $r'$
and let $c$ be its $\wcG_1$ leaf. One ray of $c$ limits to $p$, and the
other limits to a point in $I$, and this forces $c$ to intersect
$v_x$ twice, which is impossible and completes the proof. See Figure \ref{f.orient}.
\end{proof}

\subsection{Removing the orientability assumptions} 
The following result is what allows us to bring back
branching foliations associated with a lift of a finite
iterate of $f$, to 
branching foliations for $f$ in $M$ and is thus key to the proof of Theorem \ref{teo.mainorient}.

\begin{prop} \label{teo-mainorient1}
Let $f: M \to M$ be a partially hyperbolic diffeomorphism, and suppose
that there is $\hat f: \hat M \to \hat M$ a lift of an iterate to $f$ to a finite cover $\hat M$ such that $\hat f$ is an oriented partially hyperbolic diffeomorphism, and such that every pair of branching foliations $\hat f$ preserves  intersect in good position. 
Then there are unique branching foliations in $M$ tangent
to $E^{cs}$ and $E^{cu}$ for $f$, and every (complete)
curve tangent to $E^c$ in $\mt$ can be obtained
as the intersection of a leaf $L \in \wcW^{cs}$ and
a leaf $E \in \wcW^{cu}$.
These branching foliations 
in $M$ are $f$-invariant, and intersect in 
a leafwise quasigeodesic manner.
\end{prop}

\begin{proof} 
We assume without loss of generality that the finite cover $\hat M \to M$
is normal.

Theorem \ref{teo-BI1} of Burago-Ivanov shows that $\hat f$ 
preserves a pair of branching foliations.
The hypothesis of this theorem and 
Theorem \ref{teo.unique} 
imply that there is a unique pair of $\hat f$ invariant
branching foliations,
which we denote by 
 $\hat \cW^{cs}$ and $\hat \cW^{cu}$.

We now prove that 
these facts imply that $\hat \cW^{cs}$ and $\hat \cW^{cu}$ are invariant under the deck transformations of the covering $\hat M \to M$. 
To prove this statement we review the results in \cite{BI}: they proved that
if all the bundles are orientable then there are branching foliations
$\hat \cW^{cs}, \hat \cW^{cu}$
tangent to $E^{cs}$ and $E^{cu}$. We discuss $\hat \cW^{cs}$: their
construction in fact produces two possible foliations, the
uppermost and lowermost. The uppermost is obtained by patching
local surfaces tangent to $E^{cs}$ which are ``uppermost" in
the positive center direction (this uses that the center direction
is orientable and so is the unstable direction which is transverse
to $E^{cs}$). In particular the construction of the 
uppermost branching foliation is totally independent of the 
partially hyperbolic map preserving the bundles. Then they
prove that any diffeomorphism which preserves the bundles and the
orientations, then it preserves the uppermost branching
foliation (and also preserves the lowermost branching foliation).
If $\gamma$ is a deck of the cover $\hat M \to M$,
then it preserves the bundles. If it preserves the orientations
then it will preserve the uppermost branching foliation. If
it reverses orientation, then it sends local uppermost surfaces
to local lowermost surfaces. So it takes the uppermost branching
foliation to the lowermost branching foliation.

In addition $\hat f$ preserves the uppermost and lowermost branching
foliations. 
Theorem \ref{teo.unique} applied to the uppermost and lowermost $cs$-branching foliation and any (say the uppermost) $cu$-branching foliation shows that  the uppermost and lowermost $cs$-branching foliations 
must coincide. 
In particular for any such $\gamma$, then $\gamma \hat \cW^{cs} = \hat \cW^{cs}$.
It follows that $\hat \cW^{cs}$ projects down to a branching foliation $\cW^{cs}$
in $M$. The same happens for $\hat \cW^{cu}$.

In addition $\hat \cW^{cs}, \hat \cW^{cu}$ are the only branching foliations
tangent to the $E^{cs}, E^{cu}$ bundles in $\hat M$. This is
because any branching foliation tangent to $E^{cs}$ has to be between the
lowermost and the uppermost branching foliations (see \cite[\S A.2]{BFP}) which we showed must coincide. Likewise for $E^{cu}$.

Notice that the projected branching foliation must be the only branching foliation tangent to $E^{cu}$ in
$M$: if there were two, they would lift to two distinct branching
foliations in $\hat M$. 

Uniqueness then implies that the branching foliations 
$\cW^{cs}, \cW^{cu}$ are $f$-invariant: 
$f(\cW^{cs})$ is a branching foliation tangent to $E^{cs}$ 
so uniqueness implies that $f(\cW^{cs}) = \cW^{cs}$ and
similarly for $\cW^{cu}$.

Given that there are unique branching foliations for $f$
the fact that complete curves tangent to $E^c$ are obtained as the intersection of a leaf $L \in \wcW^{cs}$ and a leaf $E \in \wcW^{cu}$ follows exactly as in \cite[Proposition 10.6]{BFP}.

The property of the leafwise quasigeodesic intersection is direct since it
is true for $\hat f$ and  this condition is checked in the universal cover. 
\end{proof}

We are now ready to prove Theorem \ref{teo.mainorient}. Note that in \cite{BFeM} a construction of Anosov flows out of abstract actions on a bifoliated plane is developed that could be useful here, but we will give an elementary proof using properties of strong collapsed Anosov flows. 

\begin{proof}[Proof of Theorem \ref{teo.mainorient}.]
The previous theorem shows that $f$ preserves branching foliations
$\cW^{cs}, \cW^{cu}$ and that these are the unique branching foliations
preserved by $f$. We need to produce a (topological) Anosov flow in $M$ and show that it is associated
with $f$ to produce the strong collapsed Anosov flow property for $f$. 

Consider a finite cover
$\hat M$ of $M$ and a lift $g: \hat M \to \hat M$ of an iterate $f^k$ of $f$. Without loss of generality we can assume that $\hat M$ is a normal cover and we let $\Gamma$ the finite group of deck transformations of $\hat M$ related to the cover $\hat M \to M$. 

By Theorem \ref{teo.mainph} we know that $g$ is a strong collapsed Anosov flow with respect to the unique pair of branching foliations $\hat \cW^{cs}$ and $\hat \cW^{cu}$ which are the lifts of $\cW^{cs}$ and $\cW^{cu}$ respectively (see Proposition \ref{teo-mainorient1}). We denote by $\hat \cW^{c}$ to the center branching foliation for $g$. 

Let $\hat \varphi_t: \hat M \to \hat M$ a (topological) Anosov flow, and $\hat h: \hat M \to \hat M$ a continuous surjective map sending oriented orbits of $\hat \varphi_t$ into center curves of $g$, so that $g$ is a strong
collapsed Anosov flow with respect to $\hat \varphi_t$: there
is a self orbit equivalence $\hat \beta$ of $\hat \varphi_t$ which
satisfies $\hat h \circ \hat \beta = g \circ \hat h$.

Given a deck transformation $\gamma \in \Gamma$ we will construct a homeomorphism $\overline \gamma: \hat M \to \hat M$ in a way that the action of $\overline \Gamma = \{ \overline \gamma \ : \ \gamma \in \Gamma\}$ is a group isomorphic to $\Gamma$ acting by deck transformations in $\hat M$ and preserving the orbits of $\hat \varphi_t$. This way, we will be able to construct a (topological) Anosov flow $\varphi_t$ in the manifold $\overline{M} = \hat M/_{\overline \Gamma}$ which is diffeomorphic to $M$.

Given $\gamma \in \Gamma$ we define $\overline \gamma$ as follows: 

$$ \overline \gamma (x) = \hat h^{-1} \gamma (\hat h(x)). $$

Note that since $\hat h$ is not invertible, one needs to precise what $\hat h^{-1}$ means. However given $x \in \hat M$, it belongs to a single orbit of $\hat \varphi_t$. Then $\gamma (\hat h(x))$ belongs to the unique curve $\gamma c_x$ of $\hat \cW^{c}$ which is the image under $\gamma$ of $\hat h$ of the $\hat \varphi_t$-orbit of $x$. 
Here $c_x$ is a center curve.
We claim that there is a unique orbit of $\hat \varphi_t$ which
maps by $\hat h$ to $\gamma c_x$. To see this for instance do the
following: $x$ is in a unique 2-dimensional stable leaf $L$
of $\hat \varphi_t$. This maps by $\hat h$ to a center stable
leaf of $\hat{\cW}^{cs}$ and by $\gamma$ to another such
leaf. Then there is a unique 2-dimensional stable leaf $F$
of $\hat \varphi_t$ which is mapped by $\hat h$ to
$\gamma \hat h(L)$ \cite{BFP}. In $E$ there is a unique
flow line $\zeta$ which is mapped by $\hat h$ to $\gamma c_x$.
Now, since $\hat h$ lifts to a homeomorphism from the lift of the orbit of $\gamma x$ to the lift of $\gamma c_x$ we get that one has a well defined notion of $\hat h^{-1}$ associated to $x$ in $\gamma (\hat h(x))$. This is the choice we do.

We will check continuity of $\overline \gamma$, and that the bar operation
preserves the product, which since it is clear that $\overline{\mathrm{id}} = \mathrm{id}$, also shows that the inverse of $\overline \gamma$ exists and is continuous. 

Continuity follows from the fact that the choices are unique: if $x_n \to x$ in $\hat M$ it will follow that $\gamma (\hat h(x_n)) \to \gamma (\hat h(x))$ but also, the curves $c_{x_n} \to c_x$ and so it also holds that $\gamma c_{x_n} \to \gamma c_x$. The application of the inverse also makes the sequence converge showing continuity. 

One can also check $\overline{\gamma \circ \eta} = \overline{\gamma} \circ \overline{\eta}$ because 

$$ \overline{\gamma} \circ \overline{\eta} (x)=  \hat h^{-1} \circ \gamma \circ \hat h \circ \hat h^{-1} (\eta (\hat h(x)), $$ 

\noindent and the middle composition of $\hat h$ and $\hat h^{-1}$ are applied in the same center and same orbit, so they cancel out. 

This shows that $\overline \Gamma$ is a group action on $\hat M$.
Suppose that for some $\gamma$ deck transformation $\hat M \to M$
there is $x$ so that $\overline{\gamma}(x) = x$.
By definition $\overline{\gamma}(x)$ is a point in $\hat M$
which is mapped by $\hat h$ to $\gamma \hat h(x)$ (it may not
be the only one). So if $\overline{\gamma}(x) =x$ this
implies that  $\hat h(x) = \gamma \hat h(x)$.
Since $\gamma$ is a deck transformation this implies that $\gamma$
is the identity. Hence the action of $\overline \Gamma$
on $\hat M$ is free. Since it is finite it is obviously
properly discontinuous. Therefore the quotient
$\overline{M}$ is a manifold. It is also a $K(\pi,1)$
manifold with fundamental group naturally isomorphic
to $\pi_1(M)$, it follows that $\overline{M}$ and $M$ are
homotopy equivalent and hence, being 3-dimensional, diffeomorphic\footnote{For irreducible 3-manifolds, the fundamental group determines the topology \cite[\S 1.4.2]{BBBMP}, moreover, 3-manifolds admit a unique differentiable structure \cite{Moise}.}. We fix a differentiable
structure in $\overline{M}$ diffeomorphic to $M$.

We want to obtain a topological Anosov flow induced
in $\overline{M}$.
Any deck transformation $\gamma$ from $\hat M$ to $M$
preserves flow lines of $\hat \varphi_t$ by construction.
Hence the flow of $\hat \varphi_t$ descends to a
one dimensional  foliation $\cH$ in $\overline M$ (note that parametrizations may not behave well).
 
In addition any $\gamma$ deck transformation of $\hat M$ to $M$ preserves 
the stable and unstable foliations of $\hat \varphi_t$, so 
they descend to foliations in $\overline M$, which are
$\cH$ saturated.
In particular the foliation $\cH$ is expansive.
Therefore Proposition \ref{lem-weaktaf} implies that $\cH$ is
topologically equivalent to the flow foliation of
a topological Anosov flow $\phi_t$.

By construction the lift of the map $\hat h$ to the universal cover induces an equivariant map between the orbit space of $\varphi_t$ (which is equivariantly
the same as the
leaf space of $\cH$, and the same as the orbit space of
$\phi_t$) and the center leaf space of $f$ (which is the same as the center leaf space of $g$ in the universal cover). 
Since $f$ preserves branching foliations and $\phi_t$ is a topological Anosov flow, this means that
$f$ is a leaf space collapsed Anosov flow with respect to $\phi_t$.
Therefore one can now apply the exact same arguments as in \cite[\S 9]{BFP} to get the strong collapsed Anosov property for $f$ with
respect to the topological Anosov flow $\phi_t$, that is,
to obtain the strong collapsed Anosov flow property
from the leaf space collapsed Anosov flow property.
We note that the transverse orientability is used in \cite[\S 9]{BFP} only to produce the Anosov flow in $M$, so the arguments of \cite[\S 9]{BFP} 
apply here. This shows that $f$ is a strong collapsed Anosov flow
and finishes the proof of the theorem. 
\end{proof}

\end{document}